\documentstyle[11pt]{article}

\begin{document}
\noindent
\begin{center}
  {\LARGE Orbifold Gromov-Witten Theory}\footnote{both authors partially
supported by the National Science Foundation}
  \end{center}

  \noindent
  \begin{center}

     {\large Weimin Chen\footnote{New Address (starting Fall 2000):
Department of Mathematics, SUNY at Stony Brook, NY 11794. \\
Email: wechen@math.sunysb.edu}
             and Yongbin Ruan}\\[5pt]
      Department of Mathematics, University of Wisconsin-Madison\\
        Madison, WI 53706\\[5pt]
        Email: wechen@math.wisc.edu  and ruan@math.wisc.edu\\[5pt]
              \end{center}

              \def \x{{\bf x}}
              \def \J{{\cal J}}
              \def \M{{\cal M}}
              \def \A{{\cal A}}
              \def \B{{\cal B}}
              \def \C{{\bf C}}
              \def \Z{{\bf Z}}
              \def \R{{\bf R}}
              \def \P{{\bf P}}
              \def \I{{\bf I}}
              \def \N{{\cal N}}
              \def \T{{\cal T}}
              \def \O{{\cal O}}
              \def \Q{{\bf Q}}
              \def \D{{\cal D}}
              \def \H{{\cal H}}
              \def \S{{\cal S}}
              \def \e{{\bf E}}
              \def \CP{{\bf CP}}
              \def \U{{\cal U}}
              \def \E{{\cal E}}
              \def \F{{\cal F}}
              \def \L{{\cal L}}
              \def \K{{\cal K}}
              \def \G{{\cal G}}
              \def \z{{\bf z}}
              \def \m{{\bf m}}
              \def \n{{\bf n}}
              \def \V{{\cal V}}
              \def \W{{\cal W}}

\tableofcontents

\section{Introduction}

In 1985, Dixon, Harvey, Vafa and Witten considered  string
    theory over orbifolds (arising as global
quotients $X/G$ by a finite group $G$)\cite{DHVW}. Although an
orbifold is a singular space, orbifold string theory is
surprisingly a "smooth"  string theory. Since then, orbifold
string theory has become a rather important part of the landscape
of string theory. Although orbifold string theory has been around
for a while, it was apparently  poorly explored in mathematics.
For the last fifteen years, only a small piece of orbifold string
theory concerning the orbifold Euler number has been studied in
mathematics. This paper is one of our several efforts \cite{CR1},
\cite{Ru2}, \cite{AR} to change the situation.

    Even with a superficial understanding of orbifold string
    theory, it is obvious that the mathematics surrounding
    orbifold string theory is striking.
    In fact, it brings in much new mathematics unique to
    orbifolds. We believe that there is an emerging  new topology and geometry of orbifolds.
    The core of this new geometry and topology is the concept of
    twisted sectors. Roughly speaking, the consistency of orbifold string theory requires that
     the string
Hilbert space has to contain factors called twisted sectors.
Twisted sectors can be viewed as the contribution from
singularities. All other quantities such as correlation functions
have to contain the contribution from the twisted sectors. In
\cite{CR1}, we studied  twisted sectors in the context of
classical topology and obtained a new cohomology theory for
orbifolds (orbifold cohomology). In this paper, we continue our
work in quantum theory to construct an orbifold quantum
cohomology. Recall that ordinary quantum cohomology is a
deformation of ordinary cohomology. Orbifold quantum cohomology
can be thought as a deformation of orbifold cohomology.

    Another motivation of this paper comes from mirror symmetry.
    The current mirror symmetry is restricted to Calabi-Yau
    3-folds only. The most of known Calabi-Yau 3-folds are so
    called crepant resolution of Calabi-Yau orbifolds. In higher
    dimension, there are still plenty of Calabi-Yau orbifolds. But
    they do not have crepant resolution in general. Therefore,
    there is no hope to consider mirror symmetry for smooth higher
    dimensional Calabi-Yau manifolds. We are forced to work with
    orbifolds! This paper can be viewed as the first step towards
    higher dimensional mirror symmetry.

    The first step of our paper is to generalize the notion of stable map.
    This is a nontrivial step. Recall that orbifold is covered by orbifold charts of
    the form $U/G$, where $U$ is a smooth manifold and $G$ is a
    finite group acting on $U$. Let $G_p$ be the subgroup of  $G$ fixing a point of
    the preimage of $p$.
    $G_p$ is well-defined up to conjugation and is called the local group of $p$.
    The natural generalization of a smooth map in the orbifold category
    is the orbifold map $f: X\rightarrow Y$, where locally $f_p:
    U_p/G_p \rightarrow V_{f(p)}/G_{f(p)}$ can be lifted to a map
    $\tilde{f}: U_{p}\rightarrow V_{f(p)}$ and an injective
    homomorphism $f_{\#}: G_p\rightarrow G_{f(p)}$. Of course, the
    pair $(\tilde{f}, f_{\#})$ is only defined up to  conjugation
    of an element of $G_{f(p)}$. Suppose $X$ is a marked orbifold Riemann
    surface with orbifold point $(x_1, \cdots, x_k)$. For
    simplicity, we assume that all orbifold points are marked
    points. The local group of $x_i$ is determined by an integer
    $k_i$ representing the order of the local group. Then an orbifold
    map $f: X\rightarrow Y$ is completely determined by the map itself
    and the conjugacy classes $\x=((g_1), \cdots, (g_k))$ (We call it the
    orbifold type or twisted boundary
    condition of $f$) for $g_i\in
    G_{f(x_i)}$. The
    problem is that such a straightforward generalization is
    wrong for the purposes of quantum cohomology, where we have to
    consider the deformation theory of such a map. Suppose that $E\rightarrow Y$
    is an orbifold bundle of $Y$. The critical
    issue is that the pull-back $f^*E$ is not an
    orbifold bundle in general. If $f^*E$ has an orbifold bundle structure, it often has
    more than one orbifold bundle structure. The critical step of this paper is to
    formulate the correct structure ({\em compatible system}) such that the usual
    deformation theory is possible.
    We do this in the appendix and call an orbifold map with such a structure {\em a good map}.
    This extra structure has a far-reaching influence on the structure of orbifold cohomology.
    Once we get over this conceptual hurdle, we define an orbifold stable map as
     nodal good map
    with the usual stability condition. Let $\overline{\M}_{g,k}(X,J, A, \x)$ be the moduli space of
    orbifold stable maps of type $\x$. Then we prove
    \vskip 0.1in
    \noindent
    {\bf Theorem A (Theorem 2.3.8): }{\it Suppose that $X$ is  either
     a symplectic orbifold with tamed almost complex structure
    or a projective orbifold. The moduli space of
orbifold stable maps $\overline{\M}_{g,k}(X,J,A, \x)$ is a compact
metrizable space under a natural topology, whose ``virtual
dimension'' is $2d$, where $$ d=c_1(TX)\cdot
A+(\dim_{\C}X-3)(1-g)+k-\iota(\x).\leqno $$ Here
$\iota(\x):=\sum_{i=1}^k \iota_{(g_i)}$ for
$\x=(X_{(g_1)},\cdots,X_{(g_k)})$ and the degree shifting number
is defined in \cite{CR1}. }
    \vskip 0.1in

    For any component $\x=(X_{(g_1)},\cdots,X_{(g_k)})$, there are $k$
evaluation maps (cf. (3.5)) $$
e_i:\overline{\M}_{g,k}(X,J,A,\x)\rightarrow X_{(g_i)},
\hspace{4mm} i=1,\cdots, k. \leqno(1.1) $$ There is a map $$p:
:\overline{\M}_{g,k}(X,J,A,\x)\rightarrow \overline{\M}_{g,k}
\leqno(1.2)$$ defined by contracting the unstable components of
the domain.  For any set of cohomology classes $\alpha_i\in
H^{*-2\iota_{(g_i)}}(X_{(g_i)};\Q)\subset H^*_{orb}(X;\Q)$,
$i=1,\cdots,k$, $K\in H^*(\overline{\M}_{g,k}, \Q)$, the orbifold
Gromov-Witten invariant is defined as the pairing $$
\Psi^{X,J}_{(g,k,A,\x)}(\O_{l_1}(\alpha_1), \cdots,
\O_{l_k}(\alpha_k))=\int^{vir}_{
\overline{\M}_{g,k}(X,J,A,\x)}\prod_{i=1}^k c_1(L_i)^{l_i}e^*_i
\alpha_i,\leqno(1.3) $$ where $L_i$ is the line bundle generated
by cotangent space of the $i$-th marked point.

\vskip 0.1in

\noindent{\bf Theorem B (Proposition 3.4.1): }{\it
    \begin{enumerate}
    \item $ \Psi^{X,J}_{(g,k,A,\x)}(K; \O_{l_1}(\alpha_1), \cdots,
\O_{l_k}(\alpha_k))=0$ unless $deg K+ \sum_i
(deg_{orb}(\alpha_i)+l_i)=2C_1(A)+2n(3-g)+2k.$, where $deg_{orb}(\alpha_i)$
is orbifold degree of $\alpha_i$ obtained after degree shifting.
    \item $\Psi^{X,J}_{(g,k,A,\x)}(K; \O_{l_1}(\alpha_1), \cdots,
\O_{l_k}(\alpha_k))$ is independent of the choice of $J$ and hence
is an invariant of the symplectic structure.
\end{enumerate}}
    \vskip 0.1in
    \noindent
    {\bf Theorem C (Theorem 3.4.2): }{\it $\Psi^{X,J}_{(g,k,A,\x)}$
    satisfy the same axioms as the ordinary GW-invariants
except divisor axiom where we have to restrict ourself to divisor
class in the nontwisted sector.}
    \vskip 0.1in
    Similar to the smooth case, the genus zero orbifold GW-invariants can be
used to define an orbifold quantum cohomology which deforms the
orbifold cohomology constructed in \cite{CR1}.
    \vskip 0.1in
    \noindent
    {\bf Theorem D (Theorem 3.4.3): }{\it The orbifold quantum product is associative.}
\vskip 0.1in

    The moduli space of orbifold stable maps is constructed in section 2.
    In section 3, we carry out the construction of
    virtual cycles to define orbifold Gromov-Witten
    invariants. In this section, we adapted the technique from the
    smooth case developed by \cite{FO}, \cite{LT}, \cite{Ru1},
    \cite{S}. In particular, we follow the proof of \cite{FO} closely.
    However, the analysis is more complicated because of the existence of
    singularities.

    The results of this paper was announced in \cite{CR2}. The
    second author would like to thank R. Dijkgraaf, E. Witten and
    E. Zaslow for many stimulating discussions about orbifold
    string theory. Orbifold stable maps have been studied independently in
    algebraic setting in \cite{AV}.
     Both authors would like to thank D. Abramovich
    for clarifying the relation between our and their versions of
    orbifold stable maps.

\section{Orbifold Stable Maps}

\subsection{Almost complex orbifolds and pseudo-holomorphic maps}
    We follows the notations from the appendix. We refer readers to
    the appendix for the definitions.

    \vskip 0.1in
\noindent{\bf Definition 2.1.1: }{\it An almost complex structure
on an orbifold $X$ is a $C^\infty$ section $J$ of the orbifold
bundle $End(TX)$ of endomorphisms of $TX$ such that $J^2=-Id$. The
pair $(X,J)$ is called an {\it almost complex orbifold}. A
continuous map $f:(X,J)\rightarrow (X^\prime,J^\prime)$ is said to
be {\it pseudo-holomorphic} if there is a $C^\infty$ map
$\tilde{f}$ lifting $f$ such that $d\tilde{f}\circ J
=J^\prime\circ d\tilde{f}$.}
    \vskip 0.1in

The next lemma is obvious, we leave the details to the reader.
    \vskip 0.1in
\noindent{\bf Lemma 2.1.2: }{\it Let $(X,J)$ be an almost complex
orbifold.  Then $J$ induces an almost complex structure on
$\widetilde{\Sigma X}$ for which the canonical resolution map
$\pi: \widetilde{\Sigma X}\rightarrow X$ is pseudo-holomorphic. In
particular, the set of singularities $\Sigma X$ is a
pseudo-holomorphic subvariety in $(X,J)$.} \vskip 0.1in

When the domain is $2$-dimensional, it is convenient to use the following
slightly different version of pseudo-holomorphic map.

\vskip 0.1in \noindent{\bf Definition 2.1.3: }{\it Let $\Sigma$ be
a Riemann surface with complex structure $j$, and $(X,J)$ be an
almost complex orbifold. A continuous map $f:\Sigma\rightarrow X$
is called {\it pseudo-holomorphic} if for any point $z_0\in
\Sigma$, the following is true:
\begin{enumerate}
\item There is a disc neighborhood of $z_0$ with a branched covering map
$br:z\rightarrow z^m$. (Here $m=1$ is allowed.)
\item There is a local chart $(V_{f(z_0)},G_{f(z_0)},\pi_{f(z_0)})$ of $X$
at $f(z_0)$ and a local lifting $\tilde{f}_{z_0}$ of $f$ in the sense that
$f\circ br=\pi_{f(z_0)}\circ\tilde{f}_{z_0}$.
\item $\tilde{f}_{z_0}$ is pseudo-holomorphic, i.e.,
$d \tilde{f}_{z_0}\circ j=J\circ d\tilde{f}_{z_0}$.
\end{enumerate}}\vskip 0.1in

We shall be interested in moduli spaces of pseudo-holomorphic maps from a
Riemann surface into an almost complex orbifold. The next lemma collects
a few analytic properties of them which are very useful to keep in mind.

    \vskip 0.1in
\noindent{\bf Lemma 2.1.4: }{\it Let $f:(\Sigma,j)\rightarrow
(X,J)$ be a non-constant pseudo-holomorphic map such that
$f(\Sigma)\cap\Sigma X\neq \emptyset$. Then there are two
possibilities:

\begin{enumerate}
\item There are finitely many distinct points $z_i\in\Sigma$, $i=1,\cdots,k$,
such that $f^{-1}(\Sigma X)=\{z_i\}$ and the image $f(\Sigma)$
intersects $\Sigma X$ at each $f(z_i)$ with a finite order
tangency. More precisely, suppose $\tilde{f}_{z_i}:D\rightarrow
V_{f(z_i)}$ is a local lifting of $f$; then $\tilde{f}_{z_i}(D)$
intersects with  finite order tangency the preimage of $\Sigma X$
in $V_{f(z_i)}$, which is a union of finitely many embedded
pseudo-holomorphic submanifolds.
\item The image $f(\Sigma)$ lies entirely in $\Sigma X$. In this case,
there is a set $\{z_1,\cdots,z_k\}\subset \Sigma$ and a connected
stratum of $\Sigma X=\widetilde{\Sigma X}_{gen}$ containing
$f(\Sigma\setminus\{z_1,\cdots,z_k\})$. Therefore, after
compactification, $f$ is lifted to a pseudo-holomorphic map
$f_1:\Sigma\rightarrow \widetilde{\Sigma X}$, and $f_1(\Sigma)$
intersects the singular set of $\widetilde{\Sigma X}$ at finitely
many isolated points with  finite order tangency.
\end{enumerate}}
    \vskip 0.1in

\noindent{\bf Proof:} Let $z\in \Sigma$ be any point such that
$f(z)=p$ lies in $\Sigma X$. By Definition 2.1.3, there is a local
lifting of $f$, $\tilde{f}:D\rightarrow V_p$, satisfying
$d\tilde{f}\circ j=J\circ d\tilde{f}$, where $z$ and
$\tilde{f}(0)$ are the origins of $D$ and $V_p$, respectively. The
preimage of $\Sigma X$ in $V_p$ is a union of finitely many
embedded pseudo-holomorphic submanifolds. Let $W$ be any
component; we can choose a complex coordinate system $u_i$,
$i=1,\cdots,n$, such that $u_i(\tilde{f}(0))=0$ and $u_i$,
$i=k+1,\cdots,n$, is a complex coordinate system for $W$ near
$\tilde{f}(0)$. Then the equation $d\tilde{f}\circ j=J\circ
d\tilde{f}$ is written as $$ \frac{\partial u_i}{\partial
s}+J(u_1,\cdots,u_n)\frac{\partial u_i} {\partial t}=0,
\hspace{2mm}, i=1,\cdots,n, $$ where $z=s+it$ is the complex
coordinate on $D$. Let $\Delta=\frac{\partial^2}{\partial
s^2}+\frac{\partial^2}{\partial t^2}$ be the standard Laplacian on
$D$; then the above equation implies $$ \Delta u_i=(\partial_t
J(u))\partial_s u_i-(\partial_s J(u))\partial_t u_i, \hspace{2mm},
i=1,\cdots, n, $$ where $u=(u_1,\cdots, u_n)$. On the other hand,
if we write $J(u)$ as a matrix $$\left(\begin{array}{cc} A(u) &
B(u)\\ C(u) & D(u)
\end{array} \right)
$$ where $A(u)$ is a $k\times k$ matrix, then
$B(0,\cdots,0,u_{k+1},\cdots,u_n)=0$ ($W$ is
$J$-pseudo-holomorphic), which implies that $|\frac{\partial
B}{\partial u_i}(u)|\leq C(|u_1|+\cdots+|u_k|)$ for
$i=k+1,\cdots,n$ with a constant $C$ for all $u$ with $|u|\leq
\delta_0$. If we denote $v=(u_1,\cdots,u_k)$, then there is a
constant $C^\prime$ such that $$ |\Delta v(z)|\leq
C^\prime(|v|+|\partial_s v|+|\partial_t v|) $$ holds for all $z\in
D$ with $|z|\leq \epsilon_0$ for some $\epsilon_0$. By
Hartman-Wintner's lemma [HW], $v$ is either identically zero or
$v(z)=az^m+O(|z|^{m+1})$ for some $0\neq a\in \C$. Therefore we
prove that either $\tilde{f}(D)$ is entirely in the preimage of
$\Sigma X$ in $V_p$, or intersects it with  finite order tangency.

Now suppose that $f(\Sigma)$ is entirely in $\Sigma X$. Then by
the canonical stratification of $X$, we have $\Sigma X=
\widetilde{\Sigma X}_{gen}$. The points in $\Sigma$ are divided
into two groups, $I$ and $II$, according to the following rule:
$z$ is in $I$ if a neighborhood of $z$ is mapped into the same
stratum of the canonical stratification as $z$ under $f$; $z$ is
in $II$  otherwise. Then it follows that $II$ consists of only
finitely many isolated points. This implies that the subset of
$\Sigma$ consisting of points in $I$ is connected, so that its
image under $f$ lies in a connected open stratum of
$\widetilde{\Sigma X}_{gen}$. After compactification, we obtain a
unique pseudo-holomorphic map $f_1:\Sigma\rightarrow
\widetilde{\Sigma X}$ such that $\pi\circ f_1=f$, where
$\pi:\widetilde{\Sigma X}\rightarrow \Sigma X$ is the canonical
resolution. Moreover, $f_1$ intersects the singular set of
$\widetilde{\Sigma X}$ at finitely many points with  finite order
tangency. \hfill $\Box$

\vspace{2mm}

In Gromov-Witten theory, it is of primary importance that the moduli spaces
under consideration admit certain compactifications. We shall consider
two classes of almost complex orbifolds for which we can prove the
compactness of moduli spaces. One of them is given in the following

\vskip 0.1in\noindent
    {\bf Definition 2.1.5: }{\it A $2$-form on an
orbifold $X$ is said to be {\it non-degenerate} if each local
lifting is non-degenerate. A {\it symplectic  orbifold} is an
orbifold equipped with a closed, non-degenerate, $2$-form
$\omega$, usually denoted by $(X,\omega)$. An almost complex
structure $J$ on a symplectic orbifold $(X,\omega)$ is said to be
$\omega$-tamed if $\omega(v,Jv)>0$ for any nonzero $v\in T_p X$.
It is called $\omega$-compatible if $\omega(\cdot,J\cdot)$ is a
Riemannian metric on $X$. The space of $\omega$-tamed almost
complex structures on $(X,\omega)$ is denoted by $\J_\omega$.}
\vskip 0.1in

    The same proof as in the smooth case yields
    \vskip 0.1in
\noindent{\bf Lemma 2.1.6: }{\it $\J_\omega$ is non-empty and
contractible.}
    \vskip 0.1in

Another class of almost complex orbifolds comes from algebraic geometry.
Suppose that $(X,J)$ is an almost complex orbifold. If the local lifting
of the almost complex structure $J$ on each uniformizing system $(V,G,\pi)$
is integrable, then $(V,J)$ becomes a complex manifold and $G$ acts on it
holomorphically. The quotient space $U=\pi(V)$ inherits an
analytic structure such that the map $\pi:V\rightarrow U$ is analytic.
As a consequence, $X$ is an analytic space. We call $X$ a {\it complex
orbifold}.

A continuous map $f: (\Sigma,j)\rightarrow (X,J)$ from a Riemann
surface into a complex orbifold $(X,J)$ is called {\it analytic}
if, for any local analytic function $g$ on $X$, the pull-back
$g\circ f$ is holomorphic on $\Sigma$.

    \vskip 0.1in
\noindent{\bf Lemma 2.1.7: }{\it  A continuous map
$f:(\Sigma,j)\rightarrow (X,J)$ is analytic if and only if it is
pseudo-holomorphic in the sense of Definition 2.1.3.}
    \vskip 0.1in

\noindent{\bf Proof:}
Suppose $f$ is pseudo-holomorphic. Then at each point $z_0\in \Sigma$,
there is a disc neighborhood of $z_0$ with a branched covering map
$br:z\rightarrow z^m$ (Here $m=1$ is allowed),
a  local chart $(V_{f(z_0)},G_{f(z_0)},\pi_{f(z_0)})$ of $X$
at $f(z_0)$ and a local pseudo-holomorphic lifting $\tilde{f}_{z_0}$ of $f$
in the sense that $f\circ br=\pi_{f(z_0)}\circ\tilde{f}_{z_0}$ and
$d \tilde{f}_{z_0}\circ j=J\circ d\tilde{f}_{z_0}$. A local function $g$ is
analytic if and only if $g\circ\pi_{f(z_0)}$ is holomorphic on $V_{f(z_0)}$.
So if $g$ is a local analytic function, then $g\circ\pi_{f(z_0)}\circ
\tilde{f}_{z_0}$ is
holomorphic, which implies that  $g\circ f\circ br$ is holomorphic. So
$g\circ f$ is holomorphic, and hence $f$ is analytic.

Suppose that $f$ is analytic. For any connected component of
$\widetilde{\Sigma X}$, its image under the canonical resolution
map $\pi:\widetilde{\Sigma X}\rightarrow X$ is an analytic
subvariety. So the image of $f$ is either contained in it or
intersects it at finitely many points. So there is a connected
component $X_0$ of $\widetilde{\Sigma X}$ such that $f$ is lifted
to an analytic map $f_1:\Sigma\rightarrow X_0$ with the property
that for all $z$ except finitely many $z_1,\cdots,z_k$, $f_1(z)$
is in the set of regular points of $X_0$. It follows easily from
this that $f$ is pseudo-holomorphic. \hfill $\Box$

The second class of almost complex orbifolds we shall consider is
the so-called {\it projective orbifold}.

    \vskip 0.1in
\noindent{\bf Definition 2.1.8: }{\it A complex orbifold $X$ is
called {\it projective} if it can be realized as a projective
variety.}

\subsection{Classification of compatible systems}

In this subsection, we shall study pseudo-holomorphic maps from a
Riemann surface into an almost complex orbifold in terms of
$C^\infty$ maps between orbifolds. For definitions or notations
regarding good maps and compatible systems, the reader is referred
to section 4.4. The main result is summarized in the following

    \vskip 0.1in
\noindent{\bf Proposition 2.2.1: }{\it For any pseudo-holomorphic
map $f$ from a Riemann surface $\Sigma$ of genus $g$ with $k$
marked points $\z=(z_1,z_2,\cdots,z_k)$ into a closed almost
complex orbifold $(X,J)$, there are finitely many orbifold
structures on $\Sigma$ whose singular set is contained in $\z$,
and for each of these orbifold structures there are finitely many
pairs $(\tilde{f},\xi)$, where $\tilde{f}$ is a good map whose
underlying map is $f$, and $\xi$ is an isomorphism class of
compatible systems of $\tilde{f}$. The total number is bounded
from above by a constant $C(X,g,k)$ depending only on $X,g,k$.}

    \vskip 0.1in

We will break the proof into several steps. First recall that a
complex orbicurve is a smooth complex curve with a reduced
orbifold structure. More precisely, we have the following

    \vskip 0.1in
\noindent{\bf Definition 2.2.2: }{\it A  complex orbicurve of
genus $g$ is a triple $(\Sigma,\z,\m)$, where $\Sigma$ is a smooth
complex curve of genus $g$, $\z=(z_1,z_2,\cdots,z_k)$ is a set of
distinct marked points on $\Sigma$, $\m=(m_1,m_2,\cdots,m_k)$ is a
$k$-tuple of integers with $m_i\geq 2$. $\Sigma$ is given an
orbifold structure as follows: at each point $z_i$, a disc
neighborhood of $z_i$ is uniformized by the branched covering map
$z\rightarrow z^{m_i}$.}
    \vskip 0.1in

We will be interested in good $C^\infty$ maps from an orbicurve
$(\Sigma,\z,\m)$ into an almost complex orbifold $X$ which induce
a pseudo-holomorphic map. We make two observations: First, since
pseudo-holomorphic maps (in the classical sense) have the unique
continuity property, given an orbifold structure of $\Sigma$ with
a pseudo-holomorphic map $f$, there is a unique $C^\infty$ map
$\tilde{f}$ as a germ of liftings of $f$. Secondly, let
$\tilde{f}:(\Sigma,\z,\m)\rightarrow X$ be a good map inducing a
pseudo-holomorphic map $f$; for any isomorphism class of
compatible systems $\xi$, if the group homomorphism it defines at
$z_i$, $\lambda_i:\Z_{m_i}\rightarrow G_{f(z_i)}$,  is not
monomorphic, then we can factor $\lambda_i$ through a monomorphism
$\lambda_i^\prime: \Z_{m_i^\prime}\rightarrow G_{f(z_i)}$,
redefine the orbifold structure at $z_i$ by $z\rightarrow
z^{m_i^\prime}$, and obtain a good map
$\tilde{f}^\prime:(\Sigma,\z,\m^\prime)\rightarrow X$ with an
isomorphism class of compatible systems $\xi^\prime$ such that the
group homomorphism of $\xi^\prime$ at $z_i$ is $\lambda_i^\prime$,
and the restriction of $\xi^\prime$ to $\Sigma\setminus\{z_i\}$ is
the same as that of $\xi$. Therefore, we will  only be considering
good maps with an isomorphism class of compatible systems whose
group homomorphism at each point is monomorphic.

First we show that for any pseudo-holomorphic map $f$ from a complex curve
$\Sigma$ into $(X,J)$, there is always an orbifold structure on $\Sigma$
with respect to which $f$ admits a $C^\infty$ lifting.

    \vskip 0.1in
\noindent{\bf Lemma 2.2.3: }{\it  Let $f:(\Sigma,j)\rightarrow
(X,J)$ be a non-constant pseudo-holomorphic map such that
$f(\Sigma)$ is not entirely in $\Sigma X$; then there is an
orbifold structure on $\Sigma$ making it into a complex orbicurve
$(\Sigma,\z,\m)$, and $f$ determines a unique $C^\infty$ map
$\tilde{f}:(\Sigma,\z,\m)\rightarrow X$ as a germ of $C^\infty$
liftings of $f$, which is regular, therefore is good with a unique
isomorphism class of compatible systems. Moreover, the group
homomorphism at each point is monomorphic.}\vskip 0.1in

\noindent{\bf Proof:} When $f(\Sigma)\subset X_{reg}$, $f$ is just
a pseudo-holomorphic map in the usual sense, and the claim is
trivial. Suppose $f(\Sigma)$ is not entirely in $\Sigma X$. Then
there are finitely many $z_i\in \Sigma$ such that $f(z_i)=p_i\in
\Sigma X$. We will define the orbifold structure at each $z_i$ as
follows. There is a branched covering map $br:z\rightarrow
z^{n_i}$, a chart $(V_{p_i},G_{p_i},\pi_{p_i})$ of $X$ at $p_i$,
and a pseudo-holomorphic map $\tilde{f}_i:D\rightarrow V_{p_i}$
such that $f\circ br=\pi_{p_i}\circ \tilde{f}_i$. Let $m_i$ be the
minimum of the $n_i$'s. We take $\z=(z_1,\cdots,z_k)$ where
$z_i\in f^{-1}(\Sigma X)$ with $m_i\neq 1$, and
$\m=(m_1,\cdots,m_k)$. Then the elliptic regularity and unique
continuity property of pseudo-holomorphic maps implies that $f$
admits a $C^\infty$ lifting $\tilde{f}:(\Sigma,\z,\m)\rightarrow
X$, which is obviously regular. By Lemma 4.4.11, $\tilde{f}$ is a
good map with a unique isomorphism class of compatible systems.
The minimality of $m_i$ implies that the group homomorphism at
each point is monomorphic. \hfill $\Box$

    For the case when $f(\Sigma)\subset\Sigma X$, there is a set
$\{z_1,\cdots,z_k\}$ such that
$f(\Sigma\setminus\{z_1,\cdots,z_k\})$ is contained in a connected
stratum of $\Sigma X=\widetilde{\Sigma X}_{gen}$, i.e., for $z\in
\Sigma\setminus\{z_1,\cdots,z_k\}$, let $p=f(z)$; there is a chart
$(V_p,G_p,\pi_p)$, a branched covering map $br:z\rightarrow z^m$,
a lifting $\tilde{f}_z$ such that $\pi_p\circ\tilde{f}_z=f\circ
br$, and $\tilde{f}_z(D)$ lies in the fixed-point set of $G_p$. It
is easily seen that $br$ can be taken trivial, i.e., $m=1$.  For
each $z_i$, there is a chart at $p_i=f(z_i)$, $(V_{p_i},
G_{p_i},\pi_{p_i})$, a branched covering map $br_i$, a lifting
$\tilde{f}_i$ such that $\pi_{p_i}\circ\tilde{f}_i=f\circ br_i$.
Now it is easily seen that given a pseudo-holomorphic map
$f:\Sigma\rightarrow X$, there is always an orbifold structure on
$\Sigma$ such that $f$ is lifted to a (unique) $C^\infty$ map.

\vspace{2mm}

Next we will give a classification of compatible systems of a good
pseudo-holomorphic map from a complex orbicurve into $(X,J)$, up
to isomorphism.

Let $\tilde{f}:(\Sigma,\z,\m)\rightarrow X$ be a good
pseudo-holomorphic map. Let $\xi_0$ and $\xi$ be two isomorphism
classes of compatible systems. We will fix $\xi_0$ for a moment.
The restriction to $\Sigma\setminus\z$ gives rise to two
isomorphism classes of compatible systems of $\tilde{f}$ on
$\Sigma\setminus\z$, which are represented by geodesic compatible
systems $\{\tilde{f}_{0,UU^\prime},\lambda_0\}$ and
$\{\tilde{f}_{UU^\prime},\lambda\}$, respectively. In addition, we
fix a base point $z_0\in \Sigma\setminus\z$ such that $f(z_0)$
lies in the main stratum. We denote $G_{f(z_0)}$ by $H_{z_0}$. We
remark that since each point in $\Sigma\setminus\z$ is regular,
for each inclusion between elements in $\U$, there is only one
injection. Now we take an element $U_0$ in $\U$ containing $z_0$
which is maximal in the sense that $U_0$ is not contained in any
other element in $\U$. By the choice of $z_0$, we see that
$\tilde{f}_{0,U_0U_0^\prime}=\tilde{f}_{U_0U_0^\prime}$. Each
automorphism of the uniformizing system of $U_0^\prime$, given by
an element of $H_{z_0}$, leaves $\tilde{f}_{U_0U_0^\prime}$ fixed.
We will introduce a basic operation on the compatible system
$\{\tilde{f}_{UU^\prime},\lambda\}$, called {\it conjugation},
which will only produce an isomorphic compatible system. The
operation is this: for an element $U^\prime$ in $\U^\prime$, we
take an automorphism $\delta$ of the uniformizing system of
$U^\prime$, and change $\tilde{f}_{UU^\prime}$ to
$\delta\circ\tilde{f}_{UU^\prime}$ and for each inclusion $i:
U_1\rightarrow U$, we change $\lambda(i)$ to
$\delta\circ\lambda(i)$, and for each inclusion $i:U\rightarrow
U_2$, we change $\lambda(i)$ to $\lambda(i)\circ\delta^{-1}$. For
any $U$ in $\U$ such that $U\subset U_0$, there is a unique
conjugation done to $U$ after which we have
$\tilde{f}_{0,UU^\prime}=\tilde{f}_{UU^\prime}$ and
$\lambda_0(i)=\lambda(i)$ where $i:U\rightarrow U_0$ is the
inclusion. Then the compatibility condition of compatible systems
$(4.4.1)$ ensures that for any $i:U_2\rightarrow U_1$ where $U_1$,
$U_2\subset U_0$, we have $\lambda_0(i)=\lambda(i)$. In other
words, we show that the two compatible systems agree on $U_0$.
Next we will try  to spread this out. We take an element
$U_1\in\U$ which is maximal and has non-empty intersection with
$U_0$, which is a geodesically convex topological disc. We pick a
$U\subset U_0\cap U_1$, and let $i:U\rightarrow U_1$; then there
is a unique conjugation done to $U_1$ after which we have
$\tilde{f}_{0,U_1U_1^\prime}=\tilde{f}_{U_1U_1^\prime}$ and
$\lambda_0(i)=\lambda(i)$. The compatibility condition $(4.4.1)$
ensures that $\lambda_0(j)=\lambda(j)$ for any inclusion
$j:W\rightarrow U_1$ where $W\in \U$ and $W\subset U$, and since
$U_0\cap U_1$ is connected, we can drop the condition $W\subset
U$. We then proceed to do conjugations to all $U\subset U_1$ so
that the two compatible systems agree on $U_0\cup U_1$. We can
continue to do this as long as the new element we take in $\U$
intersects the union of the old ones in a connected subset, and
whenever this fails, we have a loop, and we are forced to do a
conjugation to an element in $\U$ that we have done  before. But
if the loop is homotopically trivial, no problem will be caused.
Now we take a loop $\gamma$ based at $z_0$, going along $\gamma$
once; we will be forced to do a conjugation twice to an element
$U\subset U_0$ so that the equation $\lambda_0(i)=\lambda(i)$ may
not hold any longer, where $i:U\rightarrow U_0$. We take
$\lambda_0(i)\circ(\lambda(i))^{-1}$ which is an element in
$H_{z_0}$, and we denote it by $g_\gamma$. One can verify:
\begin{itemize}
\item If we do a conjugation to $U_0$ with automorphism $\delta$,
then $g_\gamma$ will be changed to $\delta\circ g_\gamma \circ \delta^{-1}$.
\item If $\gamma_1$ is homotopic to $\gamma_2$, then
$g_{\gamma_1}=g_{\gamma_2}$.
\item The assignment $[\gamma]\rightarrow g_\gamma$ defines a homomorphism
$\theta_{\xi_0,\xi}:\pi_1(\Sigma\setminus\z,z_0)\rightarrow H_{z_0}$.
\item If $\theta_{\xi_0,\xi_1}=\theta_{\xi_0,\xi_2}$,  then
the restiction of $\xi_1$ to $\Sigma\setminus\z$ is
isomorphic to that of $\xi_2$.
\end{itemize}

The homomorphism $\theta_{\xi_0,\xi}$ measures the difference
between $\xi_0$ and $\xi$ restricted to $\Sigma\setminus\z$. When
$f(\Sigma)$ lies in a uniformized open set, we can get an absolute
measurement. A similar argument shows that

    \vskip 0.1in
\noindent{\bf Lemma 2.2.4: }{\it  If $f(\Sigma)$ lies in a
neighborhood uniformized by $(V,G,\pi)$, then each isomorphism
class of compatible systems determines a conjugacy class of
homomorphisms $\theta_\xi:\pi_1(\Sigma\setminus\z,z_0)\rightarrow
G$ such that $\theta_{\xi_1}=\theta_{\xi_2}$ implies that
restricting to $\Sigma\setminus \z$, $\xi_1$ is isomorphic to
$\xi_2$.}
    \vskip 0.1in

\noindent{\bf Proof:} Given a compatible system
$\{\tilde{f}_{UU^\prime},\lambda\}$. Consider the inclusion from
$\cup_{U^\prime\in \U^\prime}U^\prime$ into $\pi(V)$.  By taking
composition, we obtain a compatible system
$\{\tilde{f}_{WW^\prime},\tau\}$ which is isomorphic to
$\{\tilde{f}_{UU^\prime},\lambda\}$ and $\W^\prime$ consists of a
single element $\pi(V)$ uniformized by $(V,G,\pi)$. Moreover, each
$\tau(i)$ is an automorphism of $(V,G,\pi)$, given by an element
of $G$. Now we take a loop $\gamma$ in $\Sigma\setminus\z$ based
at $z_0$, and cover $\gamma$ by a set of discs in $\W$, say
$W_1,\cdots,W_k$, such that only two adjacent discs and $W_1$,
$W_k$ have non-empty intersection. For any $i\in\{1,\cdots,k\}$,
take a point $z_i\in W_i\cap W_{i+1}$ (if $i=k$, then $i+1=1$) and
a disc neighborhood $D_i\in\W$ of it such that $D_i\subset W_i\cap
W_{i+1}$. Let $a_i:D_i\rightarrow W_i$, $b_i:D_i\rightarrow
W_{i+1}$ be the inclusions respectively; we put
$g_i=\lambda_1(b_i)\circ\lambda_1^{-1}(a_i)$ in $G$, and
$g_\gamma=g_k\circ\cdots\circ g_1$. Then one can verify that the
assignment $\gamma\rightarrow g_\gamma$ is well-defined and gives
the required homomorphism
$\theta_\xi:\pi_1(\Sigma\setminus\z,z_0)\rightarrow G$. \hfill
$\Box$
    \vskip 0.1in
    \noindent
    {\bf Definition 2.2.5: }{\it $\theta_{\xi}$ is called the characteristic
of $\xi$. $\theta_{\xi_0, \xi_1}$ is called the difference
characteristic of the pair $(\theta_0, \theta_1)$.}
    \vskip 0.1in

The next lemma deals with compatible systems defined near an orbifold
point.

    \vskip 0.1in
\noindent{\bf Lemma 2.2.6: }{\it Suppose
$\tilde{f}:D\setminus\{0\}\rightarrow U_p$ is a good
pseudo-holomorphic map which can be extended continuously over $D$
sending $0$ to $p$, where $U_p$ is a geodesic neighborhood with a
geodesic chart $(V_p,G_p,\pi_p)$, and suppose each point in
$f(D\setminus\{0\})$ has the same orbit type given by the
conjugacy class $(H)$ of subgroups of $G_p$. For any isomorphism
class of compatible systems $\xi$ of $\tilde{f}$, let
$\theta_\xi:\pi_1(D\setminus\{0\},\ast)=\Z\rightarrow G_p$ be
given by $x\rightarrow g_\xi$ where $x$ is the positively oriented
generator; then after choosing a representative $H$, $g_\xi$ lies
in the normalizor $N_{G_p}(H)$ of $H$ in $G_p$. Let $m$ be the
order of $g_\xi$ in $G_p$. Then we can give an orbifold structure
at $\{0\}$ to $D$ by $z\rightarrow z^m$, so that $\tilde{f}$ is
extended to a good map $\tilde{f}^\prime$ from $(D,(0),(m))$ to
$U_p$ with a unique isomorphism class of compatible systems
$\xi^\prime$, such that the group homomorphism of $\xi^\prime$ at
$0$ is given by $e^{2\pi i/m}\rightarrow g_\xi$, and we have
$\theta_\xi=\theta_{\xi^\prime}$.}
    \vskip 0.1in

\noindent{\bf Proof:} Choose a representative $H$ of $(H)$. From
the construction of the homomorphism
$\theta_\xi:\pi_1(D\setminus\{0\},\ast)\rightarrow G_p$ (cf. Lemma
2.2.4), and the assumption that each point in $f(D\setminus\{0\})$
has the same orbit type given by the conjugacy class $(H)$, we see
that the element $g_\xi\in G_p$ must be in the normalizer
$N_{G_p}(H)$ of $H$ in $G_p$. Consider the preimage of
$f(D\setminus\{0\})$ in $V_p$ under $\pi_p$. Since each point in
$f(D\setminus\{0\})$ has the same orbit type given by the
conjugacy class of the subgroup $H$ of $G_p$, we consider the part
$W$ of the preimage lying in the fixed-point set of $H$, $V^H$, in
$V_p$. We take a non-trivial loop in $D\setminus\{0\}$ such that
$f$ is an embedding in a neighborhood of it; then it is easily
seen that one of the components of its preimage in $W$ cyclically
covers it under $\pi_p$, which is given by the action of $[g_\xi]$
in $N_{G_p}(H)/H$. Therefore there is a lifting of $f$, $\bar{f}:
D\setminus\{0\}\rightarrow V_p$, a branched covering map $br:
z\rightarrow z^m$ where $m$ equals the order of $g_\xi$, such that
$\pi_p\circ\bar{f}=f\circ br$, and $\bar{f}$ is equivariant under
the homomorphism  $e^{2\pi i/m}\rightarrow g_\xi$. By removability
of the singularity, $\bar{f}$ extends to a pseudo-holomorphic map
$D\rightarrow V_p$, still denoted by $\bar{f}$. Now we need to
extend the compatible system $\xi$ defined over $D\setminus\{0\}$
to a compatible system $\xi^\prime$ on $D$, including
$(\bar{f},\rho)$ where $\rho:\Z_m\rightarrow G_p$ is given by
$e^{2\pi i/m}\rightarrow g_\xi$. Let $\xi$ be given by
$\{\tilde{f}_{WW^\prime},\lambda\}$. We take an element $W\in \W$
and an injection $i$ from $W$ into $D$, and we define $\lambda(i)$
to be an injection from the uniformizing system of $W^\prime$ into
$V_p$ such that $\bar{f}\circ i=\lambda(i)\circ
\tilde{f}_{WW^\prime}$. Then we can extend out to define
$\lambda(j)$ for injections $j:W\rightarrow D$ nearby $i$. One can
verify that $\lambda(i)$ can be chosen so that we can make it
consistent when we move around a loop. This gives the extension
$\xi^\prime$ of $\xi$. The statement that $\xi^\prime$ is unique
and $\theta_\xi=\theta_{\xi^\prime}$ follows from the fact that
the compatible system on $D\setminus\{0\}$ induced by any map
$\tilde{f}: D\rightarrow V_p$ with homomorphism $e^{2\pi i/m}
\rightarrow g$ has the corresponding
$\theta:\pi_1(D\setminus\{0\},\ast) \rightarrow G_p$ given by
$x\rightarrow g$. \hfill $\Box$

With Lemma 2.2.6 understood, one can regard the characteristics
$\theta_\xi$ or $\theta_{\xi_0,\xi_1}$ as defined on the orbifold
fundamental group $\pi_1^{orb}(\Sigma,\z,\m,z_0)$ after giving a
reduced orbifold structure on $(\Sigma,\z)$ according to Lemma 2.2.6.
In the case when $f(\Sigma)$ lies in a global quotient $U=V/G$, given any
compatible system $\xi$, let $G_\xi$ be the image of $\theta_\xi$ in
$G$. Then there is a smooth Riemann surface $\Sigma_\xi$ with an action
of $G_\xi$, such that $\pi_1(\Sigma_\xi)$ is the kernel of
$\theta_\xi$ and $(\Sigma_\xi,G_\xi)$ is a uniformizing system of
$(\Sigma,\z,\m)$. Moreover, the compatible system $\xi$ is induced from
a natural equivariant map $\tilde{f}:\Sigma_\xi\rightarrow V$ (recall
that here $G_\xi$ is a subgroup of $G$).

    \vskip 0.1in
\noindent{\bf Corollary 2.2.7: }{\it  Let $f:\Sigma\rightarrow X$
be a pseudo-holomorphic map, let $\z$ be a set of finitely many
points in $\Sigma$. Then there are only finitely many orbifold
structures on $\Sigma$ having $\z$ as singular set, such that $f$
is lifted to  a (unique) good map, and for each such  orbifold
structure there are only finitely many isomorphism classes of
compatible systems. The total number of orbifold structures and
isomorphism classes of compatible systems is bounded by a number
depending only on $g,k$ and $X$, where $g$ is the genus of
$\Sigma$, and $k$ equals the order of the set $\z$.}
    \vskip 0.1in

    \vskip 0.1in
\noindent{\bf Remark 2.2.8: }{\it  In string theory, the
pseudo-holomorphic map $f:(\Sigma,\z,\m)\rightarrow X$ near a
point $z_k\in\z$ is regarded as  open string propagation on
$V_{f(z_k)}$ with the two ends of the open string satisfying
$x(2\pi)=g \cdot x(0)$ for some element $g\in G_{f(z_k)}$.
Likewise,  closed string propagation in the quotient space $X=Y/G$
may be regarded as  open string propagation in $Y$ with a boundary
condition $x(2\pi)=g\cdot x(0)$ for some $1\neq g\in G$. These are
called {\it twisted boundary conditions} in string theory. See
[DHVW]. So when $f(\Sigma)$ lies entirely in $\Sigma X$, it may
happen that although the open string propagation is the same
geometrically, it is attached with different twisted boundary
conditions. On the other hand, mathematically, the different
choices of twisted boundary conditions correspond to different
choices of compatible system, which controls the deformation of
the pseudo-holomorphic map.}

The last ingredient in the proof of Proposition 2.2.1 is provided by

    \vskip 0.1in
\noindent{\bf Lemma 2.2.9: }{\it  For any pseudo-holomorphic map
$f:\Sigma\rightarrow X$, there is an orbifold structure on
$\Sigma$ with respect to which $f$ admits a good $C^\infty$
lifting.}
    \vskip 0.1in

\noindent{\bf Proof:} We can assume the case that $f(\Sigma)$ lies
entirely in the singular set $\Sigma X$, and there is a set
$\{z_1,\cdots,z_k\}$ such that
$f(\Sigma\setminus\{z_1,\cdots,z_k\})$ is contained in a connected
stratum of $\Sigma X=\widetilde{\Sigma X}_{gen}$, so that for
$z\in \Sigma\setminus\{z_1,\cdots,z_k\}$, $f$ can be lifted to a
chart $(V_p,G_p,\pi_p)$ where $p=f(z)$ without using any branched
covering map. For each $z_i$, there is a chart $(V_{p_i},
G_{p_i},\pi_{p_i})$ at $p_i=f(z_i)$, a branched covering map
$br_i$, and a lifting $\tilde{f}_i$ such that
$\pi_{p_i}\circ\tilde{f}_i=f\circ br_i$. There is a subgroup $H_i$
of $G_{p_i}$ and an element $a_i\in N_{G_{p_i}}(H_i)/H_i$ such
that $\tilde{f}_i$ is $a_i$-equivariant. We take any $g_i\in
N_{G_{p_i}}(H_i)$ such that $[g_i]=a_i$, and let $m_i$ be the
order of $g_i$. We can let $br_i$ be given by $z\rightarrow
z^{m_i}$. Let $\tilde{f}_i$ be the lifting of $f$ that is
$\lambda_i$-equivariant,  where $\lambda_i: e^{2\pi
i/m_i}\rightarrow g_i$. On the other hand, the general
construction of geodesic compatible systems can be applied here to
produce a collection of maps $\{\tilde{f}_{UU^\prime}\}$ including
these $\tilde{f}_i$'s, but with the assignment $\lambda$ missing.
However, the homomorphisms $\lambda_i$ define a compatible system
locally near $z_i$, the idea of the proof is to extend  out over
the whole $\Sigma$ at the price of introducing some extra orbifold
points.

Here is the basic observation: Suppose a compatible system is
already built on an open subset $\Sigma_0$ of $\Sigma$, and let
$U\in\U$ which is not contained in $\Sigma_0$. We want to extend
the compatible system to $\Sigma_0\cup U$. This consists of two
steps. The first one is to assign an injection $\lambda(i)$
between the uniformizing systems of $U_1$ and $U$ to any inclusion
$i:U_1\rightarrow U$ where $U_1\subset \Sigma_0\cap U$ and
$U_1\in\U$, which extends the old assignment and satisfies the
compatibility condition $(4.4.1)$. The second one is to extend
over for any inclusion $U_1\rightarrow U$. This can be done as
long as the intersection $U\cap \Sigma_0$ is simply connected (it
could be disconnected). So after removing finitely many points, a
compatible system can be constructed on $\Sigma$, then adding
these points as orbifold points we extend the compatible system to
the whole $\Sigma$. We remark that these extra orbifold points can
be squeezed into one point, and then absorbed into one of the
$z_i$'s. \hfill $\Box$

We end this subsection with a classification of compatible systems
of a constant (pseudo-holomorphic) map. Let
$f:(\Sigma,\z)\rightarrow (X,J)$ be a constant map with
$f(\Sigma)=\{p\}$. The matter is trivial if $p\in X_{reg}$. So we
assume the case when $p$ is in $\Sigma X$. The orbifold structure
of $X$ at $p$ determines a faithful representation $\rho:
G_p\rightarrow Aut(\C^n)$ (here $2n=\dim X$). We know that an
isomorphism class of compatible systems of $f$ is determined by a
conjugacy class of homomorphisms
 $\pi_1 (\Sigma\setminus\z,\ast)\rightarrow G_p$. We will next show
 that

    \vskip 0.1in
\noindent{\bf Proposition 2.2.10: }{\it  Given a conjugacy class of
homomorphisms $\theta:\pi_1(\Sigma\setminus\z,\ast)\rightarrow
G_p$, there is a unique orbifold structure on $\Sigma$ and a
(unique) good map $\tilde{f}$ with a unique isomorphism class of
compatible systems determined by its characteristic $\theta$.} \vskip 0.1in

\noindent{\bf Proof:} We need a digression first, recalling from
[CR1]. Let $(\Sigma,\z,\m)$ be a closed 2-dimensional orbifold,
where $\z=(z_1,\cdots,z_k)$ and $\m=(m_1,\cdots,m_k)$. Let $E$ be
a complex orbifold bundle of rank $n$ over $(\Sigma,\z,\m)$. Then
at each singular point $z_i$, $i=1,\cdots,k$, $E$ determines a
representation $\rho_i:\Z_{m_i}\rightarrow Aut(\C^n)$ so that over
a disc neighborhood of $z_i$, $E$ is uniformized by
$(D\times\C^n,\Z_{m_i},\pi)$ where the action of $\Z_{m_i}$ on
$D\times\C^n$ is given by $$ e^{2\pi i/m_i}\cdot (z,w)=(e^{2\pi
i/m_i}z,\rho_i(e^{2\pi i/m_i})w) $$ for any $w\in\C^n$. Each
representation $\rho_i$ is uniquely determined by  $n$-tuple of
integers $(m_{i,1},\cdots,m_{i,n})$ with $0\leq m_{i,j}<m_i$, as
it is given by the matrix $$ \rho_i(e^{2\pi i/m_i})= diag(e^{2\pi
im_{i,1}/m_i},\cdots, e^{2\pi im_{i,n}/m_i}). $$ Over the
punctured disc $D_i\setminus\{0\}$ at $z_i$, $E$ inherits a
specific trivialization from $(D\times\C^n,\Z_{m_i},\pi)$ as
follows: We define a $\Z_{m_i}$-equivariant map
$\Psi_i:D\setminus\{0\}\times\C^n\rightarrow
D\setminus\{0\}\times\C^n$ by $$ (z,w_1,\cdots,w_n)\rightarrow
(z^{m_i}, z^{-m_{i,1}}w_1,\cdots,z^{-m_{i,n}}w_n), $$ where
$\Z_{m_i}$ acts trivially on the second
$D\setminus\{0\}\times\C^n$. Hence $\Psi_i$ induces a
trivialization $\psi_i: E_{D_i\setminus\{0\}} \rightarrow
D_i\setminus\{0\}\times\C^n$. We can extend the smooth complex
vector bundle $E_{\Sigma\setminus\z}$ over $\Sigma\setminus\z$ to
a smooth complex vector bundle over $\Sigma$ by using these
trivializations $\psi_i$. We call the resulting complex vector
bundle the {\it de-singularization} of $E$, and denote it by
$|E|$. End of digression.

We denote by $T_p$ the tangent space of $V_p$ at $p$, and by
$\Sigma^\prime$ the universal cover of $\Sigma\setminus\z$. We let
$E$ be the complex vector bundle over $\Sigma\setminus\z$ obtained
from $\Sigma^\prime\times T_p$ through
$\theta:\pi_1(\Sigma\setminus\z,\ast)\rightarrow G_p$ and the
action of $G_p$ on $T_p$. Then $E$ naturally extends to an
orbifold bundle over $(\Sigma,\z,\m)$, still denoted by $E$, where
$\m=(m_i)$ and $m_i$ equals the order of $\theta(x_i)$ in $G_p$
for the positively oriented generator $x_i$ around $z_i\in\z$.
Take a $C^\infty$ section $\tilde{s}$ of $E$ as follows. Over each
uniformizing system $(D\times C^n,\Z_{m_i},\pi)$ of $E$ at
$z_i\in\z$, $\tilde{s}$ is given by
$\tilde{s}_i(z)=(\epsilon_{i,1} z^{m_{i,1}},\cdots,\epsilon_{i,n}
z^{m_{i,n}})$. We choose a generic non-zero
$(\epsilon_{i,1},\cdots,\epsilon_{i,n})$ such that $\tilde{s}_i$
meets the fixed-point set of $G_p$ at isolated points. Through the
trivialization $\psi_i$, each $\tilde{s}_i$ is identified with the
constant section $(\epsilon_{i,1},\cdots,\epsilon_{i,n})$ in the
trivialization of $|E|$ at $z_i$. So we let $\tilde{s}$ equal a
smooth section of $|E|$ extending the constant section
$(\epsilon_{i,1},\cdots,\epsilon_{i,n})$ near each $z_i$ such that
over only finitely many points $\tilde{s}$ lies in the fixed-point
set of $G_p$ in the tangent space of $V_p$ at $p$. This can be
done since the fixed-point set is of codimension at least two. Now
we consider a $C^\infty$ map $\tilde{F}:(\Sigma,\z,\m)\times [0,1]
\rightarrow U_p$ given by $(z,t)\rightarrow Exp_p\circ
(t\cdot\tilde{s}(z))$, where $Exp_p:TX_p\rightarrow U_p$ is the
exponential map at $p$. From the construction of $\tilde{s}$, we
know that $\tilde{F}$ is  a good map with a unique isomorphism
class of compatible systems, which induces a compatible system of
$f$. One can check that it is the one determined by $\theta$.
\hfill $\Box$

    \vskip 0.1in
\noindent{\bf Remark 2.2.11: }{\it  In what follows we will
consider pseudo-holomorphic maps from a marked Riemann surface
$(\Sigma,\z)$ (or more general a marked nodal Riemann surface)
into an almost complex orbifold $(X,J)$. We introduce the
following convention: We see that in Definition 2.1.3 there are
only finitely many points in $\Sigma$ for which the multiplicity
$m$ is greater than one. It is our convention that these point are
in $\z$ (i.e. marked). On the other hand, we will call an orbifold
structure on $\Sigma$ with the set of orbifold points contained in
$\z$ together with an isomorphism class of compatible systems of
the $C^\infty$ lifting $\tilde{f}$ of $f$ with respect to this
orbifold structure a {\it twisted boundary condition} of $f$,
borrowing  terminalogy from physics. Then Proposition 2.2.1 can be
rephrased as saying that the set of twisted boundary conditions of
a pseudo-holomorphic map from a marked Riemann surface of genus
$g$ with $k$ marked points into a closed almost complex orbifold
$X$ is non-empty and its cardinality is bounded by a number
$C(X,g,k)$ depending only on $X$, $g$, $k$.} \hfill $\Box$

\subsection{The moduli space of orbifold stable maps}

We first recall
    \vskip 0.1in
\noindent{\bf Definition 2.3.1: }{\it A {\it nodal curve with $k$
marked points is a pair $(\Sigma,\z)$ of a connected topological
space $\Sigma=\bigcup\pi_{\Sigma_\nu}(\Sigma_\nu)$ where
$\Sigma_\nu$ is a smooth complex curve  and
$\pi_\nu:\Sigma_\nu\rightarrow \Sigma$ is a continuous map, and
$\z=(z_1,\cdots,z_k)$ are $k$-distinct  points in $\Sigma$ with
the following properties.
\begin{itemize}
\item For each $z\in\Sigma_\nu$, there is a neighborhood of it such that
the restriction of $\pi_\nu:\Sigma_\nu\rightarrow \Sigma$ to this neighborhood is
a homeomorphism to its image.
\item For each $z\in\Sigma$, we have $\sum_\nu \#\pi_\nu^{-1}(z)\leq 2$.
\item $\sum_\nu \#\pi_\nu^{-1}(z_i)=1$ for each $z_i\in\z$.
\item The number of complex curves $\Sigma_\nu$ is finite.
\item The set of nodal points $\{z|\sum_\nu \#\pi_\nu^{-1}(z)=2\}$ is  finite.
\end{itemize}}\vskip 0.1in}
\vskip 0.1in

A point $z\in\Sigma_\nu$ is called {\it singular} if
$\sum_\omega\#\pi_\omega^{-1}(\pi_\nu(z))=2$. A point
$z\in\Sigma_\nu$ is said to be a {\it marked point} if
$\pi_\nu(z)=z_i\in\z$. Each $\Sigma_\nu$ is called a {\it
component} of $\Sigma$. Let $k_\nu$ be the number of points on
$\Sigma_\nu$ which are either singular or marked, and $g_\nu$ be
the genus of $\Sigma_\nu$; a nodal curve $(\Sigma,\z)$ is called
{\it stable} if $k_\nu+2g_\nu\geq 3$ holds for each component
$\Sigma_\nu$ of $\Sigma$.

A map $\vartheta:\Sigma\rightarrow\Sigma^\prime$ between two nodal
curves is called an {\it isomorphism} if it is a homeomorphism and
if it can be lifted to biholomorphisms
$\vartheta_{\nu\omega}:\Sigma_\nu\rightarrow \Sigma^\prime_\omega$
for each component $\Sigma_\nu$ of $\Sigma$. If $\Sigma$,
$\Sigma^\prime$ have marked points $\z=(z_1,\cdots,z_k)$ and
$\z^\prime=(z_1^\prime,\cdots,z_k^\prime)$ then we require
$\vartheta(z_i) =z_i^\prime$ for each $i$. Let $Aut(\Sigma,\z)$ be
the group of automorphisms of $(\Sigma,\z)$.

Each nodal curve $(\Sigma,\z)$ is canonically associated with a
graph $T_\Sigma$ as follows. The vertices of $T_\Sigma$ correspond
to the components of $\Sigma$ and for each pair of components
intersecting each other in $\Sigma$ there is an edge joining the
corresponding two vertices. For each point $z\in\Sigma$ such that
$\#\pi_\nu^{-1}(z)=2$, there is an edge joining the same vertex
corresponding to $\Sigma_\nu$. For each marked point, there is a
half open edge (tail) attaching to the vertex. The graph
$T_\Sigma$ is connected since $\Sigma$ is connected. We can smooth
out all the nodal points to obtain a smooth surface. Its genus is
called arithmetic genus of $\Sigma$. The arithmetic genus can be
computed by the formula

$$g=\sum_\nu g_\nu +rank H_1(T;\Q).$$

For any real numbers $t\geq 0$, $r>0$, set
$X(t,r)=\{(x,y)\in\C^2|||x||,||y||<r,xy=t\}$. We fix an action of
$\Z_m$ on $X(t,r)$ for any $m\geq 1$ by $e^{2\pi i/m}\cdot
(x,y)=(e^{2\pi i/m}x,e^{-2\pi i/m}y)$. Observe that the (branched)
covering map $X(t,r)\rightarrow X(t^m,r^m)$ given by $(x,y)\rightarrow
(x^m,y^m)$ is $\Z_m$-invariant, hence $(X(t,r),\Z_m)$ can be regarded
as a ``uniformizing system'' of $X(t^m,r^m)$.

    \vskip 0.1in
\noindent{\bf Definition 2.3.2: }{\it A nodal orbicurve} is a
nodal marked  curve $(\Sigma,\z)$ with an orbifold structure as
follows:
\begin{itemize}
\item The singular set $\z_\nu=\Sigma\Sigma_\nu$ (in the sense of orbifolds)
of each component $\Sigma_\nu$ is contained in the set of marked
points and nodal points $z$.
\item A disc neighborhood of a marked point is uniformized by a branched
covering map $z\rightarrow z^{m_i}$.
\item A neighborhood of a nodal point is uniformized by
$(X(0,r_j),\Z_{n_j})$.
\end{itemize}
 Here $m_i$ and $n_j$ are allowed to be
equal to one, i.e., the corresponding orbifold structure is
trivial there. We denote the corresponding nodal orbicurve by
$(\Sigma,\z,\m,\n)$ where $\m=(m_1,\cdots,m_k)$ and $\n=(n_j)$.

An {\it isomorphism} between two nodal orbicurves
$\tilde{\vartheta}: (\Sigma,\z,\m,\n)\rightarrow
(\Sigma^\prime,\z^\prime,\m^\prime,\n^\prime)$ is a collection of
$C^\infty$ isomorphisms $\tilde{\vartheta}_{\nu\omega}$ between
the orbicurves $\Sigma_\nu$ and $\Sigma_\omega^\prime$ which
induces an isomorphism $\vartheta:(\Sigma,\z)\rightarrow
(\Sigma^\prime,\z^\prime)$. The {\it group of automorphisms} of a
nodal orbicurve $(\Sigma,\z,\m,\n)$ is denoted by
$Aut(\Sigma,\z,\m,\n)$. It is easily seen that
$Aut(\Sigma,\z,\m,\n)$ is a subgroup of $Aut(\Sigma,\z)$ of finite
index.

\vspace{2mm}

\noindent{\bf Definition 2.3.3:} Let $(X,J)$ be an almost complex
orbifold. An {\it orbifold stable map} is a triple
$(f,(\Sigma,\z,\m,\n),\xi)$ described as follows:
\begin{enumerate}
\item $f$ is a continuous map from a nodal
orbicurve $(\Sigma,\z)$ into $X$ such that each
$f_\nu=f\circ\pi_\nu$ is a pseudo-holomorphic map from
$\Sigma_\nu$ into $X$.
\item $\xi$ is an isomorphism class of compatible structures
defined in the same way as in the smooth case.
\item At both marked and nodal points, the induced homomorphism on the local group
is injective.
\item Let $k_\nu$ be the order of the set $\z_\nu$, namely the
number of points on $\Sigma_\nu$ which are singular (i.e. nodal or
marked ); if $f_\nu$ is a constant map, then $k_\nu+ 2g_\nu\geq
3$.
\end{enumerate}
(We will call $\xi$ a twisted boundary condition of
$f:(\Sigma,\z)\rightarrow X$; cf. Remark 2.2.11.)

We introduce an equivalence relation amongst the set of stable
maps as follows: two stable maps $(f,(\Sigma,\z),\xi)$ and
$(f^\prime,(\Sigma^\prime,\z^\prime),\xi^\prime)$ are {\it
equivalent} if there exists an isomorphism $\vartheta:
(\Sigma,\z,\m,\n)\rightarrow
(\Sigma^\prime,\z^\prime,\m^\prime,\n^\prime)$ such that
${f}^\prime\circ\vartheta={f}$, and the compatible systems defined
by $\xi^\prime$ pull back via $\vartheta$ to compatible systems
isomorphic to the ones defined by $\xi$ (we write this as
$\xi^\prime\circ\vartheta=\xi$). The {\it automorphism group} of a
stable map $(f,(\Sigma,\z,\b,\n),\xi)$, denoted by
$Aut(f,(\Sigma,\z),\xi)$, is defined by $$
Aut(f,(\Sigma,\z),\xi)=\{{\vartheta}\in Aut(\Sigma,\z,\m,\n)|
{f}\circ\vartheta={f}, \xi\circ\vartheta=\xi\}. $$
    We often drop $\m,\n$ to simplify the notation.

The proof of the following lemma is routine and is left to the
reader.

    \vskip 0.1in
\noindent{\bf Lemma 2.3.4: }{\it The automorphism group of an
orbifold stable map is finite.}\vskip 0.1in

    In the smooth case, the automorphism group is the only source of
orbifold structure of the moduli space of stable maps. In the
orbifold case, orbifold structure introduces an additional
orbifold structure to the moduli space of orbifold stable maps
(see formula (3.2.1)).

Given a stable map $(f,(\Sigma,\z), \xi)$, there is an associated
homology class $f_\ast([\Sigma])$ in $H_2(X;\Z)$ defined by
$f_\ast([\Sigma])=\sum_\nu(f\circ\pi_\nu)_\ast[\Sigma_\nu]$. On
the other hand, for each marked point $z$ on $\Sigma_\nu$, say
$\pi_\nu(z) =z_i\in\z$, $\xi_\nu$ determines, by the group
homomorphism at $z$, a conjugacy class $(g_i)$, where $g_i\in
G_{f(z_i)}$. We thus have a map $ev$ sending each (equivalence
class of) stable map into $\widetilde{X}^k$ by
$(f,(\Sigma,\z),\xi)\rightarrow
((f(z_1),(g_1)),\cdots,(f(z_k),(g_k)))$. Here
$\widetilde{X}=\{(p,(g)_{G_p})|p\in X,g\in G_p\}$, cf. [CR1]. Its
components are called twisted sectors and are indexed by the
equivalence class of the conjugacy class $(g)_{G_p}$ under
injection (See \cite{CR1} for details.  Let $\x=\prod_i X_{(g_i)}$
be a connected component in $\widetilde{X}^k$.

    \vskip 0.1in
\noindent{\bf Definition 2.3.5: }{\it A stable map
$(f,(\Sigma,\z), \xi)$ is said to be of type $\x$ if
$ev((f,(\Sigma,\z), \xi))\in \x$. Given a homology class $A\in
H_2(X;\Z)$, we let $\overline{\M}_{g,k}(X,J,A,\x)$ denote the
moduli space of equivalence classes of orbifold stable maps of
genus $g$, with $k$ marked points, and of homology class $A$ and
type $\x$, i.e., $$
\overline{\M}_{g,k}(X,J,A,\x)=\{[(f,(\Sigma,\z),
\xi)]|g_\Sigma=g,\#\z=k, f_\ast([\Sigma])=A,
ev((f,(\Sigma,\z),\xi))\in\x\}. $$}

The rest of this subsection is devoted to giving a topology on
$\overline{\M}_{g,k}(X,J,A,\x)$ and to proving that the moduli
space is compact when $(X,J)$ is a compact symplectic orbifold or
a projective orbifold.

The set of all isomorphism classes of stable curves of genus $g$
with $k$ marked points, denoted by $\overline{\M}_{g,k}$, is
called the {\it Deligne-Mumford compactification} of the muduli
space $\M_{g,k}$ of Riemann surfaces of genus $g$ with $k$ marked
points (assuming $k+2g\geq 3$). The following differential
geometric description of $\overline{\M}_{g,k}$ is standard.

The moduli space $\overline{\M}_{g,k}$ admits a stratification
which is indexed by the combinatorial types of the stable curves.
More precisely, we consider a connected graph $T$ together with a
non-negative integer $g_{\nu}$. Let  $k_{\nu}$ be the number of
edges containing $\nu$ (we count twice the edges both of whose
vertices are $\nu$).  Then the data is required to satisfy $$
k_{\nu}+2g_\nu\geq 3, \hspace{2mm}\mbox{and} \hspace{2mm} \sum_\nu
g_\nu +rank \ H_1(T;\Q)=g. $$ Let $Comb(g,k)$ be the set of all
such objects $(T,(g_\nu))$. For each element $(\Sigma,\z)\in
\overline{\M}_{g,k}$, there is an associated element of
$Comb(g,k)$ as follows: we take the graph $T=T_\Sigma$, and let
$g_\nu$ be the genus of $\Sigma_\nu$. The set of combinatorial
types $Comb(g,k)$ is known to be of finite order ([FO]).

There is a partial order $\succ$ on $Comb(g,k)$ defined as
follows. Let $(T,(g_\nu))\in Comb(g,k)$. We consider
$(T_\nu,(g_{\nu\omega})) \in Comb(g_\nu,k_\nu)$ for some of the
vertices $\nu=\nu_1,\cdots,\nu_a$ of $T$. We replace the vertex
$\nu$ of $T$ by the graph $T_\nu$, and join the edge containing
$\nu$ to the vertex $o_\nu(i)$ where $i\in \{1,\cdots,k_\nu\}$ is
the suffix corresponding to this edge. We then obtain a new graph
$\tilde{T}$. The number $\tilde{g}_\nu$ is determined from $g_\nu$
and $g_{\nu\omega}$ in an obvious way.  It is easily seen that
$(\tilde{T},(\tilde{g}_\nu), \tilde{o})$ is in $Comb(g,k)$. We
then define  $(T,(g_\nu),o)\succ
(\tilde{T},(\tilde{g}_\nu),\tilde{o})$.

The structure of the Deligne-Mumford compactification $\overline{\M}_{g,k}$
is summarized in the following

    \vskip 0.1in
\noindent{\bf Proposition 2.3.6: }{\it ([FO]) Let
$\M_{g,k}(T,(g_\nu))$ be the set of all stable curves such that
the associated object is $(T,(g_\nu))$. Then
\begin{itemize}
\item $\overline{\M}_{g,k}$ is a compact complex orbifold which admits a
stratification with finitely many strata, and each stratum is of
the form $\M_{g,k}(T,(g_\nu))$.
\item There is a fiber bundle $\U_{g,k}(T,(g_\nu))\rightarrow
\M_{g,k}(T,(g_\nu))$ which has the following property. For each
$x=(\Sigma_x,\z_x)\in \M_{g,k}(T,(g_\nu))$, there is a
neighborhood of $x$ in $\M_{g,k}(T,(g_\nu))$ of the form
$U_x=V_x/G_x$, where $G_x=Aut(\Sigma_x,\z_x)$,  such that the
inverse image of $U_x$ in $\U_{g,k}(T,(g_\nu))$ is diffeomorphic
to $V_x\times\Sigma_x/G_x$. There is a complex structure on each
fiber such that the fiber of $y=(\Sigma_y,\z_y)$ is identified
with $(\Sigma_y,\z_y)$ itself, together with a K\"{a}hler metric
$\mu_y$ which is flat in a neighborhood of the singular points and
varies smoothly in $y$.
\item $\M_{g,k}(T^\prime,(g_\nu^\prime))$ is contained in the
compactification of $\M_{g,k}(T,(g_\nu))$ in $\overline{\M}_{g,k}$
only if $(T,(g_\nu))\succ (T^\prime,(g_\nu^\prime))$.
\item Different strata are patched together in a way which is described in
the following local model of a neighborhood of a stable curve in
$\overline{\M}_{g,k}$. A neighborhood of $x=(\Sigma,\z)$ in
$\overline{\M}_{g,k}$ is parametrized by $$ \frac{V_x\times
B_r(\oplus_z T_{z_\nu}\Sigma_\nu\otimes
T_{z_\omega}\Sigma_{\omega})} {Aut(\Sigma,\z)}, $$ where
$z=\pi_\nu(z_\nu)=\pi_\omega(z_\omega)$ (Here it may happen that
$\nu=\omega$) runs over all singular points of $\Sigma$, and
$B_r(W)$ denotes the ball of radius $r$ of the vector space $W$.
Each $y\in V_x$ represents a stable curve $(\Sigma_y,\z_y)$
homeomorphic to $(\Sigma,\z)$, with a K\"{a}hler metric $\mu_y$
which is flat in a neighborhood of the singular points. Given
$y\in V_x$, for each element $\varsigma= (\sigma_z)\in\oplus_z
T_{z_\nu}\Sigma_\nu\otimes T_{z_\omega}\Sigma_\omega$ there is an
associated stable curve $(\Sigma_{y,\varsigma},\z_{y,\varsigma})$
obtained as follows. Each component $\Sigma_{\nu}$ of $\Sigma_y$
is given a K\"{a}hler metric $\mu_y$ which is flat in a
neighborhood of the singular points. This gives a Hermitian metric
on each $T_{z_\nu}\Sigma_{\nu}$. For each non-zero $\sigma_z\in
T_{z_\nu}\Sigma_{\nu}\otimes T_{z_\omega}\Sigma_{\omega}$, there
is a biholomorphic map $\Psi_{\sigma_z}: T_{z_\nu}\Sigma_{\nu}
\setminus\{0\} \rightarrow
T_{z_\omega}\Sigma_{\omega}\setminus\{0\}$ defined by
$u\otimes\Psi_{\sigma_z}(u)=\sigma_z$. Let $|\sigma_z|=R^{-2}$;
then for sufficiently large $R$, the map
$exp_{z_\omega}^{-1}\circ\Psi_{\sigma_z} \circ exp_{z_\nu}$ is a
biholomorphism between $D_{z_\nu}(R^{-1/2})\setminus
D_{z_\nu}(R^{-3/2})$ and $D_{z_\omega}(R^{-1/2})\setminus
D_{z_\omega}(R^{-3/2})$, where $D_{z_\nu}(R^{-1/2})$ is a disc
neighborhood of $z_\nu$ in $\Sigma_{\nu}$ of radius $(R^{-1/2})$
which is flat assuming $R$ is sufficiently large. We glue
$\Sigma_\nu$ and $\Sigma_\omega$ by this biholomorphism. If
$\sigma_z=0$, we do not make any change. Thus we obtain
$(\Sigma_{y,\varsigma},\z_{y,\varsigma})$. Moreover, there is a
K\"{a}hler metric $\mu_{y,\varsigma}$ on $\Sigma_{y,\varsigma}$
which coincides with the K\"{a}hler metric $\mu_y$  on $\Sigma_y$
outside a neighborhood of the singular points, and varies smoothly
in $\varsigma$. Each $\gamma\in Aut(\Sigma,\z)$ takes
$(\Sigma_y,\z_y)$ to $(\Sigma_{\gamma(y)},\z_{\gamma(y)})$
isometrically, so it acts on $\oplus_z T_{z_\nu}\Sigma_\nu\otimes
T_{z_\omega}\Sigma_\omega$. $\gamma$ induces an  isomorphism
between $(\Sigma_{y,\varsigma},\z_{y,\varsigma})$ and
$(\Sigma_{\gamma(y,\varsigma)},\z_{\gamma(y,\varsigma)})$, which
is also an isometry.
\end{itemize}}\vskip 0.1in

Now we define a topology on the moduli space $\overline{\M}_{g,k}(X,J,A,\x)$.
We put a Hermitian metric $h$ on $(X,J)$ and the  distance
function on $X$ is assumed to be induced from $h$.

\vspace{2mm}

\noindent{\bf Definition 2.3.7: }{\it A sequence of equivalence
classes of stable maps $x_n$ in $\overline{\M}_{g,k}(X,J,A,\x)$ is
said to converge to $x_0\in\overline{\M}_{g,k}(X,J,A,E)$ if there
are representatives $(f_n,(\Sigma_n,\z_n),\xi_n)$ of $x_n$ and a
representative $(f_0,(\Sigma_0,\z_0),\xi_0)$ of $x_0$ for which
the following conditions hold.
\begin{itemize}
\item For each $n$ (including $n=0$),
there is a set of distinct {\it regular} points
$\{z_{n,1},\cdots,z_{n,a}\}$ (it may happen that this set is
empty) on $\Sigma_n$ which is disjoint from the marked point set
$\z_n$ such that after adding this set to $\z_n$, we obtain a
stable curve in $\overline{\M}_{g,k+a}$, denoted by
$(\Sigma_n,\z_n)^+$. Let $(f_n^+,(\Sigma_n,\z_n)^+,\xi_n^+)$ be
the sequence of stable maps naturally obtained.
\item The sequence $(\Sigma_n,\z_n)^+$ converges to
$(\Sigma_0,\z_0)^+$ in $\overline{\M}_{g,k+a}$. This means that
for sufficiently large $n$, $(\Sigma_n,\z_n)^+$ is identified with
$(\Sigma_{y_n,\varsigma_n},\z_{y_n,\varsigma_n})$ for some
$(y_n,\varsigma_n)$ in the canonical model of a neighborhood of
$(\Sigma_0,\z_0)^+$. Let $\varsigma_n$ be given by
$(\sigma_{z,n})$ and $|\sigma_{z,n}|=R_{z,n}^{-2}$ (here $R_{z,n}$
is allowed to be $\infty$), For each $\mu>\max_z(R_{z,n}^{-1})$ we
put $$ W_{z,n}(\mu)=(D_{z_\nu}(\mu)\setminus
D_{z_\nu}(R_{z,n}^{-1}))\cup (D_{z_\omega}(\mu)\setminus
D_{z_\omega}(R_{z,n}^{-1})),\hspace{2mm}\mbox{and} \hspace{2mm}
W_n(\mu)=\cup_z W_{z,n}(\mu). $$ Then the following holds. First,
for each $\mu>0$, when $n$ is sufficiently large the restriction
of $\tilde{f}^+_n$ to $\Sigma_{y_n,\varsigma_n}\setminus W_n(\mu)$
converges to $\tilde{f}^+_0$ in the $C^\infty$ topology as a
$C^\infty$ map with an isomorphism class of compatible systems.
Secondly, $\lim_{\mu\rightarrow 0}\limsup_{n\rightarrow\infty}
Diam(f_n(W_{z,n}(\mu)))=0$ for each singular point $z$ of
$\Sigma_0$.\\
\item Suppose that $z_0\in \Sigma_0$ is a nodal point near where
$(f_0,\xi_0)$ is given by an equivariant map $(\tilde{f}_0,\lambda_0):
(X(0,r_0),\Z_m)\rightarrow (V_p,G_p)$, where
$p=f(z_0)\in X$ is an orbifold point with a uniformizing system
$(V_p, G_p, \pi_p)$. Then there exist $t_n\mapsto 0$ such that
$(f_n,\xi_n)$ can be represented locally by equivariant
maps $(\tilde{f}_n,\lambda_0):(X(t_n,r_0),\Z_m)\rightarrow (V_p,G_p)$.
\end{itemize}}\vskip 0.1in

 \noindent{\bf Proposition 2.3.8: }{\it Suppose $X$ is either
a symplectic orbifold with a symplectic form $\omega$ and an
$\omega$-compatible almost complex structure $J$, or a projective
orbifold with an integrable almost complex structure $J$. Then the
moduli space $\overline{\M}_{g,k}(X,J,A,\x)$ is compact and
metrizable.} \vskip 0.1in

 The proof follows the same lines of the proof for regular stable maps. We will
 focus on the difference and the rest is referred to [FO]. We first consider the
case when $X$ is a compact symplectic orbifold,
with a symplectic form $\omega$
and an $\omega$-compatible almost complex structure $J$. Then
$h_\omega(u,v)=\frac{1}{2}(\omega(u,Jv)+\omega(v,Ju))$ defines
a Hermitian metric on $(X,J)$. Suppose
$(f,(\Sigma,\z),\xi)$ is a stable map. The area of
$f(\Sigma)$ in $(X,J,h_\omega)$ is given by
$$
Area_{h_\omega}(f(\Sigma))=\sum_\nu\int_{\Sigma_\nu}\tilde{f}_\nu^\ast\omega.
$$

    \vskip 0.1in
\noindent{\bf Lemma 2.3.9: }{\it Let $(f,(\Sigma,\z),\xi)$ be a
stable map in $\overline{\M}_{g,k}(X,J,A,\x)$. The area of
$f(\Sigma)$, $Area_{h_\omega}(f(\Sigma))$, is equal to a constant
depending only on the homology class $A$ and the symplectic form
$\omega$.} \vskip 0.1in

\noindent{\bf Proof:} Let $\tilde{f}:(\Sigma,\z,\m)\rightarrow X$
be any good $C^\infty$ map from a complex orbicurve into $X$ which
induces a pseudo-holomorphic map $f$. Then using the pull-back of
the tangent bundle $TX$ and the exponential map $Exp$, we obtain a
deformation family of $C^\infty$ maps $\tilde{f}_t$ for $0\leq t
\leq \epsilon$ such that for each $0<t\leq \epsilon$,
${f}_t^{-1}(\Sigma X)$ is a set of finitely many points, and
$\tilde{f}_0=\tilde{f}$. Moreover, $\lim_{t\rightarrow
0}\int_\Sigma \tilde{f}_t^\ast \omega= \int_\Sigma \tilde{f}^\ast
\omega$. Each $(f_t)_\ast([\Sigma])$ defines a homology class in
the $2$nd intersection homology of $X$ with top perversity,
$IH^{\bar{t}}_2(X)$. On the other hand, the restriction of
$\omega$ to $X_{reg}$, $\omega_{reg}$, defines a functional
$<\cdot,[\omega]>$ on $IH^{\bar{t}}_2(X)$, which depends only on
the cohomology class of $\omega$, and
$\int_\Sigma\tilde{f}_t^\ast\omega=
<(f_t)_\ast([\Sigma]),[\omega]>$. Since as an oriented orbifold
$X$ is normal, i.e., $H_{2n}(X,X\setminus\{p\})=\Z$ (here $2n=\dim
X$) for any $p\in X$, we have $IH^{\bar{t}}_i(X)=H_i(X)$. So
$<\cdot,[\omega]>$ is a functional on $H_2(X)$, and
$\int_\Sigma\tilde{f}^\ast\omega= <f_\ast([\Sigma]),[\omega]>$
holds there. The lemma follows easily. \hfill $\Box$

In light of Lemma 2.3.9, we have

    \vskip 0.1in
\noindent{\bf Corollary 2.3.10: }{\it  There exists a constant
$\delta>0$ such that for any non-constant pseudo-holomorphic map
$f:\Sigma\rightarrow X$, the $h_\omega$-area of $f(\Sigma)$ is
greater than $\delta$. Here $\Sigma$ is a closed Riemann surface.}
\vskip 0.1in

The following lemma is an extended version of Lemma 11.2 in [FO]
which was used in the proof of the compactness of moduli spaces of
stable maps.

    \vskip 0.1in
\noindent{\bf Lemma 2.3.11: }{\it (cf. Lemma 11.2 in [FO]) There
exists $\epsilon$ independent of $L$ and depending only on
$(X,\omega,J)$ with the following properties. If $f: [-L,L]\times
S^1\rightarrow X$ is a pseudo-holomorphic map which admits a
$C^\infty$ lifting such that there are no orbifold points on
$[-L,L]\times S^1$, and if $Diam (f([-L,L]\times S^1))<\epsilon$,
then $$ |\frac{\partial f}{\partial \tau}(\tau,t)|+|\frac{\partial
f}{\partial t} (\tau,t)| \leq
Ce^{-\frac{1}{m}dist(\tau,\partial[-L,L])} \hspace{3mm} \forall
(\tau,t)\in [-L+1,L-1]\times S^1. $$ Here $C$ and $m$ are
independent of $L$.} \vskip 0.1in

\noindent{\bf Proof:}
When the diameter of $f([-L,L]\times S^1)$ is small enough, it lies in a
uniformized open subset $U$ of $X$, uniformized by $(V,G,\pi)$. Then there
is a finite covering map $co: [-L/m,L/m]\times S^1\rightarrow
[-L,L]\times S^1$ and a lifting $\tilde{f}:[-L/m,L/m]\times S^1\rightarrow V$
such that $\pi\circ \tilde{f}=f\circ co$. We can apply Lemma 11.2 in [FO]
to $\tilde{f}$ to get the estimate. Note that the order $m$ of the covering
map $co$ is bounded from above by a constant depending only on $X$.
\hfill $\Box$

Given a sequence of stable maps $(f_n,(\Sigma_n,\z_n),\xi_n)$, by
adding finitely many additional marked points (whose number is
independent of $n$), we may assume that $(\Sigma_n,\z_n)$ are
stable curves. Then a subsequence of $(\Sigma_n,\z_n)$ converges
in the Deligne-Mumford compactification, so that there is a
$(\Sigma_0,\z_0)$, for large $n$, such that $(\Sigma_n,\z_n)$ is
identified with $(\Sigma_{y_n,\varsigma_n},\z_{y_n,\varsigma_n})$
for a sequence $(y_n,\varsigma_n)$ in the canonical model of a
neighborhood of $(\Sigma_0,\z_0)$. Each
$(\Sigma_{y_n,\varsigma_n},\z_{y_n,\varsigma_n})$ comes with a
K\"{a}hler metric, and we can alter it in a small neighborhood of
the marked points so that it becomes a K\"{a}hler metric of the
corresponding orbicurve defined by $\xi_n$. Using this K\"{a}hler
metric and the Hermitian metric $h_\omega$ on $(X,J)$, we can talk
about the norm of the gradient of $\tilde{f}_n$, denoted by
$|d\tilde{f}_n|$.

In [FO], the bubbling phenomena is naturally interpreted as a result of
degeneration of the domain of the map, which can be easily generalized
to the present case, and is achieved in the following lemma in which
$\mu>0$ is any fixed  sufficiently small number.

    \vskip 0.1in
\noindent{\bf Lemma 2.3.12: }{\it  (cf. Proposition 11.3 in [FO])
By increasing the number of marked points and by taking a
subsequence if necessary, we may assume that $$
\sup_{\Sigma_{y_n,\varsigma_n}\setminus W_n(\mu)} |d\tilde{f}_n|<C
$$ where $C$ is independent of $n$.} \vskip 0.1in

    \vskip 0.1in
\noindent{\bf Lemma 2.3.13: }{\it  Let $(\Sigma,\z,\m)$ be a
compact complex orbicurve (with or without boundary) with a
compact family of complex structures $j_n$ and K\"{a}hler metrics
$\mu_n$. Let $\tilde{f}_n$ be a sequence of $C^\infty$ maps from
$(\Sigma,\z,\m)$ to $(X,J,\omega)$ which induces a
pseudo-holomorphic map $f_n$ with respect to $j_n$, and measured
with $\mu_n$, the norm of the gradient of $\tilde{f}_n$,
$|d\tilde{f}_n|$, is uniformly bounded by a constant $C$. Then
there is a subsequence of $\tilde{f}_n$, denoted by
$\tilde{f}_{n_i}$, such that the corresponding complex structure
$j_{n_i}$ converges to a complex structure $j_0$, and
$\tilde{f}_{n_i}$ converges in $C^\infty$ topology to a $C^\infty$
map $\tilde{f}_0$ which induces a pseudo-holomorphic map $f_0$
with respect to the complex structure $j_0$.} \vskip 0.1in

\noindent{\bf Proof:}
For any $z\in \Sigma$, suppose a subsequence $f_{n_i}(z)$ converges to a point
$p\in X$. Let $(V_p,G_p,\pi_p)$ be a local chart at $p$ with $U_p=\pi_p(V_p)$.
Then since $|d\tilde{f}_n|$ is bounded by $C$, there exists a disc
neighborhood $D$ such that for large $n_i$, $f_{n_i}(D)$ lies in $U_p$.
The issue here is to construct local liftings of $f_{n_i}$ from
the uniformizing system of a fixed neighborhood of $z$ into $V_p$.
We may assume that $D\setminus \{z\}$ contains no orbifold points. For each
$f_{n_i}$, there is a disc neighborhood $D_i$ of $z$, a branched covering map
$br$ (may be trivial) and a local lifting $\tilde{f}_i$ into $V_p$ such that
$\pi_p\circ \tilde{f}_i=f_{n_i}\circ br$, and for any other point in $D$,
$f_{n_i}$ is lifted to a map into $V_p$. It is easily seen that these local
liftings can be patched together to define a lifting $\tilde{f}_{n_i,z}$ on
a branched cover of $D$ into $V_p$. The lemma is then reduced to the classical
case.
\hfill $\Box$

    \vskip 0.1in
\noindent{\bf Lemma 2.3.14: }{\it  Let $(\Sigma,\z,\m)$ be a
compact complex orbicurve (with or without boundary) with a family
of complex structures $j_n$. Let $(\tilde{f}_n,\xi_n)$ be a
sequence of good $C^\infty$ maps from $(\Sigma,\z,\m)$ to
$(X,J,\omega)$, with isomorphism class of compatible systems
$\xi_n$, which induces a pseudo-holomorphic map $f_n$ with respect
to $j_n$. Assume $j_n$ converges to a complex structure $j_0$, and
$\tilde{f}_n$ converges in the $C^\infty$ topology to a $C^\infty$
map $\tilde{f}_0$ which induces a pseudo-holomorphic map $f_0$
with respect to $j_0$. Then $\tilde{f}_0$ is good and there is an
isomorphism class of compatible systems $\xi_0$ of $\tilde{f}_0$
such that a subsequence of $(\tilde{f}_n,\xi_n)$ converges in the
$C^\infty$ topology to $(\tilde{f}_0,\xi_0)$. Moreover, for
sufficiently large $n$, there exists a $C^\infty$ section
$\tilde{s}_n$ of the pull-back orbifold bundle $(TX)^\ast_{\xi_0}$
over $(\Sigma,\z,\m)$ determined by $\xi_0$, such that
$\tilde{f}_n$ equals $(\tilde{f}_0)_{\xi_0,s_n}$ and $\xi_n$
equals the canonically determined isomorphism class of compatible
systems of $(\tilde{f}_0)_{\xi_0,s_n}$. As a consequence,
$(\tilde{f}_0,\xi_0)$ is the unique limit of the said subsequence
of $(\tilde{f}_n,\xi_n)$.} \vskip 0.1in

\noindent{\bf Proof:} First of all, as in the construction of
induced geodesic compatible systems, we can construct a compatible
cover $\U$ of $(\Sigma,\z,\m)$ consisting of geodesic discs, a set
of geodesic neighborhoods $\U^\prime_0$ in $X$, and a collection
of local liftings $\tilde{f}_{0,UU^\prime_0}$ of $f_0$ such that
the elements of $\U$ are in 1:1 correspondence with those of
$\U^\prime_0$, and if $U_2\subset U_1$, there is an injection
between the uniformizing systems of $U^\prime_{2,0}$ and
$U^\prime_{1,0}$. Likewise, there is a set of uniformized open
subsets $\U^\prime_n$ for each $n$, and a local lifting
$\tilde{f}_{n,UU^\prime_n}$ of $f_n$ such that for each $U\in\U$,
there is $n(U)>0$, when $n>n(U)$, there is an injection
$\delta_{U_n^\prime}$ from the uniformizing systems of
$U^\prime_n$ into that of $U^\prime_0$. There is a 1:1
correspondence between elements of $\U$ and $\U^\prime_n$ such
that if $U_2\subset U_1$, there is an injection between the
uniformizing systems of $U^\prime_{2,n}$ and $U^\prime_{1,n}$. We
want to point out that there is a delicate point here. For any
collection of finitely many $C^\infty$ maps, we can construct the
sets $\U$ and $\U^\prime_n$ with $\U$ common and $\U^\prime_n$
consisting of geodesically convex and star-shaped open subsets
without any problem, but for infinitely many $C^\infty$ maps we
can not ensure  that $\U^\prime_n$ consists of geodesically convex
and star-shaped open subsets while $\U$ is common to all the maps.
This is because there is no lower bound of the injectivity radius
of points in $X$, although $X$ is compact. In the present case, we
have to exploit the fact that $\tilde{f}_n$ converges to
$\tilde{f_0}$ and $f_n$ is pseudo-holomorphic. For each point $z$
in $\Sigma$, there are two possibilities: the first one is that
there is a disc neighborhood $D$ of $z$ such that, for large $n$
$f(z)$ has the largest isotropy subgroup amongst the points in
$f_n(D)$; the second one is that there is a disc neighborhood $D$
of $z$ such that if let $z_n\in D$ be a point such that $f_n(z_n)$
has the largest isotropy subgroup amongst the points in $f_n(D)$,
then $z_n\neq z$ and $z_n\rightarrow z$ as $n\rightarrow \infty$.
In the first case, we can pick a geodesically convex and
star-shaped open subset centered at $f_n(z)$ (not necessarily a
ball) for the preliminary candidate of $U^\prime_n$, while in the
second case, we pick a geodesically convex and star-shaped open
subset centered at $f_n(z_k)$. One can verify that this will work
out.

Now we use the given isomorphism class of compatible systems
$\xi_n$ to define, for each injection $i$ from the uniformizing
system of $U_2$ into that of $U_1$, an injection $\lambda_n(i)$
from the uniformizing system of $U^\prime_{2,n}$ to that of
$U^\prime_{1,n}$ such that
$\{\tilde{f}_{n,UU^\prime_n},\lambda_n\}$ becomes a compatible
system of $\tilde{f}_n$ within the isomorphism class of $\xi_n$.
We first pick a compatible system
$\{\tilde{f}_{n,W_nW^\prime_n},\tau_n\}$ in $\xi_n$. For each
$U\in\U$ such that $U\subset \Sigma\setminus\z$, we pick a $W_n\in
\W_n$ such that $W_n\subset U$ and there is an injection
$\delta_{W_nU}$ from the uniformizing system of $W_n^\prime$ into
that of $U_n^\prime$ satisfying
$\delta_{W_nU}\circ\tilde{f}_{n,W_nW^\prime_n}=
\tilde{f}_{n,UU^\prime_n}|_{W_n}$. We will fix $\delta_{W_nU}$.
Then this can be extend out to determine an injection
$\delta_{W_{1,n}U}$ for any other element $W_{1,n}\in\W_n$ such
that $W_{1,n}\subset U$ whenever it is possible. Now for any
$U_1$, $U_2\in \U$ with an inclusion $i:U_2\subset U_1$, we take a
$W_n\in\W_n$, $W_n\subset U_2\subset U_1$, and we define
$\lambda_n(i)$ to be the unique injection such that
$\delta_{W_nU_2}=\lambda_n(i)\circ \delta_{W_nU_1}$. One can check
that $\{\tilde{f}_{n,UU^\prime_n},\lambda_n\}$ is a compatible
system of $\tilde{f}_n$ on $\Sigma\setminus\z$, which extends over
$\Sigma$ uniquely, and is in the isomorphism class of $\xi_n$.

Now we take a sequence of subsets $\U_k$ of $\U$ such that
$\U_k\subset \U_{k+1}$, $\U=\cup_k\U_k$, and each $\U_k$ is a finite cover
of $\Sigma$. By the finiteness of each $\U_k$, there is a subsequence
of $\{\tilde{f}_{n,UU^\prime_n},\lambda_n\}$ such that for any injection
$i$ associated to an inclusion $U_2\subset U_1$, $U_1,U_2\in\U_k$, we
have an injection  $\lambda_0(i)$ such that $\delta_{U_{1,n}^\prime}\circ
\lambda_n(i)=\lambda_0(i)\circ\delta_{U_{2,n}^\prime}$
and $\delta_{U^\prime_n}\circ \tilde{f}_{n,UU^\prime_n}$
converges to $\tilde{f}_{0,UU^\prime_0}$. Then it follows that
the diagonal subsequence converges to
$\{\tilde{f}_{0,UU^\prime_0},\lambda_0\}$.

Finally, we take a finite subset $\U_0$ of $\U$ such that $\U_0$
covers $\Sigma$ and for any $U_1,U_2\in\U_0$, there is a
$U_3\in\U_0$ such that $U_3\subset U_1\cap U_2$. Then for each
$n>N=\max_{U\in\U_0} n(U)$, the image of $\delta_{U^\prime_n}\circ
\tilde{f}_{n,UU^\prime_n}$ lies in the uniformizing system of
$U^\prime_0$ for any $U\in\U_0$. This gives rise to a collection
of local $C^\infty$ sections $\tilde{s}_U$ over $U$ of the
pull-back orbifold bundle $(TX)^\ast_{\xi_0}$ defined by
$\{\tilde{f}_{0,UU^\prime_0},\lambda_0\}$. On the intersection
$U_1\cap U_2$, $\tilde{s}_{U_1}=\tilde{s}_{U_2}$ restricted to
$U_3$. By the unique continuity property of pseudo-holomorphic
maps, $\tilde{s}_{U_1}=\tilde{s}_{U_2}$ on $U_1\cap U_2$, so that
they patch together to a global section $\tilde{s}_n$, and
$\tilde{f}_n$ equals $(\tilde{f}_0)_{\xi_0,s_n}$. The fact that
$\xi_n$ equals the canonically determined isomorphism class of
compatible systems of $(\tilde{f}_0)_{\xi_0,s_n}$ follows from the
description of the difference of any two isomorphism classes of
compatible systems of $\tilde{f}_n$ in terms of a homomorphism
from $\pi_1(\Sigma\setminus\z,z_0)$  to $G_{f_n(z_0)}$.

\hfill $\Box$.

Now for a closed symplectic orbifold $(X,\omega)$ with a
$\omega$-compatible almost complex structure $J$, the argument in
[FO] can be taken word by word to show that the moduli space
$\overline{\M}_{g,k}(X,J,A,\x)$ is compact, except that we have to
make sure that a limiting compatible system is determined at each
nodal point of the limiting curve and the convergence is as
specified in Definition 2.3.7. But this follows from Lemma 2.2.4,
Lemma 2.2.6 (see the remark after the proof of Lemma 2.2.6)
easily.

Now let's consider the case when $(X,J)$ is a projective orbifold
with a Hermitian metric $h$ on it. By  assumption, $X\subset \P^N$
is a subvariety of a projective space $\P^N$. A continuous map
$f:\Sigma\rightarrow X$ is analytic if and only if
$f:\Sigma\rightarrow \P^N$ is holomorphic. Now suppose we have a
sequence of stable maps $(f_n,(\Sigma_n,\z_n),\xi_n)$ in
$\overline{\M}_{g,k}(X,J,A,\x)$; forgetting the twisted boundary
conditions $\xi_n$, we have a sequence of stable maps
$(f_n,(\Sigma_n,\z_n))$ in $\overline{\M}_{g,k}(\P^N,J_0,A)$, and
by the classical Gromov compactness theorem, there is a
subsequence converging to a stable map $(f_0,(\Sigma_0,\z_0))$ in
$\overline{\M}_{g,k}(\P^N,J_0,A)$, which is measured by the
Fubini-Study metric on $\P^N$. Let's still use
$(f_n,(\Sigma_n,\z_n))$ to denote the subsequence, and assume that
additional marked points have been added. Then for sufficiently
large $n$, there are points $(y_n,\varsigma_n)$ in the canonical
model of a neighborhood of $(\Sigma_0,\z_0)$ in the
Deligne-Mumford compactification such that $(\Sigma_n,\z_n)$ are
identified with $(\Sigma_{y_n,\varsigma_n},
\z_{y_n,\varsigma_n})$, and for any small $\mu>0$, the restriction
of $f_n$ to $\Sigma_{y_n,\varsigma_n}\setminus W_n(\mu)$ converges
to $f_0$ in $C^\infty$, and $\lim_{\mu\rightarrow
0}\sup\lim_{n\rightarrow\infty} Diam(W_n(\mu))=0$. Here the
$C^\infty$ norm and the distance function are defined from the
Fubini-Study metric. We observe:

\begin{itemize}
\item Lemma 2.3.11 (cf. Lemma 11.2 in [FO]) holds for general Hermitian
orbifolds.
\item The convergence of $f_n$ to $f_0$ on
$\Sigma_{y_n,\varsigma_n}\setminus W_n(\mu)$ (although measured by
the Fubini-Study metric) implies that $$
\sup_{\Sigma_{y_n,\varsigma_n}\setminus W_n(\mu)}
|d\tilde{f}_n|_h<C $$ where $|\hspace{2mm}|_h$ is the norm defined
by the Hermitian metric $h$ on $X$, and $C$ is independent of $n$,
for otherwise, a bubble would have  developed, which contradicts
the convergence of $f_n$ to $f_0$ in $\P^N$.
\item $\lim_{\mu\rightarrow 0}\sup\lim_{n\rightarrow\infty}
Diam_h(W_n(\mu))=0$ holds where $Diam_h$ is the diameter defined
by $h$, since otherwise a bubble would have  developed.
\end{itemize}

With these understood, a similar argument shows that
the moduli space $\overline{\M}_{g,k}(X,J,A,\x)$ is compact for a projective
orbifold $(X,J)$, measured by any fixed Hermitian metric $h$.

Finally, we will show that the moduli space $\overline{\M}_{g,k}(X,J,A,\x)$
is metrizable. We recall the Metrization Theorem of Urysohn (see [Ke]):

{\it A topological space $X$ is metrizable if the following
conditions hold:
\begin{itemize}
\item $\{x\}$ is closed for each point $x\in X$, i.e., $X$ is a $T_1$-space.
\item For each $x\in X$ and each neighborhood $U$ of $x$, there is a closed
neighborhood $V$ of $x$ such that $V\subset U$ (i.e. $X$ is regular).
\item $X$ has a countable base.
\end{itemize}
}

First, the moduli space $\overline{\M}_{g,k}(X,J,A,\x)$ is Hausdorff, as shown
in [FO], so it is a $T_1$-space.

Secondly, for each point $x_0$ in $\overline{\M}_{g,k}(X,J,A,\x)$,
represented by a stable map $(f_0,(\Sigma_0,\z_0),\xi_0)$, we will
define a pseudo-distance function $d$ as follows. For any point
$x=[(f,(\Sigma,\z),\xi)]$ in $\overline{\M}_{g,k}(X,J,A,\x)$, let
$d_1(x_0,x)=\min (dist((\Sigma_0,\z_0)^+,(\Sigma,\z)^+))$ where
$(\Sigma,\z)^+$ denotes any stabilization of $(\Sigma,\z)$, and
$dist$ is the distance function on the corresponding
Deligne-Mumford compactification; let $d_2(x_0,x)=\min_{f:[f]\in
x}\max_{f^+}
\max_{(z,z_0)\in(\z^+,\z_0^+)}(dist_h(f^+(z),f_0^+(z_0)))$, and
let $d=d_1+d_2$. Then fixing a small $\mu>0$, there is an
$\epsilon(x_0)>0$ such that if $d(x_0,x)\leq \epsilon(x_0)$, there
is a representative $(f,(\Sigma,\z),\xi)$ of $x$ such that the
restriction of $\tilde{f}$ to $\Sigma\setminus W(\mu)$ equals
$(\tilde{f}_0)_{\xi_0^\prime,s}$ for some $\xi^\prime_0$ of $f_0$
and a $C^\infty$ section $\tilde{s}$. Now for any $0<\epsilon\leq
\epsilon(x_0)$, we define the set $U(x_0,\epsilon)$ to be the set
of all $x\in \overline{\M}_{g,k}(X,J,A,\x)$ such that
$d(x_0,x)<\epsilon$ and if a representative $(f,(\Sigma,\z),\xi)$
of $x$  equals $(\tilde{f}_0)_{\xi_0^\prime,s}$ on
$\Sigma\setminus W(\mu)$, then $\xi_0^\prime=\xi_0$. Then
$U(x_0,\epsilon)$ defines a family of open neighborhoods of $x_0$
and $\cap_{\epsilon}U(x_0,\epsilon)=\{x_0\}$. Observe that for any
$\epsilon_1<\epsilon$, $\overline{U(x_0,\epsilon_1)} \subset
U(x_0,\epsilon)$, hence  $\overline{\M}_{g,k}(X,J,A,\x)$ is
regular.

Finally, $\overline{\M}_{g,k}(X,J,A,\x)$ has a countable base, since it
is compact, and at each point $x_0$, there is a family of open neighborhoods
$U(x_0,\epsilon)$ of $x_0$ where $0<\epsilon\leq\epsilon(x_0)$.
Hence the moduli space $\overline{\M}_{g,k}(X,J,A,\x)$ is metrizable by
Urysohn's theorem.

\section{Virtual Cycle and Orbifold Gromov-Witten Invariants}

The orbifold Gromov-Witten invariants are defined by constructing
the virtual fundamental cycle on the moduli spaces of orbifold
stable maps using the technique of Li-Tian \cite{LT} and
Fukaya-Ono \cite{FO}. For this purpose, we first construct a
Kuranishi structure on the moduli space along the lines of [FO],
although our treatment of the analysis is different. By a general
procedure in [FO], the Kuranishi structure already gives rise to a
virtual fundamental class, which suffices to define the GW
invariants and prove its usual axioms.

\subsection{Review of Kuranishi structure}

In this subsection, we give a brief review of the notion of
Kuranishi structure; see [FO] for more details.

Let $X$ be a compact, metrizable topological space.

\noindent{\bf Definition 3.1.1:}{\it
A Kuranishi neighborhood of $p\in X$ is a system $(U_p,E_p,s_p,\psi_p)$ where

\begin{itemize}
\item $U_p=V_p/G_p$ is an orbifold and $E_p$ is an orbifold bundle on it.
\item $s_p$  is a continuous section of $E_p$.
\item $\psi_p$ is a homeomorphism from $s_p^{-1}(0)$ to a neighborhood of $p$
in $X$.
\end{itemize}
We call $E_p$ the {\it obstruction bundle} and $s_p$ the Kuranishi map. }

For any $q\in Im \psi_p$, the Kuranishi neighborhood
$(U_p,E_p,s_p,\psi_p)$ induces a  Kuranishi neighborhood of $q$ by
restricting to a neighborhood of $\psi_p^{-1}(q)$ in $U_p$. Two
Kuranishi neighborhoods of $p$, $(U_p,E_p,s_p,\psi_p)$ and
$(U_p^\prime,E_p^\prime,s_p^\prime,\psi_p^\prime)$, are {\it
isomorphic} if there exists an isomorphism $(I,J):(U_p,E_p)
\rightarrow (U_p^\prime,E_p^\prime)$ such that $J\circ
s_p=s_p^\prime\circ I$ and $\psi_p^\prime\circ I=\psi_p$. A {\it
germ} of Kuranishi neighborhood of $p$ is defined in the following
sense: $(U_p,E_p,s_p,\psi_p)$ and
$(U_p^\prime,E_p^\prime,s_p^\prime,\psi_p^\prime)$ are {\it
equivalent} if they induce isomorphic Kuranishi neighborhoods of
$p$ on a smaller neighborhood of $p$. A germ of Kuranishi
neighborhood is said to be a {\it stabilization} of another if
there are representatives
$(U_p^\prime,E_p^\prime,s_p^\prime,\psi_p^\prime)$ and
$(U_p,E_p,s_p,\psi_p)$ respectively and an embedding of orbifolds
$(\phi,\hat{\phi}):(U_p,E_p) \rightarrow (U_p^\prime,E_p^\prime)$
such that $s_p^\prime\circ\phi= \hat{\phi}\circ s_p$ and
$\psi_p^\prime\circ\phi=\psi_p$, together with an isomorphism
$\Phi:(TU_p^\prime)^\ast/TU_p\rightarrow E_p^\prime/E_p$ of
orbibundles on $U_p$, where $(TU_p^\prime)^\ast/TU_p$ is the
normal bundle of $\phi(U_p)$ in $U_p^\prime$. One can similarly
define the notion of {\it germ} of $(\phi,\hat{\phi},\Phi)$.

\noindent{\bf Definition 3.1.2: }{\it A {\it Kuranishi structure
of dimension $n$} on $X$ assigns a germ of Kuranishi neighborhood
to each point $p\in X$, such that for each representative
$(U_p,E_p,s_p,\psi_p)$ of it, $\dim U_p-rank \ E_p=n$, and the
induced germ of Kuranishi neighborhood at any $q\in Im\psi_p$ is a
stabilization of the germ of Kuranishi neighborhood of $q$.
Moreover, for each sufficiently small representative
$(U_q,E_q,s_q,\psi_q)$, there is a $(\phi_{pq},\hat{\phi}_{pq},
\Phi_{pq}): (U_q,E_q,(TU_p)^\ast/TU_q)\rightarrow
(U_p,E_p,E_p/E_q)$ satisfying the following compatibility
condition: for each $r\in Im \psi_q\cap Im \psi_p $ and any
sufficiently small representative $(U_r,E_r,s_r,\psi_r)$, the
equation $$ (\phi_{pq},\hat{\phi}_{pq})\circ
(\phi_{qr},\hat{\phi}_{qr})= (\phi_{pr},\hat{\phi}_{pr}) $$ holds
and the diagram $$
\begin{array}{ccccccccc}
0 & \rightarrow & (TU_q)^\ast/TU_r & \rightarrow & (TU_p)^\ast/TU_r &
\rightarrow & (TU_p)^\ast/TU_q |_{U_r} & \rightarrow & 0\\
  &             &  \downarrow \Phi_{qr} &         & \downarrow \Phi_{pr} &
   & \downarrow \Phi_{pq} &       &   \\
0 & \rightarrow & E_q/E_r          & \rightarrow & E_p/E_r          &
\rightarrow & E_p/E_q|_{U_r}           & \rightarrow & 0
\end{array}
$$
commutes.}
\vskip 0.1in

\noindent{\bf Remark 3.1.3: }{\it  A germ of triple
$(\phi_{pq},\hat{\phi}_{pq},\Phi_{pq})$ is well determined in the
above definition for any sufficiently two close  points $p$ and
$q$ in $X$, with a suitable compatibility condition satisfied. The
point of imposing the compatibility condition here is that
although a Kuranishi structure is not a system of trivializations
of an orbifold bundle, but defines a certain abstract object, a
sort of quasi-orbifold bundle such that $X$ is the zero set of a
continuous section of it.}

A $K$-theory can be defined over Kuranishi structures. We first give a
definition of an analogue of vector bundles over a space with a Kuranishi
structure.

\noindent{\bf Definition 3.1.4: }{\it Let $X$ be a space with a
Kuranishi structure. A {\it bundle system} on $X$ consists of the
following data:

\begin{itemize}
\item For each point $p\in X$ there exists a germ of a pair of orbifold bundles
$(F_{1,p},F_{2,p})$ on its Kuranishi neighborhood.
\item For two sufficiently  close  points $p$ and $q$ in $X$, for which the
germ of $(\phi_{pq},\hat{\phi}_{pq},\Phi_{pq})$ is defined, there exist a germ
of embeddings $(\Psi_{1,pq},\Psi_{2,pq}):(F_{1,q},F_{2,q})\rightarrow
(F_{1,p},F_{2,p})$ and a germ of isomorphisms of orbibundles
$\Psi_{pq}:F_{1,p}|_{U_q}/F_{1,q}\rightarrow F_{2,p}|_{U_q}/F_{2,q}$.
\item The following compatibility condition holds: for any triple $(r,q,p)$
such that the germs of $(\phi_{qr},\hat{\phi}_{qr},\Phi_{qr})$,
$(\phi_{pr},\hat{\phi}_{pr},\Phi_{pr})$ and
$(\phi_{pq},\hat{\phi}_{pq},\Phi_{pq})$ are defined, the equation
$$ (\Psi_{1,pq},\Psi_{2,pq})\circ (\Psi_{1,qr},\Psi_{2,qr})
=(\Psi_{1,pr},\Psi_{2,pr}) $$ holds and the diagram $$
\begin{array}{ccccccccc}
0 & \rightarrow & F_{1,q}|_{U_r}/F_{1,r} & \rightarrow & F_{1,p}|_{U_r}/F_{1,r} &
\rightarrow & F_{1,p}|_{U_r}/F_{1,q}|_{U_r} & \rightarrow & 0\\
  &             &  \downarrow \Psi_{qr} &         & \downarrow \Psi_{pr} &
   & \downarrow \Psi_{pq} &       &   \\
0 & \rightarrow & F_{2,q}|_{U_r}/F_{2,r} & \rightarrow & F_{2,p}|_{U_r}/F_{2,r} &
\rightarrow & F_{2,p}|_{U_r}/F_{2,q}|_{U_r} & \rightarrow & 0
\end{array}
$$
commutes.
\end{itemize}}

    \vskip 0.1in

One can define Whitney sum, tensor product, etc. of bundle systems
in an obvious way, as well as an isomorphism relation  between
bundle systems.

\noindent{\bf Example 3.1.5:}{\it
It is easily seen that the germ of pairs $(TU_p,E_p)$ defines a bundle
system; we call it the {\it tangent bundle} of $X$, and denote it  by $TX$.}
    \vskip 0.1in

A bundle system is said to be {\it oriented} if $F_{1,p}$, $F_{2,p}$ are
oriented for each $p\in X$ and
$\Psi_{pq}:F_{1,p}|_{U_q}/F_{1,q}\rightarrow F_{2,p}|_{U_q}/F_{2,q}$ is
orientation preserving. It is said to be {\it complex} if
if $F_{1,p}$, $F_{2,p}$ are complex and $\Psi_{1,pq}$, $\Psi_{2,pq}$ and
$\Psi_{pq}$ are complex linear. A bundle system
$((F_{1,p}, F_{2,p}),(\Psi_{1,pq}, \Psi_{2,pq}, \Psi_{pq}))$ is said to be
{\it trivial} if there exist germs of isomorphisms $F_{1,p}\cong F_{2,p}$
which are compatible with
$(\Psi_{1,pq}, \Psi_{2,pq}, \Psi_{pq})$.

Consider the free abelian group generated by the set of all isomorphism
classes of bundle systems and divide it by the relations
\begin{eqnarray*}
& [((F_{1,p}, F_{2,p}),(\Psi_{1,pq}, \Psi_{2,pq}, \Psi_{pq}))\oplus
((F_{1,p}^\prime, F_{2,p}^\prime),(\Psi_{1,pq}^\prime, \Psi_{2,pq}^\prime,
\Psi_{pq}^\prime))]\\
& =[((F_{1,p}, F_{2,p}),(\Psi_{1,pq}, \Psi_{2,pq}, \Psi_{pq}))]+
[((F_{1,p}^\prime, F_{2,p}^\prime),(\Psi_{1,pq}^\prime, \Psi_{2,pq}^\prime,
\Psi_{pq}^\prime))]\\
& [((F_{1,p}, F_{2,p}),(\Psi_{1,pq}, \Psi_{2,pq}, \Psi_{pq}))]=0 \hspace{2mm}
\mbox{if} \hspace{2mm}
((F_{1,p}, F_{2,p}),(\Psi_{1,pq}, \Psi_{2,pq}, \Psi_{pq})) \hspace{2mm}
\mbox{ is trivial}.
\end{eqnarray*}
This group is called the real $K$-group of $X$ with the Kuranishi structure,
and is denoted by $KO(X)$. By using oriented bundle systems and complex
bundle systems, the corresponding $K$-groups $KSO(X)$ and $K(X)$ are defined,
and there is an obvious map
$$
K(X)\rightarrow KSO(X) \rightarrow KO(X).
$$

A Kuranishi structure $X$ is said to be {\it stably orientable} if $[TX]$ is in
the image of $KSO(X)$. It is said to be {\it stably complex} if $[TX]$ is
in the image of $K(X)$. It is proved in [FO] that being stably orientable is
equivalent to being orientable. It is obvious that if $X$ is stably complex,
then $X$ is orientable with a canonical orientation.

\subsection{Construction of Kuranishi neighborhood}

This subsection is devoted to a local construction of Kuranishi
neighborhood for a point in the moduli space of stable maps
$\overline{\M}_{g,k}(X,J,A,\x)$. First we review briefly the basic
idea of the construction so that for those readers who are
unfamiliar with this type of argument, the big picture will not be
buried in a sea of technical details.

The construction of Kuranishi model for a nonlinear Fredholm map
is based on the following simple lemma of Banach. Let $B:
U\rightarrow U$ be a map from a ball of radius $r$ in a Banach
space to itself such that $$ ||B(x)-B(y)||\leq \epsilon ||x-y|| $$
holds for any $x,y\in U$ for an $\epsilon<1$, then the map $B$ has
a unique fixed point $x_0$ in $U$, i.e., $x_0=B(x_0)$. Moreover,
if $B$ depends on a parameter, say $u$, and $$
||B_u(x)-B_u(y)||\leq \epsilon ||x-y|| $$ holds for all $u$ for an
$\epsilon<1$, then  letting $x_u$ be the fixed point of $B_u$, we
have $$
||x_{u_1}-x_{u_2}||\leq\frac{1}{1-\epsilon}||B_{u_1}(x_{u_2})-B_{u_2}(x_{u_2})
||, $$ and when $||B_{u_1}(x)-B_{u_2}(x)||\leq C||u_1-u_2||||x||$
for any $x\in U$ for some constant $C>0$, we have $$
||x_{u_1}-x_{u_2}||<||u_1-u_2|| $$ when we restrict the maps $B_u$
to a smaller ball of radius less than $\frac{1-\epsilon}{C}$.

Now suppose $F:X\rightarrow Y$ is a $C^2$ nonlinear Fredholm map
from a neighborhood of $0$ in a Banach space $X$ into a
neighborhood of $0$ in a Banach space $Y$ such that  $F(0)=0$. Let
$L$ be the differential of $F$ at $0$, which is a linear Fredholm
map from $X$ to $Y$. Let $E$ be a finite dimensional subspace of
$Y$ such that $Y$ is spanned by $E$ and the image of $L$, $Im L$.
Let $V=\{x\in X|Lx\in E\}$; then $V$ is a finite dimensional
subspace in $X$, $\dim V-\dim E=index \ L$, and there exist
subspaces $V^\prime$ and $E^\prime$ of $X$ and $Y$ such that
$X=V\oplus V^\prime$, $Y=E^\prime \oplus E$,
$L:V^\prime\rightarrow E^\prime$ is invertible, and if
$x=v+v^\prime$ for $v\in V$, $v^\prime \in V^\prime$,
$y=e^\prime+e$ for $e^\prime\in E^\prime$, $e\in E$, then
$||v||\leq c||x||$, $||v^\prime||\leq c||x||$ and
$||e^\prime||\leq c||y||$, $||e||\leq c||y||$ hold for some
constant $c>0$. Write $F(x)=Lx+G(x)$; then $G(0)=0$ and $DG(0)=0$.
(Here $DG$ stands for the differential of $G$.) Let's consider
solving the equation $F(x)\in E$. If we write $x=v+v^\prime$, and
$PG(x)$ for the component of $G(x)$ in $E^\prime$, the equation
$F(x)\in E$ becomes $0=Lv^\prime+PG(v+v^\prime)$. If we let
$B_v(x)=-L^{-1}(PG(v+x))$, then the above equation becomes a fixed
point equation for $B_v$, i.e., $v^\prime=B_v(v^\prime)$. Let's
examine the the norms $||B_v(x_1)-B_v(x_2)||$ and
$||B_{v_1}(x)-B_{v_2}(x)||$.

\begin{eqnarray*}
||B_v(x_1)-B_v(x_2)|| & \leq & ||L^{-1}||||PG(v+x_1)-PG(v+x_2)||\\
                      & \leq & c ||L^{-1}||||G(v+x_1)-G(v+x_2)||\\
                      & \leq & c||L^{-1}||||DG(v+t_0x_1+(1-t_0)x_2)||
                                ||x_1-x_2||,
\end{eqnarray*}
and similarly,
$$
||B_{v_1}(x)-B_{v_2}(x)||\leq c||L^{-1}||||DG(x+t_0v_1+(1-t_0)v_2)||
                                ||v_1-v_2||.
$$

\noindent{\bf Lemma 3.2.1:}{\it
Consider the ball $U_{2r}$ in $X$ of radius $2r$
such that $x\in U_{2r}$ satisfies the condition
$c||L^{-1}||||DG(x)||\leq \frac{1}{3}$. Then for any $v\in V$ such that
$v\in U_r$, there is a unique $v^\prime(v)\in V^\prime\cap U_r$ such that
$F(v+v^\prime(v))\in E$, and $\psi:v\rightarrow v+v^\prime(v)$ is $1:1$.
On the other hand, for any $x\in U_{r/c}$ such that $F(x)\in E$,
there is a unique $v\in V\cap U_r$ such that $x=v+v^\prime(v)$. In particular,
let $f:V\cap\psi^{-1}(U_{r/c})\rightarrow E$ be defined by
$f(v)=F(v+v^\prime(v))$, then
$f$ is continuous and the zero set $f^{-1}(0)$ is homeomorphic to the
zero set $F^{-1}(0)\cap U_{r/c}$ by $v\rightarrow v+v^\prime(v)$.}
    \vskip 0.1in

\noindent{\bf Proof:} The condition
$c||L^{-1}||||DG(x)||\leq\frac{1}{3}$ for all $x\in U_{2r}$
implies that if $v\in V\cap U_r$, $x\in V^\prime\cap U_r$, then
$B_v(x)\in V^\prime\cap U_r$, and the inequality $$
||B_v(x_1)-B_v(x_2)||\leq \frac{1}{3}||x_1-x_2|| $$ holds. So for
any $v\in V\cap U_r$, there is a unique fixed point $v^\prime(v)
\in V^\prime\cap U_r$ of $B_v$, so that $F(v+v^\prime(v))\in E$.
That the map $v\rightarrow v+v^\prime(v)$ is $1:1$ follows from
the inequality $$ ||B_{v_1}(x)-B_{v_2}(x)||\leq
\frac{1}{3}||v_1-v_2|| $$ which implies that
$||v^\prime(v_1)-v^\prime(v_2)||\leq \frac{2}{3} ||v_1-v_2||$.

On the other hand, if $x\in U_{r/c}$ and $F(x)\in E$, and if we
write $x=v+v^\prime$, we have $v,v^\prime\in U_r$, and
$
v^\prime=B_v(v^\prime).
$
Hence $x$ lies in the image of $\psi: v\rightarrow v+v^\prime(v)$.
\hfill $\Box$

\noindent{\bf Remark 3.2.2: }{\it In this construction, the
solution set of $F(x)\in E$ in the ball $U_{r/c}$ of radius $r/c$
is identified with an open set $\psi^{-1}(U_{r/c})\cap V\cap U_r$
in a finite dimensional ball of radius $r$, and the zero set
$F^{-1}(0) \cap U_{r/c}$ is identified with the zero set of a
continuous function $f:\psi^{-1}(U_{r/c})\cap V\cap U_r
\rightarrow E$ between finite dimensional spaces. This is
generally referred to as the Kuranishi model of the map $F$ at
$0$. The size of the Kuranishi model, which is determined by the
number $r$, can be explicitly measured by the constant $c$, the
norm $||L^{-1}||$ and the second derivative $D^2F(0)$ of $F$ at
$0$, while the constant $c$ depends solely on the relative
positions of the subspace $E$ and $Im \ L$ in $Y$. (We will call
the constant $c$ the {\it ratio} of the corresponding
decompositions of Banach spaces $X$ and $Y$.) This knowledge is
important in applications because very often we need to carry out
a parametrized version of this construction, and we need to know
how the size of the Kuranishi model is changing with the
parameter.}
    \vskip 0.1in

In the construction of Kuranishi neighborhood for a stable map, we encounter
a family of nonlinear Fredholm maps $F_s:X_s\rightarrow Y_s$ parametrized by
a finite dimensional ball $s\in B$, such that $F_0(0)=0$ and $F_s(0)$ is
small. There are no explicit identifications between the $X_s$'s and the
$Y_s$'s. We will find a family of isomorphisms
$\theta_s: E_0\rightarrow E_s$ and
$\eta_s:V_0\rightarrow V_s$ such that
$h_s=(\theta_s)^{-1}\circ f_s \circ \eta_s$ is continuous in $s$ so that
we can construct a Kuranishi neighborhood $f: B\times V_0 \rightarrow E_0$
by $(s,v)\rightarrow h_s(v)$. We need to watch out for the dependence of
$c_s$, $||L^{-1}_s||$, $||D^2F_s(0)||$ and $||\theta_s||$, $||\eta_s||$
on the parameter $s$. Moreover, there are finite groups $G$ and $H$ acting
on the parameter space $B$ and the vector space $V_0$ respectively.
\hfill $\Box$

\vspace{2mm}

Now we start out with the construction. Consider the
biholomorphism $\R\times S^1\rightarrow \C^\ast\setminus \{0\}$
given by $t+is\rightarrow e^{-(t+is)}$ where $s\in S^1=\R/2\pi\Z$.
This identifies the punctured disc of radius $r$ with the
half-cylinder $[-\ln r,\infty)\times S^1$. Let
$f:(\Sigma,\z)\rightarrow (X,J,h)$ be a pseudo-holomorphic map
with a twisted boundary condition $\xi$ from a marked Riemann
surface into a Hermitian orbifold (cf. Remark 2.2.10). Let $z_0\in
\z$ be one of the marked points. We intend to think of $f$ as a
map $f_0$ from the punctured Riemann surface
$(\Sigma,\z)\setminus\{z_0\}$ into $X$, where
$(\Sigma,\z)\setminus\{z_0\}$ is understood as a Riemann surface
with a cylindrical end via the above biholomorphism $\R\times
S^1\rightarrow \C^\ast\setminus \{0\}$, as justified by Lemma
3.2.3.

Recall that we let $(TX)^\ast_{\xi}$ denote the pull-back orbifold
bundle defined by $(f,\xi)$. For  any $C^\infty$ section
$\tilde{s}$ of $(TX)^\ast_\xi$, let
$\tilde{f}_{\xi,s}=Exp\circ\bar{f}_\xi\circ \tilde{s}$,
$(TX)^\ast_{\xi,s}$ be the pull-back orbifold bundle defined by
$\tilde{f}_{\xi,s}$, and $Par_t: (TX)^\ast_\xi\rightarrow
(TX)^\ast_{\xi,ts}$ be the orbifold bundle isomorphism defined by
parallel transport along parametrized geodesics
$\tilde{\gamma}_z(t)=Exp\circ\bar{f}_\xi\circ t\tilde{s}(z)$ (cf.
Lemma 4.4.15). We define a map $F:L^p_1((TX)^\ast_\xi)\rightarrow
L^p((TX)^\ast_{\xi}\otimes \Lambda^{0,1}(\Sigma))$ (for $p>2$) by
$F(\tilde{s})=((Par_1)^{-1}\times Id)\circ\bar{\partial}
(\tilde{f}_{\xi,s})$, where $$
\bar{\partial}(\tilde{f}_{\xi,s})=\frac{1}{2}(d\tilde{f}_{\xi,s}+
J\circ d\tilde{f}_{\xi,s}\circ j)\in
C^\infty((TX)^\ast_{\xi,s}\otimes \Lambda^{0,1}(\Sigma)). $$ Then
$F$ is a nonlinear Fredholm map. Its linearization
$DF_{\tilde{s}}$ is a first order elliptic operator (in the
orbifold categary) whose index is given by the
Hirzebruch-Riemann-Roch formula: $$ index
(DF_{\tilde{s}})=2c_1(|(TX)^\ast_{\xi}|)([\Sigma])+2n(1-g_\Sigma).
$$ Here $2n$ is the dimension of $X$, and $|(TX)^\ast_{\xi}|$ is
the de-singularization of $(TX)^\ast_\xi$. (cf. Proposition 4.1.4
in [CR1].)

For the corresponding pseudo-holomorphic map $f_0$ from the
Riemann surface with a cylindrical end
$(\Sigma,\z)\setminus\{z_0\}$ into $X$, let $\xi_0$ denote the
isomorphism class of compatible systems induced by $\xi$; we have
the corresponding pull-back orbibundles $(TX)^\ast_{\xi_0}$,
$(TX)^\ast_{\xi_0,s}$, the $C^\infty$ maps
$(\tilde{f}_0)_{\xi_0,s}$, and the orbifold bundle isomorphisms
$Par_t$. Let $\delta\in (0,\epsilon_0)$ for a sufficiently small
$\epsilon_0$. We define a map $$ F_0:
L_{1,\delta}^p((TX)^\ast_{\xi_0};(TX)_{f(z_0)}^\xi)\rightarrow
L^p_\delta((TX)^\ast_{\xi_0}\otimes\Lambda^{0,1}(\Sigma\setminus\{z_0\}))
\hspace{2mm} (\mbox{ for }\hspace{2mm} p>2) $$ by
$F_0(\tilde{s})=((Par_1)^{-1}\times Id)\circ\bar{\partial}
((\tilde{f}_0)_{\xi_0,s})$, where
$L_{1,\delta}^p((TX)^\ast_{\xi_0};(TX)_{f(z_0)}^\xi)$ is the space
of local $L_1^p$ sections of $(TX)^\ast_{\xi_0}$ which
exponentially decay to an element in $(TX)_{f(z_0)}^\xi$ with a
weight $\delta$, $(TX)_{f(z_0)}^\xi$ is the linear subspace of
fixed points of the group homomorphism of $\xi$ at $z_0$ in the
fiber $TX_{f(z_0)}$. For $\tilde{s}\in
L_{1,\delta}^p((TX)^\ast_{\xi_0};(TX)_{f(z_0)}^\xi)$, let
$s_\infty$ be its limiting value at  infinity; we define its norm
by $$ ||\tilde{s}||=||\tilde{s}-s_\infty||_{L_{1,\delta}^p}+
||s_\infty||. $$

\noindent{\bf Lemma 3.2.3: }{\it When restricted to a small
neighborhood of the zero section in
$L_{1,\delta}^p((TX)^\ast_{\xi_0};(TX)_{f(z_0)}^\xi)$, the map
$F_0$ is a nonlinear Fredholm map with the same index as $F$.
Moreover, the zero set of $F_0$ is identical with the zero set of
$F$, i.e., $F_0^{-1}(0)=F^{-1}(0)$.}
    \vskip 0.1in

\noindent{\bf Proof:} In order to show that $F_0$ is Fredholm, we
need to look at the restriction of its linearization
$(DF_0)_{\tilde{s}}$ to the cylindrical end. Let $p=f(z_0)$ and
$(V,G,\pi)$ be a geodesic chart of $X$ at $p$, with complex
coordinates $u^\alpha, \alpha=1,\cdots, n$. Then in $V$ the
$\bar{\partial}$ operator is given by $$ \bar{\partial}(u^\alpha)
=  \frac{1}{2} (\frac{\partial u^\alpha}{\partial t}-
J(u^\alpha)\frac{\partial u^\alpha}{\partial s})dt+
\frac{1}{2}(\frac{\partial u^\alpha}{\partial s}+
J(u^\alpha)\frac{\partial u^\alpha}{\partial t})ds. $$ The
linearization $(D\bar{\partial})_{u^\alpha}$ at $(u^\alpha)$ is
\begin{eqnarray*}
(D\bar{\partial})_{u^\alpha}(v^\alpha)
 & = & \frac{1}{2}(\frac{\partial v^\alpha}{\partial t}-
J(u^\alpha)\frac{\partial v^\alpha}{\partial s}-(\partial_\beta J)(u^\alpha)
\frac{\partial u^\alpha}{\partial s}v^\beta)dt\\
 & + & \frac{1}{2}(\frac{\partial v^\alpha}{\partial s}+
J(u^\alpha)\frac{\partial v^\alpha}{\partial t}+(\partial_\beta J)(u^\alpha)
\frac{\partial u^\alpha}{\partial t}v^\beta)ds.
\end{eqnarray*}
As $t\rightarrow\infty$, $u^\alpha-u^\alpha_\infty$ is of
exponential decay of weight $\delta$. So
$(D\bar{\partial})_{u^\alpha}$ is a compact perturbation of an
operator whose restriction to the cylindrical end is $$
\frac{1}{2}(\frac{\partial v^\alpha}{\partial t}-
J(u^\alpha_\infty)\frac{\partial v^\alpha}{\partial s})dt+
\frac{1}{2}(\frac{\partial v^\alpha}{\partial s}+
J(u^\alpha_\infty)\frac{\partial v^\alpha}{\partial t})ds. $$ The
operator $J(u^\alpha_\infty)\frac{\partial}{\partial s}$ is a
small perturbation of $J(0)\frac{\partial}{\partial s}$ when
$u^\alpha_\infty$ is small, which has  spectrum
$\{\frac{\lambda}{m}:\lambda\in \Z\}$ where $m$ is the
multiplicity at  $z_0$ determined by $\xi$. So for any given
$\delta\in (0,\epsilon_0)$ for a sufficiently small $\epsilon_0$,
when $u^\alpha_\infty$ is small enough, $\delta$ is not in the
spectrum of $J(u^\alpha_\infty)\frac{\partial}{\partial s}$, so
that the operator whose restriction to the cylindrical end is $$
\frac{1}{2}(\frac{\partial v^\alpha}{\partial t}-
J(u^\alpha_\infty)\frac{\partial v^\alpha}{\partial s})dt+
\frac{1}{2}(\frac{\partial v^\alpha}{\partial s}+
J(u^\alpha_\infty)\frac{\partial v^\alpha}{\partial t})ds $$ is
Fredholm, which implies that $(D\bar{\partial})_{u^\alpha}$ is
Fredholm. Hence $(DF_0)_{\tilde{s}}$ is Fredholm for a given
weight $\delta\in (0,\epsilon_0)$ when $||s_\infty||$ is small. As
for the calculation of the index, we may assume that
$\tilde{s}=0$. Then $(DF_0)_{0}$ can be perturbed to an operator
whose restriction to the cylindrical end is $$
\frac{1}{2}(\frac{\partial v^\alpha}{\partial t}-
J_0\frac{\partial v^\alpha}{\partial s})dt+
\frac{1}{2}(\frac{\partial v^\alpha}{\partial s}+
J_0\frac{\partial v^\alpha}{\partial t})ds $$ where $J_0$ is the
standard complex structure on $\C^n$, while for $(DF)_0$ the
corresponding perturbed operator is of the form of the standard
$\bar{\partial}$ on a disc neighborhood of $z_0$. So they have
identical kernel and cokernel, hence the same index.

As for the zero sets, $F_0^{-1}(0)\subset F^{-1}(0)$ follows from
the removability of isolated singularites of pseudo-holomorphic
maps. Now let $\tilde{s}\in F^{-1}(0)$; restricting to
$\Sigma\setminus\{z_0\}$, $\tilde{s}$ can be thought as a section
of $(TX)^\ast_{\xi_0}$, and $F_0(\tilde{s})=0$. We need to show
that $\tilde{s}(z_0)\in (TX)_{f(z_0)}^\xi$ and
$\tilde{s}-\tilde{s}(z_0)$ exponentially decays with a weight at
least $\delta$, so that $\tilde{s}$ can be regarded as in
$L_{1,\delta}^p((TX)^\ast_{\xi_0};(TX)_{f(z_0)}^\xi)$. It is
obvious that $\tilde{s}(z_0)\in (TX)_{f(z_0)}^\xi$. On the other
hand, for any non-zero pseudo-holomorphic map $f: D\rightarrow
(\C^n,J)$ with $J(0)=J_0$ and $f(0)=0$, we have $f(z)=a\cdot
z^k+O(|z|^{k+1})$ for some $k>0$ and $a\in \C^n\setminus\{0\}$
(cf. [HW]), so that $\tilde{s}-\tilde{s}(z_0)$ decays
exponentially of  a weight at least $\frac{1}{m}>\delta$. Here $m$
is the multiplicity at $z_0$ determined by $\xi$. \hfill $\Box$

\vspace{2mm}

With this understood, we now move on to a different issue. Let
$f:(\Sigma,\z)\rightarrow X$ be a pseudo-holomorphic map with a
twisted boundary condition $\xi$, such that $f(\Sigma)$ lies
entirely in $\Sigma X$. Then after deleting finitely many points
$\z^\prime$ including $\z$, for any point $p$ in
$f(\Sigma\setminus\{\z^\prime\})$, $G_p$ is isomorphic to a fixed
group $G$. Then the isomorphism class of compatible systems $\xi$,
restricted to $\Sigma\setminus\{\z^\prime\}$, defines an
isomorphism class of fiber bundles over
$\Sigma\setminus\{\z^\prime\}$ with fiber $G$. The global
sections, which extend over  the whole $\Sigma$, form a group
isomorphic to a subgroup of $G$. We call this group the {\it
isotropy group} of $(f,\xi)$, denoted by $G_{(f,\xi)}$. Obviously,
$G_{(f,\xi)}$ acts linearly on the space of sections of the
pull-back orbifold bundle $(TX)^\ast_{f,\xi}$ defined by $(f,\xi)$
by pointwise multiplication.  Let $g\in G_{(f,\xi)}$, $z\in\z$,
and for a local representation of $\xi$, let the group
homomorphism of $\xi$ at $z$ be given by an element $\xi(z)\in
G_{f(z)}$; then within the local trivialization determined by the
local representation of $\xi$, we have $g(z)$ lying in the
centralizer of $\xi(z)$. Suppose $\sigma=(f,(\Sigma,\z),\xi)$ is a
stable map,  where $\Sigma=\cup\pi_\nu(\Sigma_\nu)$ is a
semistable curve, $\xi=(\xi_\nu)$, where $\xi_\nu$ is a twisted
boundary condition of $f_\nu=f\circ \pi_\nu$ satisfying Definition
2.3.3 (3). We let $G_\nu$ be the isotropy group of
$(f_\nu,\xi_\nu)$, and define $$ G_\sigma=\{(g_\nu)\in\prod_\nu
G_\nu|g_\nu(z_\nu)=g_\omega(z_\omega), \hspace{2mm} \mbox{ if
}\hspace{2mm} \pi_\nu(z_\nu)=\pi_\omega(z_\omega)\}, $$ and call
it the isotropy group of $\sigma$. The automorphism group
$Aut(\sigma)$ of $\sigma$ acts on $G_\sigma$ as automorphisms via
pull-backs.

Next we have a digression. Recall that a stable curve
$x=(\Sigma,\z)$ has a neighborhood in the Deligne-Mumford
compactification which is described as $$ \frac{V_x\times
B_r(\oplus_z T_{z_\nu}\Sigma_\nu\otimes
T_{z_\omega}\Sigma_{\omega})} {Aut(\Sigma,\z)}, $$ where $V_x$
parametrizes the deformation of complex structures on $x$, and
$B_r(\oplus_z T_{z_\nu}\Sigma_\nu\otimes
T_{z_\omega}\Sigma_{\omega})$ gives the parameters of resolving
the singularities. We denote $V_x$ by $V_{deform}$ and
$B_r(\oplus_z T_{z_\nu}\Sigma_\nu\otimes
T_{z_\omega}\Sigma_{\omega})$ by $V_{resolv}$. Given
$(y,\varsigma)\in V_{deform}\times V_{resolv}$, where $\varsigma=
(\sigma_z)\in\oplus_z T_{z_\nu}\Sigma_\nu\otimes
T_{z_\omega}\Sigma_\omega$, there is an associated stable curve
$(\Sigma_{y,\varsigma},\z_{y,\varsigma})$ obtained as follows.
Each component $\Sigma_{\nu}$ of $\Sigma_y$ is given a K\"{a}hler
metric $\mu_y$ which is flat in a neighborhood of the singular
points. This gives a Hermitian metric on each
$T_{z_\nu}\Sigma_{\nu}$. For each non-zero $\sigma_z\in
T_{z_\nu}\Sigma_{\nu}\otimes T_{z_\omega}\Sigma_{\omega}$, there
is a biholomorphic map $\Psi_{\sigma_z}: T_{z_\nu}\Sigma_{\nu}
\setminus\{0\} \rightarrow
T_{z_\omega}\Sigma_{\omega}\setminus\{0\}$ defined by
$u\otimes\Psi_{\sigma_z}(u)=\sigma_z$. Let $|\sigma_z|=R^{-2}$,
then for sufficiently large $R$, the map
$exp_{z_\omega}^{-1}\circ\Psi_{\sigma_z} \circ exp_{z_\nu}$ is a
biholomorphism between $D_{z_\nu}(R^{-1/2})\setminus
D_{z_\nu}(R^{-3/2})$ and $D_{z_\omega}(R^{-1/2})\setminus
D_{z_\omega}(R^{-3/2})$ where $D_{z_\nu}(R^{-1/2})$ is a disc
neighborhood of $z_\nu$ in $\Sigma_{\nu}$ of radius $(R^{-1/2})$
which is flat assuming $R$ is sufficiently large. We glue
$\Sigma_\nu$ and $\Sigma_\omega$ by this biholomorphism. If
$\sigma_z=0$, we do not make any change. Thus we obtain
$(\Sigma_{y,\varsigma},\z_{y,\varsigma})$. If we view a punctured
Riemann surface as a Riemann surface with cylindrical ends, then
the gluing in the above construction is equivalent to gluing
$(\frac{1}{2}\ln R,\frac{3}{2}\ln R)\times S^1$ to itself by
$t\rightarrow 2\ln R-t$ and $s\rightarrow -(s+\alpha)$ when
resolving $z$ with parameter $\sigma_z=R^{-2}e^{i\alpha}$. This
construction naturally gives a K\"{a}hler metric
$\mu_{y,\varsigma}$ on $\Sigma_{y,\varsigma}$ which coincides with
the K\"{a}hler metric $\mu_y$  on $\Sigma_y$ outside a
neighborhood of the cylindrical region. Each $\gamma\in
Aut(\Sigma,\z)$ takes $(\Sigma_y,\z_y)$ to
$(\Sigma_{\gamma(y)},\z_{\gamma(y)})$ isometrically, so it acts on
$\oplus_z T_{z_\nu}\Sigma_\nu\otimes T_{z_\omega}\Sigma_\omega$.
$\gamma$ induces an  isomorphism between
$(\Sigma_{y,\varsigma},\z_{y,\varsigma})$ and
$(\Sigma_{\gamma(y,\varsigma)},\z_{\gamma(y,\varsigma)})$, which
is also an isometry.

Now let $\sigma=(f,(\Sigma,\z),\xi)$ be the given stable map. We
will put a minimal number of additional marked points on each
unstable component of $(\Sigma,\z)$ so that it becomes stable, and
the images of these additional marked points under $f$ are
invariant under $Aut(\sigma)$. We denote the resulting stable
curve by $(\Sigma,\z)^+$. Then $(\Sigma,\z)^+$ has a neighborhood
described by $V_{deform}\times V_{resolv}/Aut((\Sigma,\z)^+)$.
Note that in the same way the automorphism group $Aut(\Sigma,\z)$
of $(\Sigma,\z)$, which may be of positive dimension,  also acts
on $V_{deform}\times V_{resolv}$. If
$(y,\varsigma)^\prime=\gamma(y, \varsigma)$ for some $\gamma\in
Aut(\Sigma,\z)$, then $(y,\varsigma)^\prime$ and $(y,\varsigma)$
represent the same equivalence class of stable curve, if and only
if $\gamma\in Aut((\Sigma,\z)^+)$.  If we drop  the newly added
marked points, then $(y,\varsigma)^\prime$ and $(y,\varsigma)$
represent the same equivalence class of stable curves if and only
if there is a $\gamma\in Aut(\Sigma,\z)$ such that
$(y,\varsigma)^\prime=\gamma(y,\varsigma)$. In particular, as a
subgroup of $Aut(\Sigma,\z)$, the automorphism group of $\sigma$,
$Aut(\sigma)$, acts on the space $V_{deform}\times V_{resolv}$.

Given a $(y,\varsigma)\in V_{deform}\times V_{resolv}$, we will
construct an ``almost'' pseudo-holomorphic map with a canonical
twisted boundary condition induced from $\xi$, $f_{y,\varsigma}:
(\Sigma_{y,\varsigma},\z_{y,\varsigma})\rightarrow X$, such that
$f_{y,\varsigma}=f_{\gamma(y,\varsigma)}\circ \gamma$ for any
$\gamma\in Aut(\sigma)$.  This is done as follows. Let
$\beta:\R\rightarrow [0,\infty)$ be a cut-off function such that
$\beta(t)=1$ for $t\leq 0$ and $\beta(t)=0$ for $t\geq
\frac{1}{2}$ such that $|\beta^\prime|\leq 4$. Let
$\varsigma=(\sigma_z)\in V_{resolv}$ such that
$|\sigma_z|=R^{-2}_z$ is sufficiently small for all singular point
$z$. Let $z=\pi_\nu(z_\nu)=\pi_\omega(z_\omega)$, $p_z=f(z)$, and
$(V_{p_z},G_{p_z}, \pi_{p_z})$ be a geodesic chart at $p_z$. Let
the group homomorphism of $\xi_\nu$ at $z_\nu$ be given by
$e^{2\pi i/n_\nu}\rightarrow g_\nu$ and the group homomorphism of
$\xi_\omega$ at $z_\omega$ be given by $e^{2\pi
i/n_\omega}\rightarrow g_\omega$; then $n_\nu=n_\omega=n_z$ and
$g_\nu=g_\omega^{-1}=g_z$ in $G_{p_z}$.  For large enough $R_z$,
both $f_\nu([\frac{1}{2}\ln R_z,\infty))$ and
$f_\omega([\frac{1}{2}\ln R_z,\infty))$ lie in
$\pi_{p_z}(V_{p_z})$. We define $f_{y,\varsigma}$ to be identical
to $f$ outside the neck region on $\Sigma_{y,\varsigma}$, and to
equal $\beta(t-\frac{1}{2}\ln R_z)f_\nu$ and
$\beta(t-\frac{1}{2}\ln R_z)f_\omega$ on the neck region. Note
that on the part of the gluing region $[\frac{1}{2}(\ln
R_z+1),\frac{1}{2}(3\ln R_z-1)]\times S^1$, $f_{y,\varsigma}\equiv
p_z$. $f_{y,\varsigma}$ inherits a twisted boundary condition
$\xi_{y,\varsigma}$ from $\xi$ canonically, whose representative
on the neck region can be obtained by gluing together
$(\tilde{f}_\nu,\xi_\nu)$ and $(\tilde{f}_\omega,\xi_\omega)$
suitably. Since the images of $f_\nu$ and $f_\omega$ lie in fixed
point sets which have a cone structure, so $\beta(t-\frac{1}{2}\ln
R_z)f_\nu$ and $\beta(t-\frac{1}{2}\ln R_z)f_\omega$ lie in the
same stratum of the canonical stratification
$X=X_{reg}\cup\widetilde{\Sigma X}_{gen}$ as $f_\nu$ and
$f_\omega$ do respectively, from which it follows that the
isotropy group $G_\sigma$ of $\sigma$ is naturally isomorphic to
the isotropy group $G_{y,\varsigma}$ of each
$(f_{y,\varsigma},(\Sigma_{y,\varsigma},\z_{y,\varsigma}),\xi_{y,\varsigma})$,
which is compatible with the action of $Aut(\sigma)$. Finally, it
is easy to see that $(f_{y,\varsigma},\xi_{y,\varsigma})
=(f_{\gamma(y,\varsigma)},\xi_{\gamma(y,\varsigma)})\circ \gamma$
for any $\gamma\in Aut(\sigma)$.

We define for $p>2$ $$ L^p_{1,\delta}((TX)^\ast_\sigma;(TX)^\xi)=
\{(u_\nu)\in \oplus_\nu L^p_{1,\delta}((TX)^\ast_{\xi_\nu};
(TX)^{\xi_\nu}_{f_\nu(z_\nu)})| u_{\nu,\infty}= u_{\omega,\infty},
\hspace{2mm} \mbox{if}\hspace{2mm}
\pi_\nu(z_\nu)=\pi_\omega(z_\omega)\}, $$ and $$
L^p_{\delta}((TX)^\ast_\sigma\otimes \Lambda^{0,1})= \{(u_\nu)\in
\oplus_\nu L^p_{\delta}((TX)^\ast_{\xi_\nu} \otimes
\Lambda^{0,1}(\Sigma_\nu\setminus\{z_v\})\}, $$ and we define
$L_{1,\delta}^p((TX)^\ast_{y,\varsigma};(TX)^{\xi_{y,\varsigma}})$
and $L_{\delta}^p((TX)^\ast_{y,\varsigma}\otimes \Lambda^{0,1})$
similarly, but with the understanding that the norm at each neck
region at $z$ is defined with a weight function $e^{\delta\tau}$,
where on the neck $(\frac{1}{2}\ln R_z, \frac{3}{2}\ln R_z)$ the
function $\tau$ is given by $$ \tau(t)=\left\{\begin{array}{cc} t
& t\in (\frac{1}{2}\ln R_z, \ln R_z-\frac{1}{2})\\ 2\ln R_z-t  &
t\in (\ln R_z+\frac{1}{2}, \frac{3}{2}\ln R_z)\\ \ln
R_z-\frac{1}{2} & t\in [\ln R_z-\frac{1}{2},\ln R_z+\frac{1}{2}].
\end{array} \right.\leqno (3.2.1)
$$ The group $G_\sigma$ acts linearly on
$L^p_{1,\delta}((TX)^\ast_\sigma;(TX)^\xi)$ and
$L^p_{\delta}((TX)^\ast_\sigma\otimes \Lambda^{0,1})$,  and
$G_{y,\varsigma}$ acts linearly on
$L_{1,\delta}^p((TX)^\ast_{y,\varsigma};(TX)^{\xi_{y,\varsigma}})$
and $L_{\delta}^p((TX)^\ast_{y,\varsigma}\otimes \Lambda^{0,1})$.
Moreover, the group $Aut(\sigma)$ acts linearly on
$L^p_{1,\delta}((TX)^\ast_\sigma;(TX)^\xi)$ and
$L^p_{\delta}((TX)^\ast_\sigma\otimes \Lambda^{0,1})$, covering
the action of $Aut(\sigma)$ on $(\Sigma,\z)$. If we let
$g\rightarrow \gamma^\ast(g)$ be the automorphism on $G_\sigma$
induced by pull-back via $\gamma\in Aut(\sigma)$, then for any
section $u$ in $L^p_{1,\delta}((TX)^\ast_\sigma;(TX)^\xi)$ or
$L^p_{\delta}((TX)^\ast_\sigma\otimes \Lambda^{0,1})$, we have $$
(\gamma_\ast)^{-1}\circ g \circ \gamma_\ast (u)=\gamma^\ast
(g)(u). $$ We define $\Gamma_\sigma$ to be the group generated by
$G_\sigma$ and $Aut(\sigma)$ with the above relation; then
$\Gamma_\sigma$ is a finite group and the short sequence $$
1\rightarrow G_\sigma\rightarrow \Gamma_\sigma\rightarrow
Aut(\sigma) \rightarrow 1 \leqno(3.2.2) $$ is exact. By
definition, $\Gamma_\sigma$ acts on
$L^p_{1,\delta}((TX)^\ast_\sigma;(TX)^\xi)$ and
$L^p_{\delta}((TX)^\ast_\sigma\otimes \Lambda^{0,1})$ linearly.
The group $\Gamma_\sigma$ also acts on $V_{deform}\times
V_{resolv}$ if we let the action of $G_\sigma$ be  trivial.
Finally, the maps $f_{y,\varsigma}$ are ``almost''
pseudo-holomorphic in the sense that $$
||\bar{\partial}_{\Sigma_{y,\varsigma}} \tilde{f}_{y,\varsigma}||
_{L^p_\delta(\Sigma_{y,\varsigma})} \leq C(\sum_z
e^{-1/2(m^{-1}_z-\delta)\ln R_z}+||y||) $$ where $z$ runs over all
the singular points on $\Sigma$, $m_z$ is the multiplicity at $z$,
and $\varsigma=(\sigma_z)$ with $|\sigma_z|=R^{-2}_z$.

Thus we obtain a family of nonlinear Fredholm maps
$$
F_\sigma: L^p_{1,\delta}((TX)^\ast_\sigma;(TX)^\xi)\rightarrow
L^p_{\delta}((TX)^\ast_\sigma\otimes \Lambda^{0,1})
\hspace{2mm}\mbox{and}\hspace{2mm}
F_{y,\varsigma}:
L_{1,\delta}^p((TX)^\ast_{y,\varsigma};(TX)^{\xi_{y,\varsigma}})\rightarrow
L_{\delta}^p((TX)^\ast_{y,\varsigma}\otimes \Lambda^{0,1})
$$
which are the nonlinear Cauchy-Riemann equations
associated to  $\sigma$ and
the ``almost'' pseudo-holomorphic maps $f_{y,\varsigma}$ with twisted
boundary conditions $\xi_{y,\varsigma}$, parametrized by $V_{deform}
\times V_{resolv}$. We will construct a Kuranishi model for this family.

First we show that $F_\sigma$ and $F_{y,\varsigma}$ have the same index.

\noindent{\bf Lemma 3.2.4: }{\it
Let $\sigma$ be in $\overline{\M}_{g,k}(X,J,A,\x)$. Write
$\iota(\x)=\sum_{i=1}^k\iota_{(g_i)}$ where $\x=(X_{(g_1)},\cdots,X_{(g_k)})$
and $\iota_{(g)}$ is the degree shifting number for $X_{(g)}$ (cf. [CR1]).
The nonlinear Fredholm maps $F_\sigma$ and $F_{y,\varsigma}$ have the same
index, which equals
$$
2c_1(TX)\cdot f_\ast([\Sigma])
 +2n(1-g_\Sigma)-2\iota(\x).
$$}
    \vskip 0.1in

\noindent{\bf Proof:} Without loss of generality, we may assume
that $\Sigma_{y,\varsigma}$ is obtained from $\Sigma$ by resolving
the only  singular point $z$, and assume that
$\pi_\nu(z_\nu)=\pi_\omega(z_\omega)=z$. Recall that, for a
pseudo-holomorphic map $f$ with a twisted boundary condition
$\xi_0$ from a smooth Riemann surface $\Sigma_0$, the index of the
nonlinear Cauchy-Riemann equation $F$ is given by $$
2c_1(|(TX)^\ast_{\xi_0}|)([\Sigma_0])+2n(1-g_{\Sigma_0}), $$ which
equals  $$ 2c_1(TX)\cdot
f_\ast([\Sigma_0])-2\sum_{i=1}^k\sum_{j=1}^n\frac{m_{i,j}}{m_i}
+2n(1-g_{\Sigma_0}) $$ where $m_i$ is the multiplicity of the
orbifold point $z_i$, and $m_{i,j}$ is obtained from the
representation of the group homomorphism of $\xi_0$ at $z_i$. (cf.
Proposition 4.1.4 in [CR1].) Hence we have
\begin{eqnarray*}
index F_\sigma & = & \sum_\nu index F_{\Sigma_\nu}-\dim (TX)^\xi_{f(z)}\\
               & = & \sum_\nu 2c_1(TX)\cdot (f_\nu)_\ast([\Sigma_\nu])
                     -2\sum_{i=1}^k\sum_{j=1}^n\frac{m_{i,j}}{m_i}
                     -2\sum_{j=1}^n\frac{(m_{z_\nu,j}+m_{z_\omega,j})}{m_z}\\
               &    &  + \sum_\nu 2n(1-g_{\Sigma_\nu})-\dim (TX)^\xi_{f(z)}\\
               & = & 2c_1(TX)\cdot f_\ast([\Sigma])
                     - 2\sum_{i=1}^k\sum_{j=1}^n\frac{m_{i,j}}{m_i}
                     +2n(1-g_\Sigma)
\end{eqnarray*}
where $i=1,\cdots,k$ runs over the index of the marked points
$\z$. Here we used the fact that $$
2\sum_{j=1}^n\frac{(m_{z_\nu,j}+m_{z_\omega,j})} {m_z}+\dim
(TX)^\xi_{f(z)}=2n. $$ The lemma follows from the fact that
$\iota(\x)= \sum_{i=1}^k\sum_{j=1}^n\frac{m_{i,j}}{m_i}$. \hfill
$\Box$

Let $E_\sigma$ be a subspace of
$L^p_{\delta}((TX)^\ast_\sigma\otimes \Lambda^{0,1})$ which satisfies the
following conditions:

\begin{itemize}
\item Let $L_\sigma$ be the linearization of $F_\sigma$ at $0$; then
$L^p_{\delta}((TX)^\ast_\sigma\otimes \Lambda^{0,1})$ is spanned by
$E_\sigma$ and the image of $L_\sigma$, which is a closed, finite codimensional
subspace of  $L^p_{\delta}((TX)^\ast_\sigma\otimes \Lambda^{0,1})$.
\item $E_\sigma$ is a finite dimensional
effective representation of $\Gamma_\sigma$, inherited
from the action of $\Gamma_\sigma$ on
$L^p_{\delta}((TX)^\ast_\sigma\otimes \Lambda^{0,1})$.
\item $E_\sigma$ consists of $C^\infty$ sections, and there exists a
sufficiently large number $R_\sigma>0$ such that the support of each section
in $E_\sigma$ is contained in $t\leq R_\sigma$ on the cylindrical ends.
\end{itemize}

Such a subspace $E_\sigma$ certainly exists.

\vspace{2mm}

The main result in the subsection on the construction of a
Kuranishi neighborhood is summarized in the following

\noindent{\bf Proposition 3.2.5: }{\it Let $V_\sigma^+$ be the
finite dimensional subspace in
$L^p_{1,\delta}((TX)^\ast_\sigma;(TX)^\xi)$ defined by
$V_\sigma^+=L^{-1}_\sigma(E_\sigma)$. Suppose that
$V_{deform}\times V_{resolv}$ consists of sufficiently small
$(y,\varsigma)$; in particular if $\varsigma =(\sigma_z)$, we
require that $|\sigma_z|=R^{-2}_z<R^{-2}_\sigma$. We take an
isomorphism $\theta_{y,\varsigma}:\Lambda^1(\Sigma)\rightarrow
\Lambda^1(\Sigma)$ so that $\theta_{y,\varsigma}$ induces an
isomorphism $\Lambda^{0,1}(\Sigma)|_{supp E_\sigma}\rightarrow
\Lambda^{0,1}(\Sigma_{y,\varsigma})|_{supp E_\sigma}$ which is
compatible with the action of $Aut(\sigma)$. Now
$E_{y,\varsigma}=\theta_{y,\varsigma} (E_\sigma)$ can be regarded
as a subspace of $L_{\delta}^p((TX)^\ast_{y,\varsigma}\otimes
\Lambda^{0,1}$. Let $$ Z_r=\{u|F_\sigma(u)\in E_\sigma,||u||\leq
r\}\cup_{(y,\varsigma)\in V_{deform}\times
V_{resolv}}\{u|F_{y,\varsigma}(u)\in E_ {y,\varsigma},||u||\leq
r\}. $$ Then for sufficiently small $r>0$, there is a
$\Gamma_\sigma$-invariant  open neighborhood $U^+_r$ in
$V_{deform}\times V_{resolv}\times V_\sigma^+$, and a
$\Gamma_\sigma$-equivariant one to one and onto map
$\psi_r:U^+_r\rightarrow Z_r$, which is a diffeomorphism when
restricted to the  ``slice'' over each $(y,\varsigma)$ or
$\sigma$. Let $s_r: U^+_r/\Gamma_\sigma\rightarrow
E_\sigma/\Gamma_\sigma$ be the composition of $\psi_r$ with
$F_\sigma$ and $(\theta_{y,\varsigma})^{-1}\circ F_{y,\varsigma}$.
Then $s_r$ is continuous. Finally, suppose each newly added marked
point $z$ on $\Sigma$ is chosen such that $f(z)$ and $f(z^\prime)$
have the same orbit type for nearby points $z^\prime$ and $f$ is
an embedding at $z$. For each $z$, choose a $G_{f(z)}$-invariant
subpace $W_z$ in $TX_{f(z)}$ of codimension two such that $W_z$ is
orthogonal to $Im \ df(z)$ and $\sqcup_z W_z$ is
$\Gamma_\sigma$-invariant. We define $$ V_\sigma=\{u\in
V_\sigma^+|u(z)\in W_z, \hspace{2mm} \mbox{ for all newly added
marked points $z$}\}. $$ Then $V_\sigma$ is a
$\Gamma_\sigma$-invariant subspace of $V_\sigma^+$. Let $U_r$ be
the intersection of $U_r^+$ with $V_{deform}\times V_{resolv}
\times V_\sigma$; then $\psi_r$ restricted to $U_r/\Gamma_\sigma$
induces a homeomorphism between $s_r^{-1}(0)\in U_r/\Gamma_\sigma$
and a neighborhood of $\sigma$ in the moduli space
$\overline{\M}_{g,k}(X,J,A,\x)$. The dimension of this Kuranishi
neighborhood is $$ 2c_1(TX)\cdot A
 +2(n-3)(1-g)+2k-2\iota(\x).
$$
Here $\iota(\x)$ is defined to be $\sum_{i=1}^k\iota_{(g_i)}$ for
$\x=(X_{(g_1)},\cdots,X_{(g_k)})$, where
$\iota_{(g)}$ is the degree shifting number for $X_{(g)}$ (cf. [CR1]).}
    \vskip 0.1in

\noindent{\bf Proof:} There are three technical issues in the
proof. The first one is to show that the Kuranishi model of each
individual map $F_\sigma$ or $F_{y,\varsigma}$ has a uniform size.
The second one is to show that $\psi_r$ is surjective to a
neighborhood of $\sigma$ in $\overline{\M}_{g,k}(X,J,A,\x)$. The
point here is that the  topology of
$\overline{\M}_{g,k}(X,J,A,\x)$ is defined by the weak $C^\infty$
topology plus the $C^0$ topology on the neck region, while the
topology of $U_r^+$ is given by weighted Sobolev norms. A delicate
estimate for pseudo-holomorphic maps is needed here. The third one
is to prove that when restricted to $s_r^{-1}(0)\cap
U_r/\Gamma_\sigma$, the map $\psi_r$ is a homeomorphism onto a
neighborhood of $\sigma$ in $\overline{\M}_{g,k} (X,J,A,\x)$,
which involves choosing representatives of stable maps.

We choose a $\Gamma_\sigma$-invariant decomposition
$L^p_{1,\delta}((TX)^\ast_\sigma;(TX)^\xi)= V_\sigma^+ \oplus
V^\prime_\sigma$ and a $\Gamma_\sigma$-invariant decomposition
$L_\delta^p((TX)^\ast_\sigma\otimes
\Lambda^{0,1})=E_\sigma^\prime\oplus E_\sigma$ with a
decomposition ratio $c$ such that
$L_\sigma:V_\sigma^\prime\rightarrow E_\sigma^\prime$ is an
isomorphism. Let $(y,\varsigma)$ be in $V_{deform}\times
V_{resolv}$, where $\varsigma=(\sigma_z)$ with
$|\sigma_z|=R_z^{-2}$. For each $R_z$, we define a cut-off
function $\rho_{R_z}$ such that $\rho_{R_z}\equiv 1$ for $t\leq
\ln R_z$ and $\rho_{R_z}\equiv 0$ for $t\geq \frac{1}{2}(3\ln
R_z-1)$, and $|\rho_{R_z}^\prime|\leq 4/\ln R_z$. With
$\rho_{R_z}$, for each section $u$ in
$L^p_{1,\delta}((TX)^\ast_\sigma;(TX)^\xi)$, we can construct a
section $\# u$ in $L^p_{1,\delta}((TX)^\ast_{y,\varsigma};
(TX)^{\xi_{y,\varsigma}})$ by
$\#u=\rho_{R_z}u_{z_\nu}+(1-\rho_{R_z})u_{z_\omega}$ on the neck
region at $z$. The map $\#: u\rightarrow \# u$ is
$\Gamma_\sigma$-equivariant. Now given a section $x$ in
$L_\delta^p((TX)^\ast_{y,\varsigma}\otimes \Lambda^{0,1})$, we can
think of $(\theta_{y,\varsigma})^{-1}(x)$ as a section in
$L_\delta^p((TX)^\ast_\sigma\otimes \Lambda^{0,1})$ which equals
zero on $t\geq \ln R_z$. Then $(\theta_{y,\varsigma})^{-1}(x)$
admits a decomposition as $(\theta_{y,\varsigma})^{-1}(x)=
L_\sigma u_0 + (\theta_{y,\varsigma})^{-1}(e_0)$ where $L_\sigma
u_0\in E_\sigma^\prime$ and $e_0\in E_{y,\varsigma}$, and $u_0\in
V_\sigma^\prime$. Let $L_{y,\varsigma}$ be the linearization of
$F_{y,\varsigma}$ at $0$, and consider $L_{y,\varsigma}(\# u_0)$.
Let $x_1$ be the difference $\theta_{y,\varsigma}(L_\sigma u_0)
-L_{y,\varsigma}(\# u_0)$. Then we have $x=L_{y,\varsigma}(\#
u_0)+e_0+x_1$, and $x_1$ satisfies the following estimate:
\begin{eqnarray*}
||x_1||_{L^p_\delta} &\leq & C(||y||+\sum_z(\ln R_z)^{-1})
                                  ||u_0||_{L^p_{1,\delta}}\\
                     &\leq & C(||y||+\sum_z(\ln R_z)^{-1})||L_\sigma^{-1}||
                            ||L_\sigma u_0||_{L^p_\delta}\\
                     &\leq & C(||y||+\sum_z(\ln R_z)^{-1})||L_\sigma^{-1}||c
                             ||x||_{L^p_\delta}.
\end{eqnarray*}
So when $(||y||+\sum_z(\ln R_z)^{-1})$ is sufficiently small,
we have $||x_1||\leq
\frac{1}{2}||x||$. We can continue this procedure and get
$x_1=L_{y,\varsigma}(\#u_1)+e_1+x_2$ and $||x_2||\leq\frac{1}{2}||x_1||$,
$\cdots,$ $x_n=L_{y,\varsigma}(\#u_n)+e_n+x_{n+1}$ and
$||x_{n+1}||\leq \frac{1}{2}||x_n||$, $\cdots$. Let
$\eta_1(x)=\#(\sum_{n=0}^\infty u_n)$ and $\eta_2(x)=\sum_{n=0}^\infty e_n$.
Then  we have
\begin{eqnarray*}
||\eta_1(x)||_{L^p_{1,\delta}} & \leq &
C\sum_{n=0}^\infty||u_n||_{L^p_{1,\delta}}
\leq C\sum_{n=0}^\infty ||L_\sigma^{-1}||||L_\sigma(u_n)||_{L^p_\delta}\\
                               & \leq & C\sum_{n=0}^\infty ||L_\sigma^{-1}||
                                        c||x_n||_{L^p_\delta}
                                  \leq 3Cc||L_\sigma^{-1}||||x||_{L^p_\delta}
\end{eqnarray*}
and $$ ||\eta_2(x)||_{L^p_\delta}\leq \sum_{n=0}^\infty
||e_n||_{L^p_\delta} \leq \sum_{n=0}^\infty
c||x_n||_{L^p_\delta}\leq 3c||x||_{L^p_\delta}, $$ and $$
x=L_{y,\varsigma}(\eta_1(x))+\eta_2(x). $$ It is clear that
$E_{y,\varsigma}\subset kernel(\eta_1)$ so that if we let
$E_{y,\varsigma}^\prime=Im(L_{y,\varsigma}\circ\eta_1)$, we obtain
a decomposition of
$L^p_\delta((TX)^\ast_{y,\varsigma}\otimes\Lambda^{0,1})=
E_{y,\varsigma}^\prime \oplus E_{y,\varsigma}$, which is
$\Gamma_\sigma$-equivariant and of uniformly bounded ratio. Let
$V_{y,\varsigma}^\prime=Im(\eta_1)$, $V_{y,\varsigma}^+
=L_{y,\varsigma}^{-1}(E_{y,\varsigma})$; then we have a
$\Gamma_\sigma$-equivariant decomposition of
$L_{1,\delta}^p((TX)^\ast_{y,\varsigma};(TX)^{\xi_{y,\varsigma}})
=V_{y,\varsigma}^+\oplus V_{y,\varsigma}^\prime$, and
$\eta_1=L_{y,\varsigma}^{-1}:E_{y,\varsigma}^\prime\rightarrow
V_{y,\varsigma}^\prime$ with uniformly bounded norm. The
decomposition $$
L_{1,\delta}^p((TX)^\ast_{y,\varsigma};(TX)^{\xi_{y,\varsigma}})
=V_{y,\varsigma}^+\oplus V_{y,\varsigma}^\prime $$ also has a
uniformly bounded ratio because it is given by $$
u=(u-\eta_1(L_{y,\varsigma}(u)))+\eta_1(L_{y,\varsigma}(u)), $$
and $||\eta_1(L_{y,\varsigma}(u))||\leq C||u||$. Finally, we can
easily verify that $||D^2F_{y,\varsigma}(0)||$ is uniformly
bounded.

So basically we have shown that this family of Kuranishi models
has a uniform size. In order to construct a Kuranishi neighborhood
for this family of maps, we only need to construct a family of
isomorphisms $\eta_{y,\varsigma}: V_\sigma^+\rightarrow
V_{y,\varsigma}^+$ with uniformly bounded norms. For $u\in
V_\sigma^+$, we define $\eta_{y,\varsigma}(u)= \#
u-\eta_1(L_{y,\varsigma}(\# u))$. We need to show that for
sufficiently small $(y,\varsigma)$, there is a constant $C$ such
that $||\eta_{y,\varsigma}(u)||\geq C||u||$. Suppose that this is
not the case. Then there is a sequence of $(y,\varsigma)_n$ going
to $\sigma$, a sequence of $u_n$ with $||u_n||=1$ and
$\eta_{(y,\varsigma)_n}(u_n)\rightarrow 0$. But this contradicts
the fact that $\eta_1(\theta_{(y,\varsigma)_n}(L_\sigma u_n))=0$
and
$||\theta_{(y,\varsigma)_n}(L_\sigma(u_n))-L_{(y,\varsigma)_n}(\#
u_n)|| =o(||u_n||)\rightarrow 0$. Hence we have constructed the
Kuranishi model of this family of maps. The map
$s_r:U_r^+/\Gamma_\sigma\rightarrow E_\sigma/\Gamma_\sigma$ is
obviously continuous.

Next we will show that if there is a sequence $\sigma_n$ in
$\overline{\M}_{g,k}(X,J,A,\x)$ that is convergent to $\sigma$,
which is represented by stable maps $(f_n,(\Sigma_n,\z_n),\xi_n)$
such that after adding extra marked points,
$(f_n^+,(\Sigma_n,\z_n)^+,\xi_n^+)$ converges to
$(f^+,(\Sigma,\z)^+,\xi^+)$, then for large enough $n$,
$(f_n,(\Sigma_n,\z_n),\xi_n)$ lies in the image of $\psi_r$. It is
easily seen that we may assume that $(f^+,(\Sigma,\z)^+,\xi^+)$ is
the stable map after we add the minimal number of marked points of
our choice. For large $n$, suppose $(\Sigma_n,\z_n)^+
=(\Sigma_{(y,\varsigma)_n},\z_{(y,\varsigma)_n})$; then it is
easily seen that $\tilde{f}_n$ equals
$Exp\circ\bar{f}_{(y,\varsigma)_n}\circ \tilde{s}_n$ for some
$\tilde{s}_n$ in $L_{1,\delta}^p((TX)^\ast_{(y,\varsigma)_n};
(TX)^{\xi_{(y,\varsigma)_n}})$. What we need to show is that for
sufficiently large $n$, $||\tilde{s}_n||_{L_{1,\delta}^p}<r$. For
this part, we need Lemma 2.3.11. Recall that if
$(f_n^+,(\Sigma_n,\z_n)^+,\xi_n^+)$ converges to
$(f^+,(\Sigma,\z)^+,\xi^+)$, then the following holds. First, for
each $\mu>0$, when $n$ is sufficiently large, the restriction of
$\tilde{f}^+_n$ to $\Sigma_{y_n,\varsigma_n}\setminus W_n(\mu)$
converges to $\tilde{f}^+$ in the $C^\infty$ topology as a
$C^\infty$ map with an isomorphism class of compatible systems.
Secondly, $\lim_{\mu\rightarrow 0}\limsup_{n\rightarrow\infty}
Diam(f_n(W_{z,n}(\mu)))=0$ for each singular point $z$ of
$\Sigma$. Here $$ W_{z,n}(\mu)=(D_{z_\nu}(\mu)\setminus
D_{z_\nu}(R_{z,n}^{-1}))\cup (D_{z_\omega}(\mu)\setminus
D_{z_\omega}(R_{z,n}^{-1})),\hspace{2mm}\mbox{and} \hspace{2mm}
W_n(\mu)=\cup_z W_{z,n}(\mu). $$ We pick a $\mu>0$ so that for
large $n$, $Diam(f_n(W_{z,n}(\mu)))<\epsilon$ where the $\epsilon$
here is referred to Lemma 2.3.11. Then on the neck region
$[\ln\mu^{-1},2\ln R_{z,n}-\ln\mu^{-1}]\times S^1$, we have $$
|\frac{\partial f_n}{\partial t}(t,s)|+|\frac{\partial
f_n}{\partial s} (t,s)| \leq
Ce^{-\frac{1}{m_z}(\tau(t)-\ln\mu^{-1})}, $$ where $\tau(t)$ is
defined in $(3.2.1)$. Then it follows that on the neck region $$
[2(1/m_z-\delta)^{-1}\ln\mu^{-1},2\ln R_{z,n}-
2(1/m_z-\delta)^{-1}\ln\mu^{-1}]\times S^1, $$ the
$L_{1,\delta}^p$ norm of $\tilde{s}_n$ is bounded by a term
$Ce^{-\ln\mu^{-1}}$. So we can pick $\mu>0$ small enough so that
$Ce^{-\ln\mu^{-1}}<\frac{r}{4}$. Then for this fixed $\mu$, by
weak $C^\infty$ convergence, for large enough $n$, the
$L_{1,\delta}^p$ norm of $\tilde{s}_n$ on the rest is bounded by
$\frac{r}{4}$, so that the whole $L_{1,\delta}^p$ norm of
$\tilde{s}_n$ is bounded by $r$. Hence
$(f_n,(\Sigma_n,\z_n),\xi_n)$ lies in the image of $\psi_r$ for
large $n$.

Now we prove that when restricted to $s_r^{-1}(0)\cap
U_r/\Gamma_r$, the map $\psi_r$ is a homeomorphism onto a
neighborhood of $\sigma$ in $\overline{\M}_{g,k}(X,J,A,\x)$. We
first clarify the issue that different sections of the same
orbifold bundle may define the same $C^\infty$ map under the
exponential map $Exp$. Suppose there are two sections
$\tilde{s}_1$ and $\tilde{s}_2$ such that
$Exp\circ\bar{f}_{y,\varsigma}\circ\tilde{s}_1=
Exp\circ\bar{f}_{y,\varsigma}\circ\tilde{s}_2$ gives the same
pseudo-holomorphic map $f^\prime$.  We take a geodesic compatible
system $\{\tilde{f}_{UU^\prime},\lambda\}$ of $(f,\xi)$. Then it
induces two isomorphic compatible systems
$\{\tilde{f}^\prime_{i,UU^\prime_i},\lambda_i\}$ for $i=1,2$ of
$f^\prime$ through $\tilde{s}_1$ and $\tilde{s}_2$. Now by the
unique continuity property of pseudo-holomorphic maps, we have an
isomorphism $\delta_U$ for each $U$ between the uniformizing
system $V^\prime_1$ of $U^\prime_1$ to $V_2^\prime$ of
$U^\prime_2$ such that $\tilde{f}^\prime_{2,UU^\prime_2}=
\delta_U\circ\tilde{f}^\prime_{1,UU^\prime_1}$, and $\delta_U$ can
be chosen so that
$\lambda_2=\delta\circ\lambda_1\circ(\delta)^{-1}$. Then one can
check that this means that the collection $\{\delta_U\}$ defines
an element $g$ in $G_\sigma$, and $\tilde{s}_2=g\cdot
\tilde{s}_1$. Now we prove that $\psi_r$ is a homeomorphism onto a
neighborhood of $\sigma$ in $\overline{\M}_{g,k}(X,J,A,\x)$ when
restricted to $s_r^{-1}(0)\cap U_r/\Gamma_\sigma$. For the
surjectivity, we observe that for a stable map
$(f^\prime,(\Sigma^\prime,\z^\prime),\xi^\prime)$ sufficiently
close to $\sigma$, $f^\prime(\Sigma^\prime)$ will intersect $Exp
(W_z)$ transversally at finitely many points, say
$f^\prime(z^\prime_i)$, for each newly added marked point $z$. We
mark the one, say $z^\prime$, on $\Sigma^\prime$ which is closest
to $z$ and obtain a stable curve $(\Sigma^\prime,\z^\prime)^+$. By
what we have shown in the previous paragraph, there is a
$(y,\varsigma,\tilde{s})\in V_{deform}\times V_{resolv}\times
V_\sigma^+$ such that
$((f^\prime)^+,(\Sigma^\prime,\z^\prime)^+,(\xi^\prime)^+)$ is
given by $Exp\circ\bar{f}_{y,\varsigma}\circ
\eta_{y,\varsigma}(\tilde{s})$. Note that by our choice of the new
marked point $z^\prime$, we have
$\eta_{y,\varsigma}(\tilde{s}(z))\in W_z$, but $\tilde{s}(z)$ may
not lie in $W_z$, which is the condition for $\tilde{s}$ to be in
$V_\sigma$. We fix this problem in the following way. Note that
$||\tilde{s}(z)-\eta_{y,\varsigma}(\tilde{s}(z))||\leq
C||(y,\varsigma)||$. We take a section $\tilde{s}^\prime\in
V_\sigma$ such that $||\tilde{s}^\prime(z)-\tilde{s}(z)||\leq
C||(y,\varsigma)||$, and consider the intersection of
$f^\prime(\Sigma^\prime)$ with $Exp\circ
\bar{f}_{y,\varsigma}\circ
\eta_{y,\varsigma}(t\tilde{s}^\prime(z))$, which is transversally
at finitely many points when $||(y,\varsigma)||$ is small enough.
We replace the new marked point $z^\prime$ by one of these new
intersection points which is closest to $z$, call it $z_1^\prime$,
then with $z_1^\prime$, $(\Sigma^\prime,\z^\prime)^+$ is
represented by $(\Sigma_{(y,\varsigma)_1},\z_{(y,\varsigma)_1})$
for some $(y,\varsigma)_1$ such that
$||(y,\varsigma)_1-(y,\varsigma)||\leq C ||(y,\varsigma)||^2$. Let
$((f^\prime)^+,(\Sigma^\prime,\z^\prime)^+,(\xi^\prime)^+)$ be
given by $Exp\circ\bar{f}_{(y,\varsigma)_1}\circ
\eta_{(y,\varsigma)_1}(\tilde{s}_1)$. Then
$\tilde{s}_1=(\eta_{(y,\varsigma)_1})^{-1}\circ
\eta_{y,\varsigma}(\tilde{s}^\prime)$. The upshot is that although
$\tilde{s}_1$ may not be in $V_\sigma$,  we can still find a
section $\tilde{s}_1^\prime\in V_\sigma$ such that
$||\tilde{s}_1^\prime-\tilde{s}_1||\leq C||(y,\varsigma)||^2$, an
improved estimate. We repeat this process and take the limit, we
find a $(y,\varsigma)_\infty$ and a section $\tilde{s}_\infty\in
V_\sigma$ such that
$((f^\prime)^+,(\Sigma^\prime,\z^\prime)^+,(\xi^\prime)^+)$ is
given by $Exp\circ\bar{f}_{(y,\varsigma)_\infty}\circ
\eta_{(y,\varsigma)_\infty}(\tilde{s}_\infty)$. Hence the
surjectivity of $\psi_r$. As for the injectivity of $\psi_r$, we
will first show that if $(y,\varsigma,\tilde{s})$ and
$\gamma(y,\varsigma,\tilde{s})$ represent the same stable map
$(f^\prime,(\Sigma^\prime,\z^\prime),\xi^\prime)$ for some
$\gamma\in Aut(\Sigma,\z)$ which lies in a small neighborhood of
identity (this corresponds to the effect of perturbing the newly
added marked points by the action of $\gamma$), then if both
$\tilde{s}$ and $\gamma(\tilde{s})$ lie in $V_\sigma$, we must
have $\gamma=id$. This roughly follows from the following
consideration: the difference between
$\eta_{y,\varsigma}(\tilde{s}(z))$ and
$\eta_{\gamma(y,\varsigma)}(\gamma(\tilde{s})(z))$ in the
direction transversal to $W_z$ is measured by $C||\gamma-id||$.
This must also be true for $\tilde{s}(z)$ and
$\gamma(\tilde{s})(z)$ on the one hand, since
$||\gamma(y,\varsigma)-(y,\varsigma)||\leq
C||\gamma-id||||(y,\varsigma)||$ so that the effect of
$\eta_{y,\varsigma}$ and $\eta_{\gamma(y,\varsigma)}$ can be
ignored, but on the other hand, by the assumption, both
$\tilde{s}$ and $\gamma(\tilde{s})$ lie in $V_\sigma$, a
contradiction. Now the injectivity of $\psi_r$ restricted to
$s_r^{-1}(0) \cap U_r/\Gamma_\sigma$ follows from the fact that,
up to a factor in $G_\sigma$, if $(y,\varsigma,\tilde{s})_1$ and
$(y,\varsigma,\tilde{s})_2$ represent the same stable map
$(f^\prime,(\Sigma^\prime,\z^\prime),\xi^\prime)$, then for some
$\gamma\in Aut(\sigma)$ and $\gamma_0\in Aut(\Sigma,\z)$ which
lies in a small neighborhood of the identity, we have
$(y,\varsigma,\tilde{s})_1=\gamma_0\circ\gamma
((y,\varsigma,\tilde{s})_2)$.

Finally, the dimension calculation of the Kuranishi neighborhood is
a routine business, which follows from Lemma 3.2.4.

\subsection{Construction of Kuranishi structure}

In this subsection, we will patch up the local Kuranishi
neighborhoods we constructed in the previous subsection to yield a
Kuranishi structure of the moduli space
$\overline{\M}_{g,k}(X,J,A,\x)$, which we will show to be stably
complex, therefore carrying a canonical orientation. The heart of
the construction is the fulfillment of the compatibility condition
in the definition of Kuranishi structure. For this purpose, each
local obstruction bundle $E_\sigma$ cannot be chosen arbitrarily.

In the previous subsection, we constructed for each equivalence
class of stable maps $\sigma\in \overline{\M}_{g,k}(X,J,A,\x)$ a
Kuranishi neighborhood. In particular, a neighborhood of $\sigma$
in $\overline{\M}_{g,k}(X,J,A,\x)$ is identified with
$s_r^{-1}(0)\cap U_r/\Gamma_\sigma$ via the map $\psi_r$ for some
sufficiently small $r>0$. We first cover the moduli space
$\overline{\M}_{g,k}(X,J,A,\x)$ by finitely many of this kind of
open set, say $\{O_i: i=1,\cdots,N\}$, at $\sigma_i\in
\overline{\M}_{g,k}(X,J,A,\x)$ such that there is a closed subset
$\hat{O}_i\subset O_i$ and $\overline{\M}_{g,k}(X,J,A,\x)$ is also
covered by $\{\hat{O}_i\}$. We assume that the Kuranishi
neighborhood at $\sigma_i$ is constructed with the choice of
obstruction bundle $E_i$ and a representative stable map
$(f_i,(\Sigma_i,\z_i),\xi_i)$. We denote by $U_i$ the Kuranishi
neighborhood at $\sigma_i$ such that $O_i=\psi_i(s_i^{-1}(0)\cap
U_i/\Gamma_i)$ where $\Gamma_i$ stands for $\Gamma_{\sigma_i}$ and
$s_i$ and $\psi_i$ stands for $s_{r_i}$ and $\psi_{r_i}$ for some
small $r_i>0$. Finally, for each equivalence class of stable maps
$\tau\in\overline{\M}_{g,k}(X,J,A,\x)$, we fix a representative
$(f_\tau,(\Sigma,\z)_\tau,\xi_\tau)$.

\vskip 0.1in \noindent{\bf Lemma 3.3.1: }{\it For any equivalence
class of stable maps $\tau\in\overline{\M}_{g,k}(X,J,A,\x)$, if
$\tau\in O_i$ for some $i$, then for any $(y,\varsigma,\tilde{s})$
in the Kuranishi neighborhood $U_i$ at $\sigma_i$ representing
$\tau$, there is a monomorphism $\Gamma_\tau\rightarrow \Gamma_i$
whose image is the isotropy subgroup of $(y,\varsigma,\tilde{s})$
in $U_i$.} \vskip 0.1in

\noindent{\bf Proof:} The existence of a monomorphism
$\Gamma_\tau\rightarrow \Gamma_i$ follows from the existence of
monomorphisms $Aut(\tau)\rightarrow Aut(\sigma_i)$ and
$G_\tau\rightarrow G_{\sigma_i}$. When $(\Sigma_i,\z_i)$ is a
stable curve (i.e. no new marked points are added to
$(\Sigma_i,\z_i)$), $Aut((\Sigma,\z)_\tau)$ is embedded into
$Aut(\Sigma_i,\z_i)$ as the isotropy subgroup of $(y,\varsigma)$,
and when $\tau$ is sufficiently close to $\sigma_i$, $Aut(\tau)$
as a subgroup of $Aut((\Sigma,\z)_\tau)$ has its image lying in
$Aut(\sigma_i)$, which also fixes $\tilde{s}$. When
$(\Sigma_i,\z_i)$ has an unstable component so that some new
marked points are added, then a biholomorphism between
$(\Sigma,\z)_\tau$ and $(\Sigma_{y,\varsigma},\z_{y,\varsigma})$
will add corresponding new marked points on $(\Sigma,\z)_\tau$,
and this induces a map from $Aut(\tau)$ into $Aut(\Sigma_i,\z_i)$,
whose image, when $\tau$ is sufficiently close to $\sigma_i$, lies
in $Aut(\sigma_i)$ and is a monomorphism, and it fixes
$(y,\varsigma,\tilde{s})$. As for the monomorphism
$G_\tau\rightarrow G_{\sigma_i}$, we consider the case when both
$\tau$ and $\sigma_i$ have only one component for simplicity. The
general case follows the same way. We recall the definition of the
isotropy group $G_\tau$ of
$\tau=(f_\tau,(\Sigma,\z)_\tau,\xi_\tau)$. After deleting finitely
many points, the image of $f_\tau$ has the same orbit type given
by a group, say $H_{\tau}$. Then the twisted boundary condition
$\xi_{\tau}$ defines a fiber bundle with fiber $H_{\tau}$ and
$G_{\tau}$ is just the group of global sections of this fiber
bundle. Now when $||\tilde{s}||$ is sufficiently small, all points
except for finitely many in
$Exp\circ\bar{f}_{i,(y,\varsigma)}\circ\tilde{s}$ have the same
orbit type which is given by a subgroup $H_\tau$ of the group
$H_{\sigma_i}$ defining the orbit type for $\sigma_i$. On the
other hand, the twisted boundary condition $\xi_\tau$ is also
compatible with the twisted boundary condition for $\sigma_i$ in
the sense that the corresponding fiber bundle with fiber $H_\tau$
is a subbundle of the fiber bundle with fiber $H_{\sigma_i}$. From
this it follows that $G_\tau$ is a subgroup of $G_i$. One can also
verify that if $g\in G_i$ fixes $\tilde{s}$, then $g$ lies in
$G_\tau$. Hence the lemma. \hfill $\Box$

Now suppose an equivalence class of stable maps $\tau$ is in
$O_i$, and is represented by $(y,\varsigma,\tilde{s})$ in $U_i\cap
s_i^{-1}(0)$. By parallel transport along geodesics
$Exp\circ\bar{f}_{i,(y,\varsigma)} \circ t\tilde{s}(z)$, we
transport the subspace $E_{i,(y,\varsigma)}$ of
$L^p_\delta((TX)^\ast_{i,(y,\varsigma)}\otimes\Lambda^{0,1}))$
into
$C^\infty((TX)^\ast_{(y,\varsigma,\tilde{s})}\otimes\Lambda^{0,1}(\Sigma_{y,
\varsigma}))$, where $(TX)^\ast_{(y,\varsigma,\tilde{s})}$ is the
pull-back orbifold bundle by
$Exp\circ\bar{f}_{i,(y,\varsigma)}\circ\tilde{s}$ over
$\Sigma_{y,\varsigma}$. We denote the resulting subspace by
$E_{i,(y,\varsigma,\tilde{s})}$. Then we have
$\gamma^\ast(E_{i,\gamma(y,\varsigma,\tilde{s})})
=E_{i,(y,\varsigma,\tilde{s})}$ for any $\gamma\in \Gamma_i$. Let
$\theta_\tau:(\Sigma,\z)_\tau\rightarrow
(\Sigma_{y,\varsigma},\z_{y, \varsigma})$ be a biholomorphism
sending marked points to marked points (not including new marked
points on $\Sigma_{y,\varsigma}$). Then the pull-back
$\theta_\tau^\ast(E_{i,(y,\varsigma,\tilde{s})})$ in
$C^\infty((TX)^\ast_\tau\otimes\Lambda^{0,1}(\Sigma_\tau))$ does
not depend on the choice of $(y,\varsigma,\tilde{s})$ and
$\theta_\tau$, where $(TX)^\ast_\tau$ is the pull-back orbifold
bundle by $(f_\tau,(\Sigma,\z)_\tau, \xi_\tau)$. We denote
$\theta_\tau^\ast(E_{i,(y,\varsigma,\tilde{s})})$ by $E_{\tau,i}$.

\vskip 0.1in
\noindent{\bf Lemma 3.3.2: }{\it
We can choose each obstruction bundle $E_i$ suitably so that for any $\tau
\in \cap_i O_i$, $i=i_1,\cdots,i_l$, the subspaces $E_{\tau,i}$
of $C^\infty((TX)^\ast_\tau\otimes\Lambda^{0,1}(\Sigma_\tau))$,
for $i=i_1,\cdots,i_l$, are linearly independent.}
\vskip 0.1in

\noindent{\bf Proof:}
For each $\tau\in\cap_i O_i$, $i=i_1,\cdots,i_l$,
we can choose finitely many points $w_i$ on $\Sigma_\tau$, away
from any neck or singular point region and marked points, such that in a disc
neighborhood $D_i$ of $w_i$, none of the sections in each $E_{\tau,i}$ vanishes
identically. For each $i$, we take a cut-off function $\beta_i$ which is
zero in $D_i$. We multiply $\beta_i$ to each section in $E_{\tau,j}$ for
any $j\in\{i_1,\cdots,i_l\}\setminus\{i\}$. We transport these changes on
each $E_{\tau,i}$ back to each $E_i$. Then for such a choice of $E_i$, the
lemma holds for any stable map in a neighborhood of $\tau$ in
$\overline{\M}_{g,k}(X,J,A,\x)$. By the compactness of $\overline{\M}_{g,k}
(X,J,A,\x)$, we can cover $\overline{\M}_{g,k}(X,J,A,\x)$ by finitely
many such neighborhoods. Hence the lemma.
\hfill $\Box$

We define $E_\tau$ to be the span of $E_{\tau,i}$ for those $i$'s such that
$\tau\in\cap_i \hat{O}_i$. Note that here we use the closed subset $\hat{O}_i$.
We will use $E_\tau$ for the obstruction bundle at $\tau$. We first collect
a few properties of $E_\tau$.

\begin{itemize}
\item Let $L_\tau$ be the linearization of the non-linear Cauchy-Riemann
operator $\bar{\partial}_\tau$ at $0$; then \linebreak
$L^p_{\delta}((TX)^\ast_\tau\otimes \Lambda^{0,1}(\Sigma_\tau))$
is spanned by $E_\tau$ and the image of $L_\tau$, which is a
closed, finite codimensional subspace of
$L^p_{\delta}((TX)^\ast_\tau\otimes \Lambda^{0,1}(\Sigma_\tau))$.
\item $E_\tau$ is a finite dimensional
effective representation of $\Gamma_\tau$, inherited
from the action of $\Gamma_\tau$ on
$L^p_{\delta}((TX)^\ast_\tau\otimes \Lambda^{0,1}(\Sigma_\tau))$.
\item $E_\tau$ consists of $C^\infty$ sections, and there exists a
sufficiently large number $R_\tau>0$ such that the support of each section
in $E_\tau$ is contained in $t\leq R_\tau$ on the cylindrical ends.
\end{itemize}

Next we introduce the space of equivalence classes of ``perturbed'' stable
maps. Consider a $C^\infty$ map $\tilde{h}:(\Sigma,\z)\rightarrow X$ with a
twisted boundary condition $\xi$, where $(\Sigma,\z)$ is a semistable curve,
and $\xi=(\xi_\nu)$ is a collection of twisted boundary conditions satisfying
similar compatibility conditions as in the definition of stable maps.
We also require that if $\tilde{h}$ is constant on a component $\Sigma_\nu$,
then $\Sigma_\nu$ has at least three special points (marked or singular).
We further assume that $h$ is pseudo-holomorphic in a neighborhood of each
singular point on $\Sigma$. The equivalence relation is defined as follows.
Suppose $(\tilde{h},(\Sigma,\z),\xi)$ and $(\tilde{h}^\prime,
(\Sigma^\prime,\z^\prime),\xi^\prime)$ are two ``perturbed'' stable maps.
We say they are equivalent if
there is an isomorphism $\theta:(\Sigma,\z)\rightarrow
(\Sigma^\prime,\z^\prime)$ such that $(\tilde{h},\xi)
=(\tilde{h}^\prime,\xi^\prime)\circ\theta$ and $\theta^\ast(\bar{\partial}
\tilde{h}^\prime)=\bar{\partial}\tilde{h}$, where $\theta^\ast$ is the
pull-back map from $(TX)^\ast_{\xi^\prime}\otimes\Lambda^{0,1}(\Sigma^\prime)$
to $(TX)^\ast_\xi\otimes\Lambda^{0,1}(\Sigma)$. We can define a topology on
the space of equivalence classes of ``perturbed'' stable maps in the same
way as we did for the moduli space of stable maps.

Here we insert a lemma which will play the role of the unique continuity
property for pseudo-holomorphic maps.

\vskip 0.1in
\noindent{\bf Lemma 3.3.3: }{\it
Let $D$ be the unit disc in $\C$, and $\tilde{f}_i$ for $i=1,2$ be two
$C^\infty$ maps from $D$ into a uniformizing system $(V,G,\pi)$. Then if
$\bar{\partial}\tilde{f}_1=\bar{\partial}\tilde{f}_2$ and $\tilde{f}_1-
\tilde{f}_2$ vanishes to an infinite order at $0$, then $\tilde{f}_1\equiv
\tilde{f}_2$.}
\vskip 0.1in

\noindent{\bf Proof:} Written down in local coordinates, the
condition $\bar{\partial}\tilde{f}_1=\bar{\partial}\tilde{f}_2$
means $$ \frac{\partial u^\alpha_1}{\partial
t}+J(u^\alpha_1)\frac{\partial u^\alpha_1} {\partial s}=
\frac{\partial u^\alpha_2}{\partial t}+J(u^\alpha_2)\frac{\partial
u^\alpha_2} {\partial s} $$ and $$ \frac{\partial
u^\alpha_1}{\partial s}-J(u^\alpha_1)\frac{\partial u^\alpha_1}
{\partial t}= \frac{\partial u^\alpha_2}{\partial
s}-J(u^\alpha_2)\frac{\partial u^\alpha_2} {\partial t} $$ from
which we can derive the condition needed in the Hartman-Wintner
lemma $$ ||\Delta v||\leq C(||\partial_t v||+||\partial_s
v||+||v||) $$ where $v=(u^\alpha_1-u^\alpha_2)$. \hfill $\Box$

\vskip 0.1in \noindent{\bf Remark 3.3.4: }{\it There is a
technical issue involved in the construction. For any two sections
$\tilde{s}_1$ and $\tilde{s}_2$ near a stable map $\tau$ that
represent the same $C^\infty$ map, we need to show that there is
an element $g\in G_\tau$ such that
$\tilde{s}_2=g\cdot\tilde{s}_1$. For this we need a certain unique
continuity property, which fails for general $C^\infty$ maps. This
lemma tells us that we can overcome this problem by controlling
the image under the $\bar{\partial}$ operator, which lies in the
obstruction bundle $E$. So we need to require that any obstruction
bundle under consideration consists of sections such that if any
local representatives of two sections are related by an
isomorphism $\delta$, it must be induced from a global one, i.e.,
an element in  $G_\tau$.} \hfill $\Box$ \vskip 0.1in

For each equivalence class of ``perturbed'' stable maps $\kappa$,
we fix a representative
$(h_\kappa,(\Sigma,\z)_\kappa,\xi_\kappa)$. Suppose $\kappa$ is
close to $\sigma_i$ so that there is a $(y,\varsigma,\tilde{s})$
such that $\kappa$ is represented by
$h_{(y,\varsigma,\tilde{s})}=Exp\circ\bar{f}_{i,(y,\varsigma)}\circ\tilde{s}$.
We add $(\theta_{y,\varsigma}\circ Par)^{-1}(\bar{\partial}
h_{(y,\varsigma,\tilde{s})})$ to $E_i$ and denote by
$\hat{E}_{i,\kappa}$ the finite dimensional $\Gamma_i$-invariant
subspace generated by them. We repeat the local construction in
the previous subsection with the choice of $\hat{E}_{i,\kappa}$,
and we obtain finitely many, $\Gamma_i$-invariant representatives
$(y,\varsigma,\tilde{s})_j$ for $\kappa$. Then we can do the same
thing to $\kappa$ as we did to each stable map $\tau$. We use
parallel transport along geodesics
$Exp\circ\bar{f}_{i,(y,\varsigma)_j} \circ t\tilde{s}_j(z)$ to
transport the subspace $E_{i,(y,\varsigma)_j}$ of
$L^p_\delta((TX)^\ast_{i,(y,\varsigma)_j}\otimes\Lambda^{0,1}))$
into
$C^\infty((TX)^\ast_{(y,\varsigma,\tilde{s})_j}\otimes\Lambda^{0,1}(\Sigma_{y,
\varsigma})_j)$, where $(TX)^\ast_{(y,\varsigma,\tilde{s})_j}$ is
the pull-back orbifold bundle by
$Exp\circ\bar{f}_{i,(y,\varsigma)_j}\circ\tilde{s}_j$ over
$\Sigma_{(y,\varsigma)_j}$. We denote the resulting subspace by
$E_{i,(y,\varsigma,\tilde{s})_j}$. Then we have
$\gamma^\ast(E_{i,\gamma(y,\varsigma,\tilde{s})_j})
=E_{i,(y,\varsigma,\tilde{s})_j}$ for any $\gamma\in \Gamma_i$.
Let $\theta_{j,\kappa}:(\Sigma,\z)_\kappa\rightarrow
(\Sigma_{(y,\varsigma)_j},\z_{(y,\varsigma)_j})$ be a
biholomorphism sending marked points to marked points (not
including new marked points on $\Sigma_{(y,\varsigma)_j}$). Then
the pull-back
$\theta_{j,\kappa}^\ast(E_{i,(y,\varsigma,\tilde{s})_j})$ in
$C^\infty((TX)^\ast_\kappa\otimes\Lambda^{0,1}(\Sigma_\kappa))$
does not depend on the choice of $j$ and $\theta_{j,\kappa}$,
where $(TX)^\ast_\kappa$ is the pull-back orbifold bundle by
$(h_\kappa,(\Sigma,\z)_\kappa,\xi_\kappa)$. We denote
$\theta_\kappa^\ast(E_{i,(y,\varsigma,\tilde{s})_j})$ by
$E_{\kappa,i}$. The above choice of representatives
$(y,\varsigma,\tilde{s})_j$ of $\kappa$ are canonical in the sense
that both of the following conditions are satisfied: they are
$\Gamma_i$-invariant, and when $\kappa$ converges to a stable map
$\tau\in O_i$, the representatives of $\kappa$ converge to the
representatives of $\tau$. However, these representatives depend
on the choice of $(y,\varsigma,\tilde{s})$ at the beginning.

Now for each equivalence class of stable maps $\tau$, we consider
all the equivalence classes of ``perturbed'' stable maps $\kappa$
which are sufficiently close to $\tau$, and there is a
representative $(y,\varsigma, \tilde{s})$ centered at $\tau$. We
denote this set by $\B_\tau$ (we can think of it in the sense of
germs). For each $\kappa\in\B_\tau$, we construct an obstruction
space $E_{\kappa,\tau}$ as a finite dimensional subspace in
$C^\infty((TX)^\ast_\kappa\otimes\Lambda^{0,1}(\Sigma_\kappa))$ as
follows. Suppose $\tau\in \cap_i \hat{O}_i$ for
$i=i_1,\cdots,i_l$. Then we define $E_{\kappa,\tau}$ to be the
span of $E_{\kappa,i}$ for $i=i_1,\cdots,i_l$. Note that
$E_{\kappa,i}$ are linearly independent when $\kappa$ is
sufficiently close to $\tau$. As an immediate consequence, which
is crucial to the fulfillment of the compatibility requirement  in
the construction of Kuranishi structure, we have
$E_{\kappa,\tau_1}\subset E_{\kappa,\tau_2}$ when $\tau_1$ is
sufficiently close to $\tau_2$.

\vspace{2mm}

For each $\tau\in\overline{\M}_{g,k}(X,J,A,\x)$, we define $Z_\tau$ to be
the solution set of $Z_\tau=\{\kappa\in \B_\tau|\bar{\partial}h_\kappa\in
E_{\kappa,\tau}\}$. The main result of this subsection on the construction
of Kuranishi structure is summarized in

\vspace{2mm}

\noindent{\bf Proposition 3.3.5: }{\it Each $Z_\tau$ is an
orbifold and $E_\tau/\Gamma_\tau$ is an orbifold bundle over
$Z_\tau$ with a continuous section $s_\tau$ and a map $\psi_\tau$
such that $\psi_\tau$ is a homeomorphism from $s_\tau^{-1}(0)$
onto a neighborhood of $\tau$ in $\overline{\M}_{g,k}(X,J,A,\x)$,
i.e., $(Z_\tau,E_\tau/\Gamma_\tau,s_\tau, \psi_\tau)$ is a
Kuranishi neighborhood at $\tau$. Moreover, for any
$\sigma\in\overline{\M}_{g,k}(X,J,A,\x)$ which is sufficiently
close to $\tau$, we have an orbifold embedding
$(\phi_{\tau\sigma}, \hat{\phi}_{\tau\sigma})$ from
$(Z_\sigma,E_\sigma/\Gamma_\sigma)$ into
$(Z_\tau,E_\tau/\Gamma_\tau)$, where $\phi_{\tau\sigma}$ is given
by the natural inclusion. There is also an orbifold bundle
isomorphism $\Phi_{\tau\sigma}:TZ_\tau/TZ_\sigma\rightarrow
(E_\tau/\Gamma_\tau)/(E_\sigma/\Gamma_\sigma)$, and the collection
$$ \{(Z_\tau,E_\tau/\Gamma_\tau,
s_\tau,\psi_\tau,\phi_{\tau\sigma}, \hat{\phi}_{\tau\sigma},
\Phi_{\tau\sigma}):\tau \in\overline{\M}_{g,k}(X,J,A,\x)\} $$
forms a stably complex Kuranishi structure on
$\overline{\M}_{g,k}(X,J,A,\x)$ of dimension $2d$ for $$
d=c_1(TX)\cdot A +(\frac{1}{2}\dim_\R X-3)(1-g)+k-\iota(\x). $$
Here $\iota(\x)=\sum_{i=1}^k \iota_{(g_i)}$ for
$\x=(X_{(g_1)},\cdots, X_{(g_k)})$ where $\iota_{(g)}$ is the
degree shifting number for $X_{(g)}$.} \vskip 0.1in

\noindent{\bf Proof:} We first show that $Z_\tau$ is an orbifold.
Suppose that $\tau$ is in $O_i$. For any $(y,\varsigma,\tilde{s})$
centered at $\tau$, we construct a sequence
$(y,\varsigma,\tilde{s}_n)$ by cutting down $\tilde{s}$ on the
cylindrical ends such that $\tilde{s}_n$ converges to $\tilde{s}$.
Then the ``perturbed'' stable maps defined by
$(y,\varsigma,\tilde{s}_n)$ can be represented by some
$(y^\prime,\varsigma^\prime,\tilde{s}^\prime)_n$ centered at
$\sigma_i$ when $\tilde{s}$ is small, and by repeating the local
construction we obtain some canonical representatives which are
$\Gamma_i$-invariant. When $n\rightarrow\infty$, these
representatives converge to a set of $\Gamma_i$-invariant
representatives for $(y,\varsigma,\tilde{s})$. Similar to what we
did to a ``perturbed'' stable map, we obtain a finite dimensional
subspace in $C^\infty((TX)^\ast_{(y,\varsigma,\tilde{s})}\otimes
\Lambda^{0,1})$; call it $E_{(y,\varsigma,\tilde{s}),\tau}$. If
$(y,\varsigma,\tilde{s})$ represents a ``perturbed'' stable map
$\kappa$, then $E_{(y,\varsigma,\tilde{s}),\tau}$ is just the
pull-back of $E_{\kappa,\tau}$. This system of subspaces is also
$\Gamma_\tau$-equivariant. We want to consider the solution set of
$\{(y,\varsigma,\tilde{s})|
\bar{\partial}h_{(y,\varsigma,\tilde{s})}\in
E_{(y,\varsigma,\tilde{s})}\}$, where
$h_{(y,\varsigma,\tilde{s})}$ is the $C^\infty$ map defined by
$(y,\varsigma,\tilde{s})$. The technique of the local construction
can be modified to show that this solution set (or a slice of it
when $(\Sigma,\z)_\tau$ has an unstable component) is identified
with a $\Gamma_\tau$-invariant finite dimensional open ball
$U_\tau$, and $U_\tau/\Gamma_\tau$ is homeomorphic to $Z_\tau$.
Moreover, $E_\tau/\Gamma_\tau$ is an orbifold bundle over $Z_\tau$
with a continuous section $s_\tau$ and a map $\psi_\tau$ such that
$\psi_\tau$ is a homeomorphism from $s_\tau^{-1}(0)$ onto a
neighborhood of $\tau$ in $\overline{\M}_{g,k}(X,J,A,\x)$. So we
have shown that $(Z_\tau,E_\tau/\Gamma_\tau,s_\tau,\psi_\tau)$ is
a Kuranishi neighborhood at $\tau$. The dimension of the Kuranishi
structure is easily seen to be independent of $\tau$.

As for the maps $\phi_{\tau\sigma}, \hat{\phi}_{\tau\sigma}$ and
$\Phi_{\tau\sigma}$, we observe that if $\kappa\in\B_\sigma$, this
means that $\kappa$ is represented by some
$(y,\varsigma,\tilde{s})$ with $\tilde{s}$ sufficiently small. If
$\sigma$ is close enough to $\tau$, $\sigma$ is also represented
by a $C^\infty$ section, so that $\kappa$ can be also represented
by some $(y^\prime,\varsigma^\prime,\tilde{s}^\prime)$ centered at
$\tau$. This means that $\B_\sigma$ is naturally included in
$\B_\tau$ when $\sigma$ is sufficiently close to $\tau$. Then
$E_{\kappa,\sigma}\subset E_{\kappa,\tau}$ implies that $Z_\sigma$
is naturally included in $Z_\tau$, which is defined to be
$\phi_{\tau\sigma}$. In order to show that $\phi_{\tau\sigma}$ is
an orbifold embedding, we observe that the local group
$\Gamma_\sigma$ is isomorphic to the isotropy subgroup of any
representative of $\sigma$ in $U_\tau$ (Lemma 3.3.1). This implies
the embedding between the germs of the corresponding uniformizing
systems, because they are totally determined by the local groups
via geodesic charts. Similarly $\hat{\phi}_{\tau\sigma}$ is
obtained via natural inclusions $E_{\kappa,\sigma}\subset
E_{\kappa,\tau}$. The map $\Phi_{\tau\sigma}$ is obtained through
the following consideration. Suppose that $F$ is a nonlinear
Fredholm map with linearization $L$ at $0$, and a Kuranishi model
is constructed with a choice of obstruction space $E_1$ and $E_2$
respectively, on a small ball of $0$ in $V_1=L^{-1}(E_1)$ and
$V_2=L^{-1}(E_2)$ respectively. When $E_1$ is a subspace of $E_2$,
$V_1$ is a subspace of $V_2$ too, and the complement of $V_1$ in
$V_2$ is isomorphic to the complement of $E_1$ in $E_2$ via $L$.
However, there is a slight complication here. When
$(\Sigma,\z)_\tau$ has an unstable component, new marked points
are added, and only a subspace of $V$, which is complementary to
the direction of parametrization of the new marked points,  is
used. But when $(\Sigma,\z)_\sigma$ is a stable curve, there is no
need for new marked points, and the whole space  $V$ is used. The
remedy to this problem is this: We add new marked points to
$\sigma$ accordingly, and one can verify that the new Kuranishi
model is canonically isomorphic to the old one, because after
adding new marked points, the space of $V_{deform}\times
V_{resolv}$ is ``thickened'' by the direction of parametrization
of new marked points, which is isomorphic to the direction in $V$
which is complementary to the subspace in $V$ used in the new
Kuranishi model. The compatibility requirement is fulfilled
automatically because all of the maps are obtained through natural
inclusions.

Finally, we show that the Kuranishi structure thus obtained is
stably complex. The key point in this fact is that the symbol of
the linearization of the Cauchy-Riemann equation is complex linear
and the space $V_{deform} \times V_{resolv}$ is naturally a
complex manifold. Let's first look at a toy model. Suppose we have
an orbifold $U$ and a smooth family of Fredholm operators
$F_s:X_s\rightarrow Y_s$ for $s\in U$ which can be lifted to local
charts of $U$ equivariantly. If $E\rightarrow U$ is a finite
dimensional orbifold bundle over $U$ and for each $s\in U$ the
fiber $E_s\subset Y_s$ such that $F_s:X_s\rightarrow Y_s/E_s$ is
surjective, then $V=\{V_s\}$ is an orbifold bundle over $U$, where
$V_s=F_s^{-1}(E_s)$. If we deform $F=\{F_s\}$ through a homotopy
$F^t=\{F_s^t\}$, and let $V^t=(F^t)^{-1}(E)$, then $V^t$ are
isomorphic. In the case we don't have an orbifold $U$, but a
Kuranishi structure in general, we can do this on each Kuranishi
neighborhood compatibly. In the present case, we take a complex
linear vector space for each $E_i$, and use the complex linear
part of the geodesic parallel transport to transport each $E_i$ to
nearby points, so we can assume each $E_{\kappa,\tau}$ is complex
linear. Over each Kuranishi neighborhood $Z_\tau$, the tangent
bundle $TZ_\tau$ is given by $T(V_{deform}\times
V_{resolv}/\Gamma_\tau)\oplus V_\tau$, where $V_\tau$ is an
orbifold bundle obtained from the family of Fredholm operators
$D\bar{\partial}h_\kappa$ at $([y,\varsigma],\kappa)$, which can
be deformed to a family of complex linear Fredholm operators since
the symbols are complex linear. Hence $TZ_\tau$ is stably
isomorphic to a complex orbifold bundle. However, there is a
slight complication here which does not effect the conclusion,
namely, when $(\Sigma,\z)_\tau$ is a stable curve, we take
$V_\tau$ to be the orbifold bundle
$\{V_{([y,\varsigma],\kappa)}=(D\bar{\partial}
h_\kappa)^{-1}(E_{\kappa,\tau})\}$; but when $(\Sigma,\z)_\tau$
has a unstable component, we take a subbundle of it which is
stably complex. The complement is made up through the
``thickening'' in  the space $V_{deform}\times V_{resolv}$ by the
parametrization of new marked points. See [FO] for more details.
\hfill $\Box$

\subsection{Orbifold Gromov-Witten invariants and axioms}
    Once we construct the Kuranishi structure with the necessary patching properties,
we can use Fukaya-Ono's abstract argument to construct a virtual fundamental
    cycle of the moduli space of orbifold stable maps
in $H_*(\overline{\M}_{g,k}(X,J,A,\x),\Q)$,
    denoted by $[\overline{\M}_{g,k}(X,J,A,\x)]^{vir}$. The degree of
$[\overline{\M}_{g,k}(X,J,A,\x)]^{vir}$ is given by the index
formula $2C_1(A)+2n(3-g)+2k-2\iota(\x)$.

    For any component $\x=(X_{(g_1)},\cdots,X_{(g_k)})$, there are $k$
evaluation maps (cf. (3.5)) $$
e_i:\overline{\M}_{g,k}(X,J,A,\x)\rightarrow X_{(g_i)},
\hspace{4mm} i=1,\cdots, k. \leqno(5.1) $$
    We also have a map
    $$p: \overline{\M}_{g,k}(X,J,A,\x)\rightarrow \overline{\M}_{g,k}$$
    where $p$ contracts the unstable components of the domain to obtain
    a stable Riemann surface in $\overline{\M}_{g,k}$.
    For any set of
cohomology classes $\alpha_i\in
H^{*-2\iota_{(g_i)}}(X_{(g_i)};\Q)\subset H^*_{orb}(X;\Q)$,
$i=1,\cdots,k$ and $K\in H^*(\overline{\M}_{g,k}, \Q)$, the
orbifold Gromov-Witten invariant is defined as the pairing
 $$ \Psi^{X,J}_{(g,k,A,\x)}(K; \O_{l_1}(\alpha_1), \cdots,
\O_{l_k}(\alpha_k))=p^*K\cup \prod_{i=1}^k c_1(L_i)^{l_i}e^*_i
\alpha_i[\overline{\M}_{g,k}(X,J,A,\x)]^{vir}. $$ where $L_i$ is
the line bundle generated by the cotangent space of the $i$-th
marked point.
    When $l_i<0$, we define it to be zero.
    The same argument as in the smooth case yields

    \vskip 0.1in
    \noindent
    {\bf Proposition 3.4.1: }{\it
    \begin{enumerate}
    \item $ \Psi^{X,J}_{(g,k,A,\x)}(\O_{l_1}(K; \alpha_1), \cdots,
\O_{l_k}(\alpha_k))=0$ unless $deg K+ \sum_i
(deg_{orb}(\alpha_i)+l_i)=2C_1(A)+2(n-3)(1-g)+2k.$, where
$deg_{orb}(\alpha_i)$ is the orbifold degree of $\alpha_i$
obtained after degree shifting.
    \item $\Psi^{X,J}_{(g,k,A,\x)}(K; \O_{l_1}(\alpha_1), \cdots,
\O_{l_k}(\alpha_k))$ is independent of the choice of $J$ and hence
is an invariant of symplectic structures.
\end{enumerate}}
    We can drop $J$ from the notation.

    \vskip 0.1in
    In the smooth case, Gromov-Witten invariants satisfy a set of
    axioms which serves as the guide line to construct general
    Gromov-Witten invariants. The deformation invariants axiom was first proposed by
    Ruan \cite{Ru}. The others came from physics and
    were first proposed by Witten in \cite{W}. The same
    argument as that of Fukaya-Ono also yields the same axioms for
    orbifold Gromov-Witten invariants except for a slightly restrictive divisor axiom.
We shall list them here and
    refer for the proof to \cite{FO}. Here, we list the axioms for the
    invariant without descendants only. The invariant with
    descendants ($l_i>0$) satisfies slightly modified axioms which is the same as
    in the smooth case. We leave it to the readers.

    \vskip 0.1in
    \noindent
    {\bf Theorem 3.4.2: }{\it $\Psi^{X}_{(g,k,A,\x)}(K; \alpha_1, \cdots,
\alpha_k)$ satisfies following axioms:
    \begin{description}
    \item[Deformation Invariance Axiom] $\Psi^{X}_{(g,k,A,\x)}(K; \alpha_1, \cdots,
\alpha_k)$ is independent of the smooth deformation of symplectic
structures.

    There is a map $\pi: \overline{\M}_{g,k+1}\rightarrow \overline{\M}_{g,k}$
    for $(g,k)\neq (0,2),(1,0)$. This map yields two axioms.

    \item[Point Axiom] For any $\alpha _1, \cdots , \alpha _{k}$ in $H^*_{orb}(X, \Q)$,
we have $$\Psi ^X_{(A,g,k+1)}(K; 1, \alpha _1, \cdots,\alpha_k)~=~
\Psi ^X_{(A,g,k)}( \pi_*K; \alpha _1, \cdots,\alpha _{k})$$

   \item[Divisor Axiom] Let $\beta\in H^2(X, \Q)$ (not $H^2_{orb}(X, \Q)$).
     $$\Psi ^X_{(A,g,k+1)}(\pi^*K; \beta, \alpha_1, \cdots,\alpha _k)~=~\beta (A) \Psi
^X_{(A,g,k)}(K; \alpha_1, \cdots,\alpha _k).$$
      There are two maps
      $$\theta_1: \overline{\M}_{g_1, k_1+1}\times
      \overline{\M}_{g_2, k_2+1}\rightarrow
      \overline{\M}_{g_1+g_2, k_1+k_2}, $$
      $$\theta_2: \overline{\M}_{g, k+2}\rightarrow
      \overline{\M}_{g+1, k}$$
      by joining the two extra marked points. $\theta_1, \theta_2$
      gives two axioms collectively called the splitting axiom. Let
      $\Delta\subset \widetilde{\Sigma X}\times \widetilde{\Sigma
      X}$ be the graph of involution map $I$. Using K\"{u}nneth
      formula, we can express its Poincar\'{e} dual as
      $$\Delta^*=\sum_{a,b}\eta^{ab}\beta_a\otimes \beta_b,$$
      where $\beta_a$ is a basis of $H^*_{orb}(X, \Q)$. Note that
      $\eta^{ab}=<\beta_a, \beta_b>_{orb}$.

    \item[Splitting Axiom I] Let $K_i\in H^*(\overline{\M}_{g_i,
    k_i+1}, \Q)$. Then,
    $$\begin{array}{rl}
&\Psi ^X_{(A,g_1+g_2,k_1+k_2)}((\theta _{1})_{*}(K_1\cup
K_2)\{\alpha _i\})\\ =(-1)^{deg(K_2)\sum^{k_1}_{i=1} deg
(\alpha_i)} ~& \sum \limits _{A=A_1+A_2} \sum \limits_{a,b} \Psi
^X_{(A_1,g_1,k_1+1)}(K_1;\{\alpha _{i}\}_{i\le k}, \beta _a) \eta
^{ab} \Psi ^Y_{(A_2,g_2,k_2+1)}(K_2;\beta _b, \{\alpha _{j}\}_{j>
k}) \\
\end{array}$$
    \item[Splitting Axiom II] Let $K\in H^*(\overline{\M}_{g,
    k+2}, \Q)$.
$$ \Psi ^X_{(A,g+1,k)}((\theta_2)_*(K);\alpha _1,\cdots, \alpha
_k) =\sum _{a,b} \Psi ^X_{(A,g-1,k+2)}( K;\alpha _1,\cdots, \alpha
_k, \beta _a,\beta _b) \eta ^{ab}.$$
    \end{description}}
    \vskip 0.1in

    It is well-known that Gromov-Witten invariants satisfying the previous
axioms yield an associative quantum multiplication over $H^*_{orb}
(X, \bigwedge_{\omega})$ where $\bigwedge_{\omega}$ is the Novikov
ring. Define $$\Psi^X(\alpha, \beta, \gamma)=\sum_A
\Psi^X_{(A,0,3)}(1; \alpha, \beta, \gamma) q^A.$$
    Then the quantum product $\alpha\times_{Q} \beta$ is defined by the
relation
    $$<\alpha\times_Q\beta, \gamma>_{orb}=\Psi^Y(\alpha, \beta, \gamma).$$
This is so called the "small" quantum product. One can also define
the big quantum product which depends on a parameter $w\in
H^*_{orb}(X, \C)$. The definition is identical to the smooth case.
We omit the details.
    \vskip 0.1in
\noindent
{\bf Theorem 3.4.3: }{\it Quantum product (small or big) is associative.}
\vskip 0.1in

\section{Appendix: An Introduction to Orbifolds}

In this section we give some relevant background material on
orbifolds. Some of it are scattered in the literature or folklore
theorems; some are completely new (e.g. the notion of good map) or
have no precise formulation available in the literature. It seems
to us that there is no good reference on orbifold theory
available, starting from the basics, so it is worth while to write
up such a self-contained exposition in which almost everything has
a detailed proof. We also provide many examples to assist in
understanding the material. Our formulation is along the lines of
[FO] in contrast to the more traditional ones in [S] or [K1].

\subsection{Basic notions of orbifold}

Primary examples of orbifolds are quotient spaces of smooth
manifolds by a smooth finite group action. Here we consider that
the quotient space is uniformized (or modeled) by the manifold
with the finite group action. Hence a notion of smoothness for the
quotient space is inherited from the manifold through those
objects which are invariant under the group action. We require
that any element of the group either acts trivially or has
fixed-point set  of codimension at least two. This is the case,
for example, when the action is orientation-preserving. This
requirement has a consequence that the non-fixed-point set is
locally connected.

\vspace{2mm}

Let $U$ be a connected topological space, $V$ be a connected
n-dimensional smooth manifold and $G$ be a finite group acting on
$V$ smoothly. {\it An n-dimensional uniformizing system} of $U$ is
a triple $(V,G,\pi)$, where $\pi: V \rightarrow U$ is a continuous
map inducing a homeomorphism between $V/G$ and $U$. Two
uniformizing systems $(V_i,G_i,\pi_i)$, $i=1,2$, are {\it
isomorphic} if there is a diffeomorphism $\phi: V_1\rightarrow
V_2$ and an isomorphism $\lambda: G_1\rightarrow G_2$ such that
$\phi$ is $\lambda$-equivariant, and $\pi_2\circ\phi=\pi_1$. It is
easily seen that if $(\phi,\lambda)$ is an automorphism of
$(V,G,\pi)$, then there is a  $g\in G$ such that $\phi(x)=g\cdot
x$ and $\lambda(a)=g\cdot a\cdot g^{-1}$ for any $x\in V$ and
$a\in G$. $g$ is unique iff the action of $G$ on $V$ is effective.
We use $ker \ G$ to denote the subgroup of $G$ acting trivially.

 Let
$i:U^\prime\hookrightarrow U$ be a connected open subset of $U$,
and $(V^\prime,G^\prime, \pi^\prime)$ be a uniformizing system of
$U^\prime$. We say that $(V^\prime,G^\prime,\pi^\prime)$ is
induced from  $(V,G,\pi)$ if there is a monomorphism $\lambda:
G^\prime\rightarrow G$ and a $\lambda$-equivariant open embedding
$\phi: V^\prime\rightarrow V$ such that $\lambda$ induces an
isomorphism from $ker \ G'$ to $ker \ G$ and $i
\circ\pi^\prime=\pi\circ\phi$. We follow Satake [S] and call
$(\phi,\lambda): (V^\prime,G^\prime,\pi^\prime)\rightarrow
(V,G,\pi)$ an {\it injection}. Two injections $(\phi_i,\lambda_i):
(V^\prime_i,G^\prime_i,\pi^\prime_i)\rightarrow (V,G,\pi)$,
$i=1,2$, are {\it isomorphic} if there is an isomorphism
$(\psi,\tau)$ between $(V^\prime_1,G^\prime_1,\pi^\prime_1)$ and
$(V^\prime_2,G^\prime_2,\pi^\prime_2)$, and an automorphism
$(\bar{\psi}, \bar{\tau})$ of $(V,G,\pi)$ such that
$(\bar{\psi},\bar{\tau})\circ
(\phi_1,\lambda_1)=(\phi_2,\lambda_2)\circ (\psi,\tau)$.

\vspace{2mm}

\noindent{\bf Lemma 4.1.1: }{\it Let $(V,G,\pi)$ be a uniformizing
system of $U$. For any connected open subset $U^\prime$ of $U$,
$(V,G,\pi)$ induces a unique isomorphism class of uniformizing
systems of $U^\prime$.}

\vspace{2mm}

\noindent{\bf Proof:}

{\it Existence}: Consider the preimage $\pi^{-1}(U^\prime)$ in
$V$. $G$ acts as permutations on the set of connected components
of $\pi^{-1}(U^\prime)$. Let $V^\prime$ be one of the connected
components of $\pi^{-1}(U^\prime)$, $G^\prime$ be the subgroup of
$G$ which fixs the component $V^\prime$ and
$\pi^\prime=\pi|_{V^\prime}$. Then
$(V^\prime,G^\prime,\pi^\prime)$ is an induced uniformizing system
of $U^\prime$.

{\it Uniqueness}: First of all, different choices of the connected
components of $\pi^{-1}(U^\prime)$ induce isomorphic uniformizing
systems. Secondly, let $(V^\prime_1,G^\prime_1,\pi^\prime_1)$ be
any induced uniformizing system of $U^\prime$ and $(\psi,\tau)$ be
the injection into $(V,G,\pi)$. We will show that $(\psi,\tau)$
induces an isomorphism between
$(V^\prime_1,G^\prime_1,\pi^\prime_1)$ and the induced
uniformizing system given by a connected component of
$\pi^{-1}(U^\prime)$. Suppose $\psi( V^\prime_1)$ lies in the
connected component $V^\prime$. We can show that
$\psi(V^\prime_1)$ is closed in $V^\prime$. Let
$\psi(x_n)\rightarrow y_0$ in $V^\prime$, $x_n\in V^\prime_1$,
then there exists a $z_0\in V^\prime_1$ such that
$\pi^\prime_1(z_0)=\pi(y_0)$, and $z_n\in V^\prime_1$ such that
$z_n\rightarrow z_0$,
 $\pi^\prime_1(z_n)=\pi(\psi(x_n))=\pi^\prime_1(x_n)$. So  there exist
$a_n\in G^\prime_1$ such that $a_n(z_n)= x_n $. Since $G^\prime_1$
is finite, it follows that for large $n$, $a_n=a$ is a constant.
So $x_n\rightarrow a(z_0)$ in $V^\prime_1$ and $y_0=\psi(a(z_0))$,
i.e., $\psi(V^\prime_1)$ is closed in $V^\prime$.
$\psi(V^\prime_1)$ is also open in $V^\prime$. So $\psi$ induces a
diffeomorphism between $V^\prime_1$ and $V^\prime $. From this we
can easily see that $(V^\prime_1,G^\prime_1,\pi^\prime_1)$ and
$(V^\prime, G^\prime,\pi^\prime)$ are isomorphic. \hfill $\Box$

\vspace{2mm}

Let $U$ be a connected and locally connected topological space. For any point
$p\in U$, we can define the {\it germ} of uniformizing systems at $p$ in the
following sense. Let $(V_1, G_1,\pi_1)$ and $(V_2,G_2,\pi_2)$ be uniformizing
systems of neighborhoods $U_1$ and $U_2$ of $p$. We say that
$(V_1, G_1,\pi_1)$ and $(V_2,G_2,\pi_2)$
 are {\it equivalent} at $p$ if they induce
isomorphic uniformizing systems for a neighborhood $U_3$ of $p$.

\vspace{2mm}

\noindent{\bf Definition 4.1.2: }{\it Let $X$ be a Hausdorff,
second countable topological space. An {\it n-dimensional orbifold
structure} on $X$ is given by the following data: for any point
$p\in X$, there is a neighborhood $U_p$ and an n-dimensional
uniformizing system $(V_p, G_p, \pi_p)$ of $U_p$ such that for any
point $q\in U_p$, $(V_p, G_p, \pi_p)$ and $(V_q,G_q,\pi_q)$ are
equivalent at $q$ (i.e., define the same germ at $q$). The {\it
germ} of orbifold structures on $X$ is defined in the following
sense: two orbifold structures $\{(V_p,G_p,\pi_p): p\in X\}$ and
$\{(V^\prime_p,G^\prime_p,\pi^\prime_p): p\in X\}$ are {\it
equivalent} if for any $p\in X$, $(V_p,G_p,\pi_p)$ and
$(V^\prime_p,G^\prime_p,\pi^\prime_p)$ are equivalent at $p$. With
a given germ of orbifold structure on it, $X$ is called a {\it
orbifold}. We call each $U_p$ a {\it uniformized neighborhood } of
$p$, and $(V_p,G_p,\pi_p)$ a {\it chart} at $p$. An open subset
$U$ of $X$ is called an {\it uniformized open set} if it is
uniformized by $(V,G,\pi)$ such that for each $p\in U$,
$(V,G,\pi)$ defines the same germ as $(V_p,G_p,\pi_p)$ at $p$.
Later on we will show and therefore always assume that for each
$p\in X$, $(V_p,G_p,\pi_p)$ can be taken such that $V_p$ is a ball
in $\R^n$,  and $\pi_p$ sends the origin of $\R^n$ to $p$. A point
$p\in X$ is called {\it regular or smooth} if $G_p$ is trivial;
otherwise, it is called {\it singular}. The set of smooth points
is denoted by $X_{reg}$, and the set of singular points is denoted
by $\Sigma X$. If $G_p$ acts effectively, we call $X$ a reduced
orbifold.}

\vspace{2mm}
    It is clear that every orbifold $X$ induces a reduced orbifold
    $X_{red}$     by redefining the local group. We call $X_{red}$
    the {\em reduced associate of $X$}.

\noindent{\bf Remark 4.1.3: }{\it There is a notion of {\it
orbifold with boundary}, in which we allow the uniformizing
systems to be smooth manifolds with boundary, with a finite group
action preserving the boundary. If $X$ is an orbifold with
boundary, then it is easily seen that the boundary $\partial X$
inherits an orbifold structure from $X$ and becomes an orbifold.}
\vskip 0.1in

\noindent{\bf Example 4.1.4:}{\it Let's consider the 2-dimensional
sphere $S^2$. Let $D_s$, $D_n$ be open disc neighborhoods of the
south pole and the north pole such that $S^2=D_s\cup D_n$. Let
$D_s$ be uniformized by $(\tilde{D}_s,\Z_2,\pi_s)$, and $D_n$ be
uniformized by $(\tilde{D}_n,\Z_3,\pi_n)$ where $\Z_2$, $\Z_3$ act
on $\tilde{D}_s$ and $\tilde{D}_n$ by rotations. For any point in
$S^2$ other than the south pole and the north pole, we take a
chart at it induced by either $(\tilde{D}_s,\Z_2,\pi_s)$ or
$(\tilde{D}_n,\Z_3,\pi_n)$. It is easily seen that this defines a
2-dimensional orbifold structure on $S^2$. Note that as an open
subset of both $D_s$  and $D_n$, $D_s\cap D_n$ has non-isomorphic
induced uniformizing systems from $(\tilde{D}_s,\Z_2,\pi_s)$ and
$(\tilde{D}_n,\Z_3,\pi_n)$, although they define the same germ at
each point in $D_s\cap D_n$. This also shows that although both
$D_s$ and $D_n$ are uniformized, their union $S^2$ cannot be
uniformized, therefore is not a global quotient.} \vskip 0.1in

The notion of orbifold was first introduced by Satake in \cite{S}, where a
different name, $V$-manifold, was used. In [S], an orbifold structure on a
topological space $X$ is given by an open cover $\U$ of $X$
satisfying the following conditions:
\begin{itemize}
\item {(4.1.1a)} Each element $U$ in $\U$ is uniformized, say by $(V,G,\pi)$.
\item {(4.1.1b)} If $U^\prime\subset U$, then there is a collection of
injections $(V^\prime,G^\prime,\pi^\prime)\rightarrow (V,G,\pi)$.
\item {(4.1.1c)} For any point $p\in U_1\cap U_2$, $U_1,U_2\in\U$,
there is a $U_3\in\U$ such that $p\in U_3\subset U_1\cap U_2$.
\end{itemize}
One can show that our definition is equivalent to Satake's. On the
one hand, it is easy to see that an open cover of $X$ satisfying
$(4.1.1a-c)$ gives rise to an orbifold structure on $X$ in the
sense of Definition 4.1.2. (We will call such a cover of an
orbifold $X$ a {\it compatible cover} if it gives rise to the same
germ of orbifold structures on $X$.) On the other hand, given an
orbifold $X$ as in Definition 4.1.2, we can construct a compatible
cover $\U$ of $X$ as follows: Take a locally finite refinement of
$\{U_p\}$, denoted by $\{U_i\}$. Each $q$ in $X$ is contained in
finitely many, say $\alpha(q)$, of $U_i$'s, so there is a
connected open neighborhood $W_q$ such that $W_q$ is contained in
these $\alpha(q)$ $U_i$'s and inherits a unique isomorphism class
of uniformizing systems from them. Each connected open subset of
$W_q$ containing $q$ inherits a uniformizing system from $W_q$. We
call $\alpha(q)$ the order of $W_q$ or any of its open subsets
containing $q$. Then we take $\U$ to be the collection of all of
the $W_q$'s and its connected open subsets containing $q$. We need
to show that for any $U_1$ and $U_2$ in $\U$ such that $U_1\subset
U_2$, there is an injection between the uniformizing systems
$(V_1,G_1,\pi_1)\rightarrow (V_2,G_2,\pi_2)$. But this is easily
seen from the fact that the order of $U_1$ is always greater than
or equal to the order of $U_2$.

\vspace{2mm}

Next we consider a class of continuous maps between two orbifolds
which carry an additional structure of differentiability with
respect to the orbifold structures. Let $U$ be uniformized by
$(V,G,\pi)$ and $U^\prime$ by $(V^\prime,G^\prime,\pi^\prime)$,
and $f: U\rightarrow U^\prime$ be a continuous map. A {\it $C^l$
lifting, $0\leq l\leq\infty$}, of $f$ is a $C^l$ map
$\tilde{f}:V\rightarrow V^\prime$ such that
$\pi^\prime\circ\tilde{f}=f\circ\pi$, and for any $g\in G$, there
is $g^\prime \in G^\prime$ satisfying
$g^\prime\cdot\tilde{f}(x)=\tilde{f}(g\cdot x)$ for any $x\in V$.
Two liftings $\tilde{f}_i: (V_i,G_i,\pi_i)\rightarrow
(V^\prime_i,G^\prime_i,\pi^\prime_i)$, $i=1,2$, are {\it
isomorphic} if there exist isomorphisms
$(\phi,\tau):(V_1,G_1,\pi_1)\rightarrow (V_2,G_2,\pi_2)$ and
$(\phi^\prime,\tau^\prime):(V_1^\prime,G_1^\prime,\pi_1^\prime)\rightarrow
(V_2^\prime,G_2^\prime,\pi_2^\prime)$ such that $\phi^\prime\circ
\tilde{f}_1= \tilde{f}_2\circ \phi$. Let $p\in U$, for any
uniformized neighborhood $U_p$ of $p$ and uniformized neighborhood
$U_{f(p)}$ of $f(p)$ such that $f(U_p)\subset U_{f(p)}$; a lifting
$\tilde{f}$ of $f$ will induce a lifting $\tilde{f}_p$ for
$f|_{U_{p}}:U_p\rightarrow U_{f(p)}$ as follows: For any injection
$(\phi,\tau):(V_p,G_p,\pi_p)\rightarrow (V,G,\pi)$, consider the
map $\tilde{f}\circ\phi: V_p\rightarrow V^\prime$, and observe
that $\pi^\prime\circ\tilde{f}\circ\phi(V_p)\subset U_{f(p)}$
implies $\tilde{f}\circ\phi(V_p)\subset
(\pi^\prime)^{-1}(U_{f(p)})$. Therefore there is an injection
$(\phi^\prime,\tau^\prime):(V_{f(p)}, G_{f(p)},\pi_{f(p)})
\rightarrow (V^\prime, G^\prime, \pi^\prime)$ such that
$\tilde{f}\circ\phi(V_p)\subset \phi^\prime(V_{f(p)})$. We define
$\tilde{f}_p=(\phi^\prime)^{-1}\circ\tilde{f}\circ\phi$. In this
way we obtain a lifting $\tilde{f}_p: (V_p,G_p,\pi_p)\rightarrow
(V_{f(p)}, G_{f(p)},\pi_{f(p)})$ for $f|_{U_{p}}:U_p\rightarrow
U_{f(p)}$. We can verify that different choices give isomorphic
liftings. We define the {\it germ} of liftings as follows: two
liftings are {\it equivalent at $p$} if they induce isomorphic
liftings on a smaller neighborhood of $p$.

Now consider orbifolds $X$ and $X^\prime$ and a continuous map
$f:X\rightarrow X^\prime$. A {\it lifting} of $f$ consists of the
following data: for any point $p \in X$, there exist charts
$(V_p,G_p,\pi_p)$ at $p$ and $(V_{f(p)},G_{f(p)}, \pi_{f(p)})$ at
$f(p)$ and a lifting $\tilde{f}_p$ of $f_{\pi_p(V_p)}:
\pi_p(V_p)\rightarrow \pi_{f(p)}(V_{f(p)})$ such that for any
$q\in \pi_p(V_p)$, $\tilde{f}_p$ and $\tilde{f}_q$ induce the same
germ of liftings of $f$ at $q$. We can define the {\it germ} of
liftings in the sense that two liftings of $f$ $\{\tilde{f}_{p,i}:
(V_{p,i},G_{p,i},\pi_{p,i}) \rightarrow
(V_{f(p),i},G_{f(p),i},\pi_{f(p),i}): p\in X\}$, $i=1,2$, are {\it
equivalent} if for each $p\in X$, $\tilde{f}_{p,i}, i=1,2$, induce
the same germ of liftings of $f$ at $p$.

    \vskip 0.1in
\noindent{\bf Definition 4.1.5: }{\it A {\it $C^l$ map} ($0\leq
l\leq\infty$) between orbifolds $X$ and $X^\prime$ is a germ of
$C^l$ liftings of a continuous map between $X$ and $X^\prime$. We
denote by $\tilde{f}$ a $C^l$ map which is a germ of liftings of a
continuous map $f$.

A sequence of $C^l$ maps $\tilde{f}_n$ is said to converge to a
$C^l$ map $\tilde{f}_0$ in the $C^l$ topology if there exists a
sequence of liftings $\tilde{f}_{p,n}: (V_p,G_p,\pi_p)\rightarrow
(V_{f_n(p)},G_{f_n(p)},\pi_{f_n(p)})$ defining the germs
$\tilde{f}_n$ such that for each $p\in X$, there exists a chart
$(V_{f_0(p)},G_{f_0(p)},\pi_{f_0(p)})$ and an integer $n(p)>0$
with the following property: for each $n\geq n(p)$, there is an
injection $(\psi_{p,n},\tau_{p,n}):
(V_{f_n(p)},G_{f_n(p)},\pi_{f_n(p)})\rightarrow
(V_{f_0(p)},G_{f_0(p)},\pi_{f_0(p)})$ such that
$\psi_{p,n}\circ\tilde{f}_{p,n}$ converges in $C^l$ to
$\tilde{f}_{p,0}$ which defines the germ $\tilde{f}_0$.}\vskip
0.1in

\noindent{\bf Example 4.1.6a: }{\it The real line $\R$ as a smooth
manifold is trivially an orbifold. A $C^l$ map from an orbifold
$X$ to $\R$ is called a {\it $C^l$ function} on $X$. The set of
all $C^l$ functions on $X$ is denoted by $C^l(X)$. A $C^l$
function is essentially a continuous function which locally can be
lifted to an invariant $C^l$ function on a local chart. On the
other hand, a $C^l$ map from $\R$ (or an interval $I$) into $X$ is
called a {\it $C^l$ path} in $X$. In subsection 4.1.2, we will
discuss geodesics on a Riemannian orbifold and define the
exponential map.}\vskip 0.1in

\noindent{\bf Example 4.1.6b: }{\it Let $X=\R\times \C$, and be
given an orbifold structure by $(\R\times\C,\Z_4, \pi)$ where
$\Z_4$ acts only on the factor $\C$ by multiplication of
$\sqrt{-1}$. Define $C^1$ maps $\tilde{f}_1:\R\rightarrow
(\R\times\C,\Z_4,\pi)$ by $t\rightarrow (t,t^2)$ and
$\tilde{f}_2:\R\rightarrow (\R\times\C,\Z_4,\pi)$ by $t\rightarrow
(t,t^2)$ for $t\leq 0$ and $(t,\sqrt{-1}t^2)$ for $t\geq 0$. Then
$\tilde{f}_1,\tilde{f}_2$ induce the same continuous map
$f:\R\rightarrow X$, but they are {\it not} isomorphic as $C^1$
maps.}\vskip 0.1in

Next we describe the notion of orbifold bundle which corresponds
to the notion of smooth vector bundle. We begin with local
uniformizing systems for orbibundles. Given a uniformized
topological space $U$ and a topological space $E$ with a
surjective continuous map $pr:E\rightarrow U$, a {\it uniformizing
system of rank $k$ orbifold bundle} for $E$ over $U$ consists of
the following data:

\begin{itemize}
\item A uniformizing system $(V,G,\pi)$ of $U$.
\item A uniformizing system $(V\times\R^k,G,\tilde{\pi})$ for $E$.
The action of $G$ on $V\times\R^k$ is an extension of the action
of $G$ on $V$ given by $g(x,v)=(gx,\rho(x,g)v)$, where
$\rho:V\times G\rightarrow Aut(\R^k)$ is a smooth map satisfying:
$$ \rho(gx,h)\circ\rho(x,g)=\rho(x,h\circ g), \hspace{3mm} g,h\in
G, x\in V. $$
\item The natural projection map $\tilde{pr}:V\times\R^k\rightarrow V$
satisfies $\pi\circ\tilde{pr}=pr\circ\tilde{\pi}$.
\end{itemize}

We can similarly define {\it isomorphisms} between uniformizing
systems of orbifold bundle for $E$ over $U$. The only additional
requirement is that the diffeomorphisms between $V\times\R^k$ are
linear on each fiber of $\tilde{pr}:V\times\R^k\rightarrow V$.
Moreover, for each connected open subset $U^\prime$ of $U$, we can
similarly prove that there is a unique isomorphism class of
induced uniformizing systems of orbifold bundle for
$E^\prime=pr^{-1}(U^\prime)$ over $U^\prime$. The {\it germ} of
uniformizing systems of orbifold bundle at a point $p\in U$ can
also be  similarly defined.

    \vskip 0.1in
\noindent{\bf Definition 4.1.7: }{\it Let $X$ be an orbifold and
$E$ be a topological space with a surjective continuous map
$pr:E\rightarrow X$. A {\it rank $k$ orbifold bundle  structure}
on $E$ over $X$ consists of the following data: For each point
$p\in X$, there is a uniformized neighborhood $U_p$ and a
uniformizing system of rank $k$ orbifold bundle for $pr^{-1}(U_p)$
over $U_p$ such that for any $q\in U_p$, the uniformizing systems
of orbifold bundle over $U_p$ and $U_q$ define the same germ at
$q$. The {\it germ} of rank $k$ orbifold bundle structures on $E$
over $X$ can be similarly defined. The topological space $E$ with
a given germ of orbifold bundle structures becomes an orbifold and
is called an {\it orbifold bundle} over $X$. Each chart
$(V_p\times\R^k,G_p,\tilde{\pi}_p)$ is called a {\it local
trivialization} of $E$. At each point $p\in X$, the fiber
$E_p=pr^{-1}(p)$ is isomorphic to $\R^k/G_p$. It contains a linear
subspace $E^p$ of fixed points of $G_p$. Two orbibundles
$pr_1:E_1\rightarrow X$ and $pr_2:E_2\rightarrow X$ are {\it
isomorphic} if there is a $C^\infty$ map
$\tilde{\psi}:E_1\rightarrow E_2$ given by $\tilde{\psi}_p:
(V_{1,p}\times\R^k,G_{1,p},\tilde{\pi}_{1,p})\rightarrow
(V_{2,p}\times\R^k,G_{2,p},\tilde{\pi}_{2,p})$ which induces an
isomorphism between $(V_{1,p},G_{1,p},\pi_{1,p})$ and
$(V_{2,p},G_{2,p},\pi_{2,p})$, and is a linear isomorphism between
the fibers of $\tilde{pr}_{1,p}$ and $\tilde{pr}_{2,p}$. By
replacing $\R^k$ with $\C^k$, we have the definition of  complex
orbifold bundle}. \vskip 0.1in

\noindent{\bf Remark 4.1.8a: }{\it There is a notion of orbifold
bundle over an orbifold with boundary. One can easily verify that
if $pr:E\rightarrow X$ is an orbifold bundle over an orbifold with
boundary $X$, then the restriction to the boundary $\partial X$,
$E_{\partial X}=pr^{-1}(\partial X)$, is an orbifold bundle over
$\partial X$.}\vskip 0.1in

\noindent{\bf Remark 4.1.8b: }{\it One can define an orbifold
bundle
 with fiber a general space in the same vein.}
    \vskip 0.1in

A $C^l$ map $\tilde{s}$ from $X$ to an orbifold bundle
$pr:E\rightarrow X$ is called a {\it $C^l$ section} if locally
$\tilde{s}$ is given by $\tilde{s}_p: V_p\rightarrow
V_p\times\R^k$ where $\tilde{s}_p$ is $G_p$-equivariant and
$\tilde{pr}\circ\tilde{s}_p=Id$ on $V_p$. We observe that

\begin{itemize}
\item For each point $p$, $s(p)$ lies in $E^p$, the linear subspace of fixed
points of $G_p$.
\item The space of all $C^l$ sections of $E$, denoted by $C^l(E)$, has a
structure of vector space over $\R$ (or $\C$) as well as
a $C^l(X)$-module structure.
\item The $C^l$ sections $\tilde{s}$ are in $1:1$ correspondence with the
underlying continuous maps $s$.
\end{itemize}

Orbibundles are more conveniently described by transition maps
(see [S]). More precisely, an orbifold bundle over an orbifold $X$
can be constructed from the following data: A compatible cover
$\U$ of $X$ such that for any injection
$i:(V^\prime,G^\prime,\pi^\prime)\rightarrow (V,G,\pi)$, there is
a smooth map $g_i:V^\prime\rightarrow Aut(\R^k)$ giving an open
embedding $V^\prime\times\R^k\rightarrow V\times\R^k$ by
$(x,v)\rightarrow (i(x),g_i(x)v)$, and for any composition of
injections $j\circ i$, we have $$ g_{j\circ i}(x)=g_j(i(x))\circ
g_i(x), \forall x\in V. \leqno (4.1.2) $$ Two collections of maps
$g^{(1)}$ and $g^{(2)}$ define isomorphic orbibundles if there are
maps $\delta_V:V\rightarrow Aut(\R^k)$ such that for any injection
$i: (V^\prime,G^\prime,\pi^\prime) \rightarrow (V,G,\pi)$, we have
$$ g^{(2)}_i(x)=\delta_V(i(x))\circ g^{(1)}_i(x)\circ
(\delta_{V^\prime}(x))^{-1}, \forall x\in V^\prime. \leqno (4.1.3)
$$ Since $(4.1.2)$ behaves naturally under constructions of vector
spaces such as tensor product, exterior product, etc. we can
define these constructions for orbibundles.

    \vskip 0.1in
\noindent{\bf Example 4.1.10: }{\it For an orbifold $X$, the
tangent bundle $TX$ can be constructed because the differential of
any injection satisfies $(4.1.2)$. Likewise, we define the
cotangent bundle $T^\ast X$, the bundles of exterior power or
tensor product. The $C^\infty$ sections of these bundles give us
vector fields, differential forms or tensor fields on $X$. There
exists a de Rham cohomology theory for orbifolds, which is
isomorphic to the singular cohomology theory of the underlying
topological space. Observe also that if $\omega$ is a differential
form on $X^\prime$ and $\tilde{f}:X\rightarrow X^\prime$ is a
$C^\infty$ map, then there is a pull-back form
$\tilde{f}^\ast\omega$ on $X$.} \vskip 0.1in

Let $U$ be an open subset of an orbifold $X$ with an orbifold
structure $\{(V_p,G_p,\pi_p):p\in X\}$; then
$\{(V^\prime_p,G^\prime_p,\pi^\prime_p): p\in U\}$ is an orbifold
structure on $U$, where $(V^\prime_p,G^\prime_p,\pi^\prime_p)$ is
a uniformizing system of $\pi_p(V_p)\cap U$ induced from
$(V_p,G_p,\pi_p)$. Likewise, let $pr:E\rightarrow X$ be an
orbifold bundle and $U$ an open subset of $X$; then
$pr:E_{U}=pr^{-1}(U)\rightarrow U$ inherits a unique germ of
orbifold bundle structures from $E$, called the {\it restriction
of $E$ over $U$}. When $U$ is a uniformized open set in $X$, say
uniformized by $(V,G,\pi)$, then there is a smooth vector bundle
$E_V$ over $V$ with a smooth action of $G$ such that
$(E_V,G,\tilde{\pi})$ uniformizes $E_U$. This is seen as follows:
We first take a compatible cover $\U$ of $U$, fine enough so that
the preimage under $\pi$ is a compatible cover of $V$. Let $E_U$
be given by a set of transition maps with respect to $\U$
satisfying $(4.1.2)$; then the pull-backs under $\pi$ form a set
of transition maps with respect to $\pi^{-1}(\U)$ with an action
of $G$ by permutations, also satisfying $(4.1.2)$, so that it
defines a smooth vector bundle over $V$ with a compatible smooth
action of $G$. Any $C^l$ section of $E$ on $X$ restricts to a
$C^l$ section of $E_U$ on $U$, and when $U$ is a uniformized open
set by $(V,G,\pi)$, it lifts to a $G$-equivariant $C^l$ section of
$E_V$ on $V$.

We end this subsection with a result which is analogous to the homotopy
invariance of
vector bundles. Let $I=[0,1]$. Then if $X$ is an orbifold, $X\times I$ is
an orbifold with boundary, with $\partial(X\times I)=X\times\{0\}\cup
X\times\{1\}$.

    \vskip 0.1in
\noindent{\bf Proposition 4.1.11: }{\it Let $pr: E\rightarrow
X\times I$ be an orbifold bundle over $X\times I$. Then there is a
$C^\infty$ map $\tilde{\Psi}: E\rightarrow E$ covering the
$C^\infty$ map $\tilde{\psi}: X\times I \rightarrow X\times I$,
given by $\tilde{\psi}(x,t)=(x,1)$, such that each local lifting
of $\tilde{\Psi}$ is an isomorphism on each fiber of $\tilde{pr}$.
In particular, the two orbifold bundles $E_{X\times\{0\}}$ and
$E_{X\times\{1\}}$ are isomorphic.}\vskip 0.1in

\noindent{\bf Proof:}
First of all, by the compactness of $I$ and the local triviality of $E$, for
any $p\in X$, there exists a finite open covering $\{I_i\}$ of $I$ and a set
of neighborhoods $\{U_{p,i}\}$ of $p$ such that $E_{U_{p,i}
\times I_i}$ is trivial. Take a neighborhood $U_p\subset \cap_i U_{p,i}$,
then $E_{U_p\times I_i}$ is trivial for all $I_i$'s.
Next we show that if $I_i\cap I_j\neq \emptyset$, we can construct a
trivialization of $E_{U_p\times (I_i\cup I_j)}$, which is a uniformizing
system of $E_{U_p\times (I_i\cup I_j)}$, and
a trivialization of $E_{U_p\times I}$ successively.
Without loss of generality, we assume that $I_i=[0,b)$ and $I_j=(a,1]$,
for some $a<b$.

Let $(V_p\times [0,b)\times \R^k,G_p,\tilde{\pi}_b)$ be a
trivialization of $E_{U_p\times [0,b)}$ and $(V_p\times
(a,1]\times \R^k,G_p,\tilde{\pi}_a)$ a trivialization of
$E_{U_p\times (a,1]}$. We let $(V_p\times (a,b)\times
\R^k,G_p,\tilde{\pi}_{a,b})$ be the trivialization of
$E_{U_p\times (a,b)}$ induced from $(V_p\times [0,b)\times
\R^k,G_p, \tilde{\pi}_b)$, and $\psi: (V_p\times (a,b)\times
\R^k,G_p,\tilde{\pi}_{a,b})\rightarrow (V_p\times (a,1]\times
\R^k,G_p, \tilde{\pi}_a)$ be the injection covering  the injection
$(V_p\times (a,b),G_p,\pi_{a,b})\rightarrow (V_p\times
(a,1],G_p,\pi_a)$, which is identical on the $V_p$ factor. We
remark that one can assume that the action of $G_p$ on $\R^k$ is
independent of $t\in I$.  Now we can define a trivialization
$(V_p\times [0,1]\times \R^k,G_p,\tilde{\pi})$ of $E_{U_p\times
[0,1]}$ as follows: pick a $c$ satisfying $a<c<b$ and a $C^\infty$
diffeomorphism $\beta: (a,1]\rightarrow (a,b)$ which is the
identity on $(a,c)$; then we define $$
\tilde{\pi}(x,t,v)=\left\{\begin{array}{cc} \tilde{\pi}_b(x,t,v) &
\mbox{if   } t\in [0,b)\\ \tilde{\pi}_a\circ
(id\times\beta^{-1}\times id) \circ \psi^{-1}\circ
(id\times\beta\times id)(x,t,v) & \mbox{if  } t\in (a,1].
\end{array} \right.
$$

Now we pick a smaller neighborhood $U_p^\prime$ of $p$ such that
$\overline{U_p^\prime}\subset U_p$, and a $G_p$-equivariant
$C^\infty$ function $\theta (x,t)$ from $V_p\times [0,1]$ to
$[0,1]$ such that $\theta(x,t)=t$ for $x$ outside a neighborhood
of $V_p^\prime\subset V_p$ and $\theta(x,t)=1$ for $x\in
V_p^\prime$. We define $\tilde{f}_p: E_{U_p\times I}\rightarrow
E_{U_p\times I}$ to  be the $C^\infty$ map given by the
$G_p$-equivariant map $(x,t,v)\rightarrow (x,\theta(x,t),v)$,
which is the identity outside a neighborhood of $V_p^\prime\times
[0,1]\times \R^k$. We can extend it to a $C^\infty$ map
$E\rightarrow E$, still denoted by $\tilde{f}_p$. Now we cover $X$
by a locally finite covering of $U_{p_i}^\prime\subset U_{p_i}$,
and define $\tilde{\Psi}: E\rightarrow E$ to be the product of the
$\tilde{f}_{p_i}$'s. \hfill $\Box$

\subsection{Riemannian metric and exponential map}

A {\it partition of unity} on an orbifold $X$ is a collection $\{\rho_i:
i\in \Lambda\}$ of $C^\infty$ functions on $X$ such that
\begin{itemize}
\item The collection of supports $\{supp\rho_i: i\in\Lambda\}$ is locally
finite.
\item $\sum_{i\in\Lambda}\rho_i(p)=1$ for all $p\in X$, and
$\rho_i(p)\geq 0$ for all $p\in X$ and $i\in \Lambda$.
\end{itemize}

    \vskip 0.1in
\noindent{\bf Lemma 4.2.1: }{\it Any cover of $X$ consisting of
uniformized open subsets admits a partition of unity on $X$
subordinate to it.}\vskip 0.1in

\noindent{\bf Proof:}
Let $\{U_i:i\in \Lambda\}$ be a locally finite refinement of the cover.
By the paracompactness of $X$, such a refinement exists. Each $U_i$ is also
a uniformized open set, say by $(V_i,G_i,\pi_i)$, so that we can take a
$G_i$-invariant, non-negative smooth function $\tilde{f}_i$ on $V_i$ such that
$f_i$ can be extended  over $X$ by zeros and
$\{supp f_i:i\in\Lambda\}$ also covers $X$, i.e.,
$f(p)=\sum_{i\in\Lambda}f_i(p)\neq 0$ for all $p\in X$.
We define $\rho_i=f_i/f$. Then $\{\rho_i:i\in\Lambda\}$ is the required
partition of unity.
\hfill $\Box$

    \vskip 0.1in
\noindent{\bf Definition 4.2.2: }{\it  A {\it Riemannian metric}
on an orbifold $X$ is a positive, symmetric $(0,2)$ tensor field
on $X$.}\vskip 0.1in

Let $U$ be a uniformized open subset of $X$, uniformized by $(V,G,\pi)$.
Then any Riemannian metric on $X$ induces a compatible $G$-invariant
metric on $V$.

    \vskip 0.1in
\noindent{\bf Lemma 4.2.3: }{\it The space of all Riemannian
metrics on $X$ is a non-empty cone in the vector space of
$C^\infty$ sections of the orbifold bundle of the symmetric square
of $T^\ast X$.}\vskip 0.1in

\noindent{\bf Proof:} The cone structure is obvious, we only need
to show the existence. But this follows easily from the existence
of partitions of unity on $X$. \hfill $\Box$

Let $p\in X$ and $(V_p,G_p,\pi_p)$ be a local chart at $p$. Let
$G^p$ be the isotropy subgroup of $G_p$ at $p$ (which is the germ
of the groups $G_p$). Since $G_p$ takes geodesics in $V_p$ to
geodesics, there is a convex geodesic ball $B_p(r)$ of radius $r$
at $p$ which is invariant under $G^p$ and $(B_p(r),G^p,\pi_p)$ is
a uniformizing system of $\pi_p(B_p(r))$ in $\pi_p(V_p)$. Via the
exponential map, we can think of $B_p(r)$ as a ball of radius $r$
in $\R^n$ and $G^p$ acts as a subgroup of $O(n)$. We call
$(B_p(r),G^p,\pi_p)$ a {\it geodesic chart of radius} $r$ at $p$
and $\pi_p(B_p(r))$ a {\it geodesic neighborhood of radius} $r$.
An open subset $U$ of a geodesic neighborhood is called {\it
star-shaped} if the induced uniformizing system from the geodesic
chart is a star-shaped domain with respect to the origin.

    \vskip 0.1in
\noindent{\bf Lemma 4.2.4: }{\it Let $U_p$ be a geodesic
neighborhood of radius $r$ at $p$ and $(V_i,G_i,\pi_i)$, $i=1,2$,
be two uniformizing systems of $U_p$ defining the same germ at
$p$; then $(V_i,G_i,\pi_i)$, $i=1,2$, are isomorphic. As a
consequence, for any uniformized open set $U$ in $X$, uniformized
by $(V,G,\pi)$, if a geodesic neighborhood $U_p$ is contained in
$U$, then there is an injection from the corresponding geodesic
chart into $(V,G,\pi)$. The same conclusion also holds  for any
star-shaped open subset in a geodesic neighborhood.}\vskip 0.1in

\noindent{\bf Proof:} The uniformizing system $(V_i,G_i,\pi_i)$ is
isomorphic to the geodesic chart $(B_p(r),G_{p,i},\pi_{p,i})$,
$i=1,2$. Since the two uniformizing systems define the same germ
at $p$, $G_{p,1}=G_{p,2}$ and $\pi_{p,1}=\pi_{p,2}$, hence
$(V_i,G_i,\pi_i)$, $i=1,2$, are isomorphic. \hfill $\Box$

    \vskip 0.1in
\noindent{\bf Definition 4.2.5: }{\it A $C^\infty$ path
$\tilde{\gamma}:I\rightarrow X$ in a Riemannian orbifold $X$ is
called a {\it parametrized geodesic} if for each $t\in I$,
$\tilde{\gamma}_t$ is a geodesic in some chart
$(V_{\gamma(t)},G_{\gamma(t)},\pi_{\gamma(t)})$ at $\gamma(t)$.
The image $\gamma(I)$  is called a {\it geodesic} in $X$.}\vskip
0.1in

Since on a smooth manifold a geodesic is uniquely determined by the tangent
vector at one point on it, it is easily seen that the parametrized geodesic
$\tilde{\gamma}$ is in $1:1$ correspondence with the underlying continuous
map $\gamma$.

The length of a parametrized geodesic $\tilde{\gamma}$ is defined to be
$$
length(\tilde{\gamma})=\int_{I}|\dot{\gamma}(t)|dt,
$$
where $\dot{\gamma}(t)$ is the tangent vector of $\tilde{\gamma}$ at $t\in I$,
as a vector in the fiber of the tangent bundle $TX$ at $\gamma(t)$. It is
easily seen that the norm $|\dot{\gamma}(t)|$ is constant in $t$.

    \vskip 0.1in
\noindent{\bf Lemma 4.2.6: }{\it  Let $p\in X$ be any point,
$TX_p=\R^n/G_p$ ($n=\dim X$) be the fiber of the tangent bundle
$TX$ at $p$. For any vector $v\neq 0$ in $TX_p$, there is a unique
maximal parametrized geodesic $\tilde{\gamma}:
(-\epsilon_1,\epsilon_2) \rightarrow X$ such that $\gamma(0)=p$
and $\dot{\gamma}(0)=v$. When $X$ is compact,
$\epsilon_1=\epsilon_2=\infty$ for any $p\in X$.}

    \vskip 0.1in

\noindent{\bf Definition 4.2.7: }{\it $X$ is called {\it
geodesically complete} if $\epsilon_1=\epsilon_2=\infty$ for any
$p\in X$.}\vskip 0.1in

\noindent{\bf Proof:}
Parametrized geodesics satisfying $\gamma(0)=p$ and $\dot{\gamma}(0)=v$
exist locally around $p$. So we can take the union of all parametrized
geodesics satisfying the said conditions to get the maximal one.
We  only need to show that a parametrized geodesic
$\tilde{\gamma}: (a,b)\rightarrow X$
is extendable in the following sense: if there is a sequence $t_i\in (a,b)$
such that $t_i\rightarrow b$ and $\lim\gamma(t_i)$ exists in $X$, then there
is a parametrized geodesic $\tilde{\gamma}^\prime: (a,b^\prime)\rightarrow
X$, $b^\prime>b$, such that $\tilde{\gamma}=\tilde{\gamma}^\prime$ on
$(a,b)$. This is seen as follows: Let $q=\lim\gamma(t_i)$. We take a
geodesic chart $(V_q,G_q,\pi_q)$ at $q$. For large $i$,
$\gamma(t_i)$ lies in $\pi_q(V_q)$. Since $V_q$ is
geodesically convex, all of $\gamma(t)$
lies in $\pi_q(V_q)$ when $t\rightarrow b$ and $\lim\gamma(t)=q$
as $t\rightarrow b$. The extendability follows from local existence.
\hfill $\Box$

Now we are ready to define the exponential map. We assume that $X$ is
geodesically complete. We define the exponential map $exp: TX\rightarrow X$ by
$$
exp(p,v)={\gamma}_{(p,v)}(1),
$$
where $\tilde{\gamma}_{(p,v)}$ is the unique maximal parametrized geodesic
satisfying $\gamma(0)=p$ and $\dot{\gamma}(0)=v$.

    \vskip 0.1in
\noindent{\bf Proposition 4.2.8: }{\it The exponential map $exp$
is continuous, and there exists a canonically defined $C^\infty$
map, denoted by $Exp$, as a germ of $C^\infty$ liftings of
$exp$.}\vskip 0.1in

\noindent{\bf Remark 4.2.9: }{\it We will introduce a notion of
{\it good map} and {\it compatible system} in subsection 4.3. The
proof here also shows that $Exp$ is a good map with a unique
isomorphism class of compatible systems.}\vskip 0.1in

\noindent{\bf Proof:}
We first prove that $exp$ is continuous. Suppose $(p_i,v_i)\in TX$ converges to
$(p,v)$. Let $I$ be the subset of $(-\infty,\infty)$ consisting of $t$
such that $\lim\gamma_{(p_i,v_i)}(t)$ exists and equals $\gamma_{(p,v)}(t)$.
Then an open neighborhood of $0$ is in
$I$ and $I$ is both open and closed. So $I=(-\infty,\infty)$ and
$\gamma_{(p,v)}(1)=\lim\gamma_{(p_i,v_i)}(1)$. Hence $exp$ is continuous.

The germ of $C^\infty$ liftings $Exp$ is characterized as follows:
Each geodesic chart of $TX$ at $(p,v)$ is a product of a geodesic chart $B_p$
of $X$ at $p$ with a ball $B_v$ in $\R^n$ centered at a preimage of $v$
in $\R^n$. The group action is given by the subgroup
$G_{p,v}$ of $G_p$ which fixes $v$. We denote such a chart of $TX$ at $(p,v)$
by $(V_{p,v},G_{p,v},\pi_{p,v})$. Let $q=exp(p,v)$.
Then any local lifting $\tilde{f}:(V_{p,v},G_{p,v},\pi_{p,v})
\rightarrow (V_q,G_q,\pi_q)$ compatible with $Exp$ has the property
that the restriction of $\tilde{f}$ to the intersection of
each 1-dim linear subspace  of $\R^n$ with $B_v$
(we will call these intersections {\it lines} in $V_{p,v}$) is a
parametrized geodesic in $V_q$.

Before we construct such local liftings of $exp$, we prove that
they are unique up to isomorphism. In other words, let
$\tilde{f}_i$, $i=1,2$, be two local liftings of $exp$ compatible
with $Exp$ from $(V_{p,v},G_{p,v},\pi_{p,v})$ to
$(V_q,G_q,\pi_q)$; then there is an automorphism $(\psi,\tau)$
(which is actually unique) of $(V_q,G_q,\pi_q)$ such that
$\psi\circ\tilde{f}_1=\tilde{f}_2$. This is seen as follows: For
each line $l$ in $V_{p,v}$, $\tilde{f}_1(l)$ and $\tilde{f}_2(l)$
are two parametrized geodesics in $(V_q,G_q,\pi_q)$ which are
mapped to the same geodesic in $X$ under $\pi_q$. So there is an
element $g(l)\in G_q$ such that
$g(l)\circ\tilde{f}_1(l)=\tilde{f}_2(l)$. When $\tilde{f}_1(l)$
and $\tilde{f}_2(l)$ do not lie entirely in the singular set in
$V_q$, $g(l)$ is unique. This is true when $l$ is not a singular
line in $V_{p,v}$. The set of singular lines in $V_{p,v}$ is of
codimension at least two in the space of all lines in $V_{p,v}$,
so that there is an element $g\in G_q$ such that
$g\circ\tilde{f}_1(l)=\tilde{f}_2(l)$ for any non-singular line
$l$. By continuity, $g\circ\tilde{f}_1=\tilde{f}_2$ on $V_{p,v}$,
which proves the claim.

Finally, we show that these local liftings do exist. we proceed as follows:
It suffices to show that
for each $p_0\in X$, there is a geodesic chart $B_{p_0}$ such that for any
$p\in B_{p_0}$, these local liftings of $exp$ exist in a neighborhood of
$(p,v)$ for any $v\in TX_{p}$. This is certainly true for
those $v$'s with small norms. We will show that
for any $r>0$ there is a $\epsilon_r>0$
such that if it is true for any $(p,v)$ where $p \in B_{p_0}$
and $|v|<r$, then it is true for any $(p, v)$ where $|v|<r+\epsilon_r$.
This goes as follows:
For any $(p,v)$, $|v|=r$, let $q=exp(p,v)$, and let
$(V_{p,v},G_{p,v},\pi_{p,v})$ be a geodesic chart at
$(p,v)$ and $(V_q,G_q,\pi_q)$ be any chart at $q$ such that
$exp\circ \pi_{p,v}(V_{p,v})\subset \pi_q (V_q)$. For any $(p^\prime,v^\prime)
\in V_{p,v}$ such that $|v^\prime|<r$, the local liftings defined in a
local chart at $(p^\prime,v^\prime)$ will induce local liftings via injections
into $(V_{p,v},G_{p,v},\pi_{p,v})$ on the images of injections in $V_{p,v}$.
By the same reason as stated  in the previous paragraph, these induced local
liftings fit together to give a local lifting defined on the subset $V$ of
$V_{p,v}$ such that $\pi_{p,v}(V)$ consists of those $(p^\prime,v^\prime)$
with $|v^\prime|<r$. We simply extend this lifting over the whole $V_{p,v}$
by parametrized geodesics. Now the existence of $\epsilon_r$ for each $r$
follows from the precompactness of the set $\{(p,v):p\in B_{p_0},|v|=r\}$
in $TX$ so that for each such a $(p,v)$, there is a chart at $exp(p,v)$
such that a geodesic ball of radius $2\epsilon_r$ centered at a preimage
of $exp(p,v)$ is contained in the chart, for some $\epsilon_r>0$.
This proves the existence of $Exp$.
\hfill $\Box$.

    \vskip 0.1in
\noindent{\bf Corollary 4.2.10: }{\it Suppose $X$ is geodesically
complete. If for $p, q\in X$, there is a parametrized geodesic
connecting $p$ and $q$, then there is a parametrized geodesic of
minimal length connecting $p$ and $q$.}\vskip 0.1in

\noindent{\bf Proof:}
The assumption that there is a parametrized geodesic connecting $p$ and $q$
means that there is a $v\in TX_p$ such that $exp(p,v)=q$. Let $v_0$ be in
$TX_p$ such that $|v_0|=\inf|v|$
where $\inf$ is taken over all $v\in TX_p$ such that $exp(p,v)=q$, then
by continuity of $exp$, $exp(p,v_0)=q$ and the parametrized geodesic
$\tilde{\gamma}_{p,v_0}$ is of minimal length.
\hfill $\Box$

    \vskip 0.1in
\noindent{\bf Definition 4.2.11: }{\it A parametrized geodesic
connecting $p$ and $q$ which is of minimal length is called a {\it
minimal geodesic}. An open subset $U$ of $X$ is called {\it
geodesically convex} if for any $p,q\in U$, every minimal geodesic
connecting $p$ and $q$ lies entirely in $U$. It is easily seen
that a geodesic neighborhood of sufficiently small radius is
geodesically convex, and the intersection of two geodesically
convex sets is geodesically convex.}\vskip 0.1in

\subsection{Differential geometry of orbifold vector bundle}

In this subsection, we discuss the Chern-Weil theory for orbifold
bundles. Because a general orbifold vector bundle may not have any
local section, the usual constructions of connection and curvature
do not work for such an orbifold vector bundle. We call an
orbifold vector bundle {\em good} if the local group of the base
and total space have the same kernel. Equivalently, an orbifold
vector bundle  $E\rightarrow X$ is good iff $E_{red}\rightarrow
X_{red}$ is an orbifold vector bundle of $X_{red}$. In this case,
the local sections of $E$ are in one-to-one correspondence to
local sections of $E_{red}$. The differential geometry of $E$ is
identical to the differential geometry of $E_{red}$. Hence, we can
assume that $X$ is reduced for this purpose. Examples of good
orbifold vector bundles include the tangent and cotangent bundles
as well as their tensor and exterior products.  Throughout this
subsection, we assume that our orbifold vector bundle is good.

Let $U$ be a connected topological space uniformized by
$(V,G,\pi)$, with an orbifold bundle of rank $k$, $pr:
E\rightarrow U$ uniformized by $(V\times \R^k,G,\tilde{\pi})$. We
consider $G$-equivariant connections $\nabla$ on $V\times \R^k$,
i.e., for any smooth section $u$ and vector field $v$ on $V$ and
any element $g\in G$, we have $g(\nabla_v u)=\nabla_{gv}gu$. Let
$(V_i\times \R^k,G_i,\tilde{\pi}_i)$, $i=1,2$, be two isomorphic
uniformizing systems of $pr: E\rightarrow U$, with
$G_i$-equivariant connections $\nabla_i$ on it. We say $\nabla_1$
and $\nabla_2$ are {\it isomorphic} if for an isomorphism $\psi$
from $(V_1\times \R^k,G_1,\tilde{\pi}_1)$ to $(V_2\times
\R^k,G_2,\tilde{\pi}_2)$, the equation $\psi^\ast\nabla_2
=\nabla_1$ holds. Then one can easily check that given a connected
open subset $U^\prime$ of $U$, each isomorphism class of $\nabla$
over $U$ induces a unique isomorphism class of connections
$\nabla^\prime$ over $U^\prime$. One can define the {\it germ} at
point $p$ in the sense that $\nabla_1$ and $\nabla_2$ are {\it
equivalent} at $p$ if they induce isomorphic connections over a
neighborhood of $p$.

    \vskip 0.1in
\noindent{\bf Definition 4.3.1: }{\it Let $pr: E\rightarrow X$ be
an orbifold bundle with an orbifold bundle structure
$\V=\{(V_p\times\R^k,G_p,\tilde{\pi}_p): p\in X\}$. A {\it
connection} $\nabla$ on $pr: E\rightarrow X$  is a collection of
connections $\{\nabla_p: p\in X\}$, with each $\nabla_p$ being a
$G_p$-equivariant connection on
$(V_p\times\R^k,G_p,\tilde{\pi}_p)$ such that for any $q\in
U_p=\pi_p(V_p)$, $\nabla_p$ and $\nabla_q$ are equivalent at $q$.
Two connections $\nabla_i$ in the reference of orbifold bundle
structures $\V_i$, $i=1,2$, are {\it equivalent} if they induce
isomorphic connections over a neighborhood of each point $p\in X$.
We will call each equivalence class a {\it connection } on $pr:
E\rightarrow X$. Note that if $u$ is a $C^\infty$ section of $E$
and $v$ is a $C^\infty$ section of $TX$, then $\nabla_v u$ is a
$C^\infty$ section of $E$.}\vskip 0.1in

Given two connections $\nabla_1$ and $\nabla_2$ on $pr:
E\rightarrow X$, say in the reference to orbifold bundle
structures $\V_1$ and $\V_2$. For each $p\in X$, there is a
uniformized neighborhood $U_p$ such that $\V_1$ and $\V_2$ induce
isomorphic uniformizing systems of $E_{U_p}=pr^{-1}(U_p)$ over
$U_p$, say $(V_p\times\R^k,G_p,\tilde{\pi}_p)$.  $\nabla_1$ and
$\nabla_2$ then induce two $G_p$-equivariant connections
$\nabla_{1,p}$ and $\nabla_{2,p}$ on
$(V_p\times\R^k,G_p,\tilde{\pi}_p)$, whose difference is a
$G_p$-equivariant, $End(\R^k)$-valued smooth $1$-form on $V_p$.
From this consideration it is easily seen that the difference of
$\nabla_1$ and $\nabla_2$ is a $C^\infty$ section of the orbifold
bundle $T^\ast X\otimes End(E)$. On the other hand, a connection
on $pr: E\rightarrow X$ added to a $C^\infty$ section of $T^\ast
X\otimes End(E)$ gives rise to another connection on $E$.

    \vskip 0.1in
\noindent{\bf Lemma 4.3.2: }{\it The set of all connections on
$pr: E\rightarrow X$ forms a non-empty affine space modelled on
the space of $C^\infty$ sections of the orbifold bundle $T^\ast
X\otimes End(E)$.}\vskip 0.1in

\noindent{\bf Proof:} We only need to show that the set of all
connections is non-empty. Suppose $X$ is covered by a collection
of uniformized open sets $\U=\{U_i:i \in \Lambda\}$, each $U_i$
uniformized by $(V_i,G_i,\pi_i)$, and let $E_i$ be a smooth vector
bundle over $V_i$ which uniformizes $pr^{-1}(U_i)$. Moreover we
assume there is a $C^\infty$ partition of unity
$\{\rho_i:i\in\Lambda\}$ subordinate to the cover $\U$. Let
$\tilde{\rho}_i$ be the corresponding $G_i$-invariant function on
$V_i$.  Let $\nabla_i$ be a $G_i$-equivariant connection on $E_i$.
Now suppose $\V=\{(V_p\times\R^k,G_p,\tilde{\pi}_p): p\in X\}$ is
an orbifold bundle structure of $pr: E\rightarrow X$; we will
construct a connection $\nabla$ on $E$, that is, a collection of
connections $\nabla_p$ on $V_p\times\R^k$ which is
$G_p$-equivariant and for any $q\in \pi_p(V_p)$, $\nabla_p$ and
$\nabla_q$ are equivalent at $q$. We first define the connection
$\nabla_p$. This amounts  for any smooth section $u$ of
$V_p\times\R^k$ and any vector field $v$ on $V_p$, to define
$(\nabla_p)_v u$, at any point $\tilde{x}\in V_p$. Let
$x=\pi_p(\tilde{x})$. We take a small geodesic neighborhood $U_x$
of $x$ such that if $x$ is in $U_i$ for any $i\in\Lambda$,
$U_x\subset U_i$, and if $U_x\cap U_j\neq \emptyset$ but $U_x$ is
not contained in $U_j$ for some $j\in\Lambda$, then
$\overline{U_x\cap U_j}\cap supp \ \rho_j=\emptyset$. Let
$(V_x\times\R^k,G_x,\tilde{\pi}_x)$ be a uniformizing system of
$pr^{-1}(U_x)$ and $\psi:
(V_x\times\R^k,G_x,\tilde{\pi}_x)\rightarrow
(V_p\times\R^k,G_p,\tilde{\pi}_p)$ an injection covering
$\tilde{x}$.  For any inclusion $U_x\subset U_i$, we pick an
injection $\psi_i: (V_x\times\R^k,G_x,\tilde{\pi}_x)\rightarrow
E_i$. We define $$ (\nabla_p)_v u|_{\tilde{x}}
=[\psi_\ast\sum_{i\in\Lambda}(\psi_i^{-1})_\ast(\tilde{\rho}_i(\nabla_i)_
{(\psi_i)_\ast(\psi^{-1})_\ast(v)}(\psi_i)_\ast(\psi^{-1})_\ast(u))]_
{\tilde{x}}. $$ One can verify that $(\nabla_p)_v u|_{\tilde{x}}$
is independent of the choices of $U_x$, $\psi$ and $\psi_i$, and
is smooth in $\tilde{x}$, so that we have a well-defined
connection on $V_p\times\R^k$, and such defined $\nabla_p$ is also
$G_p$-equivariant, and for any $q\in \pi_p(V_p)$, $\nabla_p$ and
$\nabla_q$ are equivalent at $q$. Hence the lemma. \hfill $\Box$

    \vskip 0.1in
\noindent{\bf Remark 4.3.3: }{\it One can put a Riemannian metric
on an orbifold bundle, or a Hermitian metric on a complex orbifold
bundle, and define the corresponding notion of $O(k)$-connection
or $U(k)$-connection, and prove that the  set of
$O(k)$-connections or $U(k)$-connections is a non-empty affine
space.}\vskip 0.1in

For each connection $\nabla$ on $pr:E\rightarrow X$, we can define
the curvature $F(\nabla)$ of $\nabla$ as a $C^\infty$ section of
the orbifold bundle $\Lambda^2(X)\otimes End(E)$. This is done as
follows: let $\nabla$ be a connection for orbifold bundle
structure $\V=\{(V_p\times\R^k,G_p,\tilde{\pi}_p): p\in X\}$, by
$\{\nabla_p:p\in X\}$. Then the curvature $F(\nabla_p)$ of
$\nabla_p$ is a $G_p$-equivariant smooth 2-form on $V_p$ with
values in $End(\R^k)$. The collection of $\{F(\nabla_p):p\in X\}$
defines a $C^\infty$ section of $\Lambda^2(X)\otimes End(E)$. We
denote it by $F(\nabla)$. The Chern-Weil construction applied to
each $F(\nabla_p)$ yields the following

    \vskip 0.1in
\noindent{\bf Proposition 4.3.4: }{\it To each invariant
polynomial, there is a characteristic class in the following
sense. To each (isomorphism class of) orbifold bundle
$pr:E\rightarrow X$, there is an associated cohomology class in
the deRham cohomology group $H^\ast_{de}(X)$ of $X$ as follows:
each connection on $E$ is assigned  a closed differential form on
$X$, and for two different connections on $E$, the associated
closed forms differ by an exact form on $X$. As a consequence, the
Chern classes, the Pontrjagin classes and the Euler class are all
well-defined in the category of orbibundles.}\vskip 0.1in

\noindent{\bf Remark 4.3.5: }{\it The same thing holds for the
more general notion of orbifold bundles.}\vskip 0.1in

\noindent{\bf Example 4.3.6a: }{\it  Let $(\Sigma,j, z_1,...,z_k)$
be a smooth complex curve with a set of finitely many distinct
points $(z_1,...,z_k)$. An orbifold structure on $(\Sigma,j)$ can
be given by a collection of uniformizing systems $\tilde{D}_i$ at
each $z_i$, i.e, a collection of disk neighborhoods $D_i$ of $z_i$
and branched covering maps $\tilde{D}_i\rightarrow D_i$ by
$z\rightarrow z^{m_i}$. The numbers $m_i$ are called the {\it
multiplicities} at $z_i$. The orbifold $(\Sigma,j)$ thus defined
is an example of a complex orbifold. Its cotangent bundle, or the
canonical bundle, denoted by $K_\Sigma$, is an example of a
holomorphic line orbifold bundle. The first Chern number
$c_1(K_\Sigma)([\Sigma])$ is equal to $$ 2g_\Sigma-2+\sum_{i=1}^k
(1-\frac{1}{m_i}), $$ where $g_\Sigma$ is the genus of $\Sigma$.
This example shows that the characteristic classes of an orbifold
bundle are rational classes in general.}\vskip 0.1in

\noindent{\bf Example 4.3.6b: }{\it  Let $\tilde{E}$ be a
holomorphic line bundle over a smooth complex curve
$\tilde{\Sigma}$, with a finite group $G$ acting holomorphically
on it. Let $\Sigma=\tilde{\Sigma}/G$ and
$\pi:\tilde{\Sigma}\rightarrow \Sigma$ be the natural branched
covering map. Then $E=\tilde{E}/G$ is a holomorphic line orbifold
bundle over $\Sigma$. We have the following relation amongst the
Chern classes of $\tilde{E}$ and $E$:
$
\pi^\ast(c_1(E))=c_1(\tilde{E}),
$
or equivalently, in terms of Chern numbers:
$c_1(E)([\Sigma])=\frac{1}{|G|}c_1(\tilde{E})([\tilde{\Sigma}])$.}\vskip
0.1in

\subsection{Pull-back orbifold bundles and good maps}

Let $pr: E\rightarrow Y$ be a vector bundle over a topological
space $Y$. Then for any continuous map $f: X\rightarrow Y$ from  a
topological space $X$, the pull-back vector bundle $f^\ast E$ over
$X$ is well-defined. However, this is no longer the case for
orbifold bundles. Let $pr:E\rightarrow X^\prime$ be an orbifold
bundle over $X^\prime$, and $\tilde{f}:X\rightarrow X^\prime$ a
$C^\infty$ map. By a {\it pull-back orbifold bundle of $E$ over
$X$ via $\tilde{f}$ } we mean an orbifold bundle
$\pi:E^\ast\rightarrow X$ together with a $C^\infty$ map
$\bar{f}:E^\ast\rightarrow E$ such that each local lifting of
$\bar{f}$ is an isomorphism restricted to each fiber, and
$\bar{f}$ covers the $C^\infty$ map $\tilde{f}$ between the bases.

Let $\tilde{f}:X\rightarrow X^\prime$ be a $C^\infty$ map between
orbifolds $X$ and $X^\prime$ whose underlying continuous map is
denoted by $f$. Suppose there is a compatible cover $\U$ of $X$,
and a collection of open subsets $\U^\prime$ of $X^\prime$
satisfying $(4.1.1a-c)$ and the following condition: There is a
1:1 correspondence between elements of $\U$ and $\U^\prime$, say
$U\leftrightarrow U^\prime$, such that $f(U)\subset U^\prime$, and
an inclusion $U_2\subset U_1$ implies an inclusion
$U_2^\prime\subset U_1^\prime$. Moreover, there is a collection of
local $C^\infty$ liftings  $\{\tilde{f}_{UU^\prime}\}$ of $f$,
where $\tilde{f}_{UU^\prime}:(V,G,\pi)\rightarrow
(V^\prime,G^\prime,\pi^\prime)$ satisfies the following condition:
each  injection $i:(V_2,G_2,\pi_2)\rightarrow (V_1,G_1,\pi_1)$ is
assigned with an injection
$\lambda(i):(V_2^\prime,G_2^\prime,\pi^\prime_2)\rightarrow
(V^\prime_1,G_1^\prime,\pi_1^\prime)$ such that
$\tilde{f}_{U_1U^\prime_1}\circ
i=\lambda(i)\circ\tilde{f}_{U_2U^\prime_2}$, and for any
composition of injections $j\circ i$, the following compatibility
condition holds: $$ \lambda(j\circ i)=\lambda(j)\circ \lambda(i).
\leqno (4.4.1) $$ Observe that when the injection
$i:(V,G,\pi)\rightarrow (V,G,\pi)$ is an automorphism of
$(V,G,\pi)$, the assignment of $\lambda(i)
:(V^\prime,G^\prime,\pi^\prime)\rightarrow
(V^\prime,G^\prime,\pi^\prime)$ to $i$ satisfying $(4.4.1)$ is
equivalent to a homomorphism $\lambda_{UU^\prime}: G\rightarrow
G^\prime$. We call $\lambda_{UU^\prime}: G\rightarrow G^\prime$
the {\it group homomorphism} of
$\{\tilde{f}_{UU^\prime},\lambda\}$ on $U$. Such a collection of
maps clearly defines a $C^\infty$ lifting of the continuous map
$f$. If it is in the same germ of $\tilde{f}$, we call
$\{\tilde{f}_{UU^\prime},\lambda\}$ a {\it compatible system} of
$\tilde{f}$.

    \vskip 0.1in
\noindent{\bf Definition 4.4.1: }{\it A $C^\infty$ map is called
{\it good} if it admits a compatible system.}\vskip 0.1in

\noindent{\bf Example 4.4.2a: }{\it Not every $C^\infty$ map is
good, as shown in the following example: consider an effective
linear representation of a finite group $(\R^n,G)$. Let $H^g$ be
the linear subspace of fixed points of an element $1\neq g\in G$.
Then the centralizer $C_G(g)$ of $g$ in $G$ acts on $H^g$, and
$(H^g,C_G(g)/K_g)$ is an effective linear representation, where
$K_g\subset C_G(g)$ is the isotropy subgroup of $H^g$. Suppose
$H^g\neq \{0\}$ and there is no homomorphism $\lambda:
C_G(g)/K_g\rightarrow C_G(g)$ such that $\pi\circ\lambda=id$,
where $\pi: C_G(g)\rightarrow C_G(g)/K_g$ is the projection; then
the inclusion $H^g\rightarrow \R^n$ is a $C^\infty$ map but not a
good one.}\vskip 0.1in

\noindent{\bf Example 4.4.2b: }{\it  There can be essentially
different compatible systems of the same $C^\infty$ map, as shown
in the following example: Let $X=\C\times\C/G$ where $G=\Z_2
\oplus\Z_2$ acts on $\C\times\C$ in the standard way. For the
$C^l$ map $\tilde{f}:(\C,\Z_2)\rightarrow (\C\times\C,G)$ defined
by the inclusion of $\C\times\{0\}$, there are two compatible
systems $(\tilde{f},\lambda_i): (\C,\Z_2)\rightarrow
(\C\times\C,G)$, $i=1,2$, for $\lambda_1(1)=(1,0)$ and
$\lambda_2(1)=(1,1)$, which are apparently different.}\vskip 0.1in

\noindent{\bf Lemma 4.4.3: }{\it Let $pr:E\rightarrow X^\prime$ be
an orbifold bundle over $X^\prime$. For any $C^\infty$ compatible
system $\xi=\{\tilde{f}_{UU^\prime},\lambda\}$ of a good
$C^\infty$ map $\tilde{f}:X\rightarrow X^\prime$, there is a
canonically constructed pull-back orbifold bundle of $E$ via
$\tilde{f}$: an orbifold bundle $pr:E^\ast_\xi\rightarrow X$
together with a $C^\infty$ map $\bar{f}_\xi:E^\ast_\xi \rightarrow
E$ covering $\tilde{f}$. Suppose that $E$ is good. $E^*_{\xi}$ is
obviously good as well.  Let $c$ be a universal characteristic
class defined by the Chern-Weil construction; then
$\tilde{f}^\ast(c(E))=c(E^\ast_\xi)$.}\vskip 0.1in

\noindent{\bf Proof:} Without loss of generality, we assume that
$E_{U^\prime}$ over $U^\prime$ is uniformized by $(V^\prime\times
\R^k,G^\prime,\tilde{\pi}^\prime)$. Then we have a collection of
pull-back bundles $\tilde{f}_{UU^\prime}^\ast (V^\prime \times
\R^k)$ over $V$ which has the form of $V\times\R^k$. Let
$\{g^\prime\}$ be a collection of transition maps of $E$ with
respect to $\U^\prime$, we define a set of transition maps $\{g\}$
on $X$ with respect to the cover $\U$ by pull-backs, i.e., we set
$g_i=g^\prime_{\lambda(i)}\circ\tilde{f}_{U_2U^\prime_2}$ for any
injection $i:(V_2,G_2,\pi_2)\rightarrow (V_1,G_1,\pi_1)$, where
$\lambda(i):(V_2^\prime,G_2^\prime,\pi_2^\prime)\rightarrow
(V_1^\prime,G_1^\prime,\pi_1^\prime)$ is the injection assigned to
$i$. Then the compatibility condition $(4.4.1)$ implies that the
set of maps $\{g\}$ satisfies the equation $(4.1.2)$, which
defines an orbifold bundle over $X$. We denote it by
$pr:E^\ast\rightarrow X$. The existence of a $C^\infty$ map
$\bar{f}:E^\ast\rightarrow E$ is obvious from the construction. On
the other hand, for any connection $\nabla$ on $E$, there is a
pull-back connection $\tilde{f}^\ast(\nabla)$ on $E^\ast$, so that
the equation $\tilde{f}^\ast(c(E))=c(E^\ast)$ holds for any
universal characteristic class $c$ defined by the Chern-Weil
construction. \hfill $\Box$

    \vskip 0.1in
\noindent{\bf Definition 4.4.4: }{\it  Two compatible systems
$\xi_i$ for $i=1,2$ of a $C^\infty$ map $\tilde{f}:X\rightarrow
X^\prime$ are isomorphic if for any orbifold bundle $E$ over
$X^\prime$, there is an isomorphism $\psi$ between the
corresponding pull-back orbibundles $E^\ast_i$ with
$\bar{f}_i:E_i^\ast\rightarrow E$, $i=1,2$, such that
$\bar{f}_1=\psi\circ\bar{f}_2$.}\vskip 0.1in

For technical reasons, we need to construct some standard compatible systems
from a given compatible system. We first give the definition.

    \vskip 0.1in
\noindent{\bf Definition 4.4.5: }{\it A compatible system
$\{\tilde{f}_{UU^\prime},\lambda\}$ is called a {\it geodesic
compatible system} if $\U$ and $\U^\prime$ satisfy the following
conditions: There is a countable dense subset $\Lambda$ of $X$ and
for each point $p$ in $\Lambda$, there is a sequence of geodesic
neighborhoods $U_{p,i}$ of radius $r_i\searrow 0$.  Moreover,
$\U=\{U_{p,i}:p\in\Lambda, i\in\Z_+\}$, and each element
$U_{p,i}^\prime$ in $\U^\prime$ is a star-shaped neighborhood of
$p$ whose diameter goes to zero as $i\rightarrow \infty$.}\vskip
0.1in

The technique in the proof of the following lemma will be frequently used.

    \vskip 0.1in
\noindent{\bf Lemma 4.4.6: }{\it  For any compatible system
$\{\tilde{f}_{UU^\prime},\lambda\}$, there is a geodesic
compatible system $\{\tilde{f}_{WW^\prime},\lambda_1\}$ such that
for each pair $(W,W^\prime)$, there is $(U,U^\prime)$ such that
$\tilde{f}_{WW^\prime}$ is induced from $\tilde{f}_{UU^\prime}$
for some injections $W\rightarrow U$ and $W^\prime\rightarrow
U^\prime$, and $\lambda_1$ is also induced from $\lambda$ in a
sense that will be clear from the proof.}\vskip 0.1in

Such a geodesic compatible system is called an {\it induced}
compatible system.

\noindent{\bf Proof:}
We take a countable subset $\U_0=\{U_i:i\in\Z_+\}$ of $\U$ which covers $X$.
Then we construct a countable dense subset $\Lambda$ of $X$ which is given as
the union of the following sequence of subsets:
$$
\Lambda_1,\Lambda_2,\cdots,
$$
where $\Lambda_k$ is at most countable and locally finite, and each point $p$
in $\Lambda_k$ comes with a geodesic neighborhood $W_p$ of radius less than
$\frac{1}{k}$ such that if $U_i\in\U_0$ is the first element containing $p$,
then $W_p\subset U_i$ and $f(W_p)$ is contained in a geodesic neighborhood
$W_{f(p)}\subset U^\prime_i$, and  any point in $X$
is contained in $W_p$ for some point $p$ in $\Lambda_k$ for each $k$. This
is done as follows: For each integer $k\geq 1$, for any point $p\in X$,
we take a geodesic neighborhood $W_p$ of radius less than $\frac{1}{k}$ and
if $U_i\in\U_0$ is the first element containing $p$, then $W_p\subset U_i$
and $f(W_p)$ is contained in a geodesic neighborhood
$W_{f(p)}\subset U^\prime_i$. Then we obtain a cover $\{W_p:p\in X\}$ of $X$.
We take a subcover of it which is locally finite, and we obtain the set
$\Lambda_k$.

Now for each point $p\in \Lambda\cap U_1$, we take a sequence of geodesic
neighborhoods $\{U_{p,i}:i\in\Z_+\}$ of radius $r_i\searrow 0$ with
$U_{p,1}=W_p$, and a sequence of corresponding geodesic neighborhoods
$\{U_{p,i}^\prime:i\in\Z_+\}$ of radius $r_i^\prime\searrow 0$ with
$U_{p,1}^\prime=W_{f(p)}$ and $f(U_{p,i})\subset U^\prime_{p,i}$. It is
easily seen that $\W_1=\{U_{p,i}:i\in\Z_+,p\in\Lambda\cap U_1\}$ forms a
compatible cover of $U_1$. However, an inclusion $U_{q,j}\subset U_{p,i}$
may not imply an inclusion $U_{q,j}^\prime\subset U_{p,i}^\prime$. We
overcome this problem by changing $U_{q,j}^\prime$ to
$U_{q,j}^\prime\cap U_{p,i}^\prime$, which is a geodesically convex and
star-shaped neighborhood of $f(q)$. This can be done because there are
only finitely many $U_{p,i}$'s satisfying $U_{q,j}\subset U_{p,i}$.
Now for each $U_{p,i}\in\W_1$, the corresponding $U_{p,i}^\prime$ is a
geodesically convex and star-shaped neighborhood of $f(p)$. We let
$\W_1^\prime=\{U_{p,i}^\prime:i\in\Z_+,p\in\Lambda\cap U_1\}$.
For each inclusion $U_{p,j}\subset U_{p,i}$, there is a canonical injection
between the corresponding geodesic charts $(V_{p,j},G_p,\pi_{p,j})
\rightarrow (V_{p,i},G_p,\pi_{p,i})$ which induces identity map on the tangent
space at $p$ and identity homomorphism $G_p\rightarrow G_p$. Similar
injections exist between $(V_{p,i}^\prime, G_{f(p)},\pi_{p,i}^\prime)$'s,
the geodesic uniformising systems of $U_{p,i}^\prime$'s. We call these
canonical injections {\it geodesic injections}. On the other hand, for each
$U_{p,i}\in \W_1$, we fix an induced injection $i_{p,i}:(V_{p,i},G_p,\pi_{p,i})
\rightarrow (V_1,G_1,\pi_1)$ and an induced injection
$i_{p,i}^\prime:(V_{p,i}^\prime,G_{f(p)},\pi_{p,i}^\prime)
\rightarrow (V_1^\prime,G_1^\prime,\pi_1^\prime)$ such that
$\tilde{f}_{U_1U_1^\prime}\circ i_{p,i}(V_{p,i})\subset i_{p,i}^\prime
(V_{p,i}^\prime)$ so that we can induce a lifting $\tilde{f}_{p,i}=
(i_{p,i}^\prime)^{-1}\circ\tilde{f}_{U_1U_1^\prime}\circ i_{p,i}$ from
$(V_{p,i},G_p,\pi_{p,i})$ to $(V_{p,i}^\prime,G_{f(p)},\pi_{p,i}^\prime)$.
We can further require that $i_{p,j}=i_{p,i}\circ g_{i,j}$ and
$i_{p,j}^\prime=i_{p,i}^\prime\circ g_{i,j}^\prime$ where $g_{i,j}$ and
$g_{i,j}^\prime$ are the geodesic injections.
We call the injections $i_{p,i}$ and $i_{p,i}^\prime$ {\it home injections},
and the open set $U_1\in \U_0$ the {\it home} of each element $U_{p,i}$ in
$\W_1$.

To each injection $i: (V_{q,j},G_q,\pi_{q,j})\rightarrow
(V_{p,i},G_p,\pi_{p,i})$, we need to assign an injection $\lambda_1(i):
(V_{q,j}^\prime,G_{f(q)},\pi_{q,j}^\prime)\rightarrow
(V_{p,i}^\prime,G_{f(p)},\pi_{p,i}^\prime)$ such that
$\lambda_1(i)\circ\tilde{f}_{q,j}=\tilde{f}_{p,i}\circ i$ and the
compatibility condition
$(4.4.1)$ holds for $\lambda_1$. First observe that for each injection
$i: (V_{q,j},G_q,\pi_{q,j})\rightarrow
(V_{p,i},G_p,\pi_{p,i})$, there is a unique automorphism $g(i)\in G_1$ of
$(V_1,G_1,\pi_1)$ such that $i_{p,i}\circ i=g(i)\circ i_{q,j}$. The uniqueness
of $g$ ensures $g(j\circ i)=g(j)\circ g(i)$. The image of $g(i)$ under
the group homomorphism $\lambda_{U_1U_1^\prime}:G_1\rightarrow G_1^\prime$,
denoted by $g^\prime(i)$, determines a unique injection $i^\prime:
(V_{q,j}^\prime,G_{f(q)},\pi_{q,j}^\prime)\rightarrow
(V_{p,i}^\prime,G_{f(p)},\pi_{p,i}^\prime)$ by the equation
$i_{p,i}^\prime\circ i^\prime=g^\prime(i)\circ i_{q,j}^\prime$. We define
$\lambda_1(i)=i^\prime$, which satisfies $(4.4.1$.

Now we move to the next stage of construction. For each point $p\in\Lambda
\cap(U_2\setminus U_1)$, we similarly take a sequence of geodesic neighborhoods
$U_{p,i}$ with corresponding $U_{p,i}^\prime$ which are geodesically convex
and star-shaped. We define $\W_2=\{U_{p,i}\}$ and $\W_2^\prime
=\{U_{p,i}^\prime\}$. We can similarly fix a home injection for each element
in $\W_2$ and $\W_2^\prime$ and define local liftings $\tilde{f}_{p,i}$
and the assignment $\lambda_1$, just as what we have done for the points
in $\Lambda\cap U_1$. We call $U_2$ the {\it home} of each element of $\W_2$.

There is a possibility that for an inclusion $W_1\subset W_2$ with
different homes, here $W_1\in\W_1$ and $W_2\in\W_2$, there may not
be an inclusion $W_1^\prime\subset W_2^\prime$. We overcome this
problem again by changing $W_1^\prime$ to the intersection
$W_1^\prime\cap W_2^\prime$ which remains to be geodesically
convex and star-shaped. With this done, we take
$\W_{12}=\W_1\cup\W_2$ and
$\W_{12}^\prime=\W_1^\prime\cup\W_2^\prime$. It is easily seen
that $\W_{12}$ is a compatible cover of $U_1\cup U_2$. The last
thing we need to do in order to complete the construction of a
geodesic compatible system on $U_1\cup U_2$ is to assign to each
injection $i: (V_{q,j},G_q,\pi_{q,j})\rightarrow
(V_{p,i},G_p,\pi_{p,i})$, where $U_{q,j}$ has home $U_1$ and
$U_{p,i}$ has home $U_2$, an injection $\lambda_1(i):
(V_{q,j}^\prime,G_{f(q)},\pi_{q,j}^\prime)\rightarrow
(V_{p,i}^\prime,G_{f(p)},\pi_{p,i}^\prime)$ such that
$\lambda_1(i)\circ\tilde{f}_{q,j}=\tilde{f}_{p,i}\circ i$ and
$(4.4.1)$ holds for $\lambda_1$. Since the set of injections $i:
(V_{q,j},G_q,\pi_{q,j})\rightarrow (V_{p,i},G_p,\pi_{p,i})$ is in
1:1 correspondence with the set of injections $i_1:
(V_{q,j_1},G_q,\pi_{q,j_1})\rightarrow (V_{p,i},G_p,\pi_{p,i})$
through the geodesic injection $g_{j,j_1}$, it suffices to define
$\lambda_1(i)$ for those $U_{q,j}$ with $j$ large.  This is done
as follows: First observe that $q\in U_1\cap U_2$, so that there
is a $U\in\U$ such that $U\subset U_1\cap U_2$ and $U_{q,j}\subset
U$ by taking large $j$. For each injection $i:
(V_{q,j},G_q,\pi_{q,j})\rightarrow (V_{p,i},G_p,\pi_{p,i})$, there
are injections $i_{UU_1}:(V,G,\pi) \rightarrow (V_1,G_1\pi_1)$ and
$i_{UU_2}:(V,G,\pi)\rightarrow (V_2,G_2,\pi_2)$ such that
$i=i_{p,i}^{-1}\circ i_{UU_2}\circ i_{UU_1}^{-1} \circ i_{q,j}$.
We define $\lambda_1(i)=(i_{p,i}^\prime)^{-1}\circ
\lambda(i_{UU_2})\circ \lambda(i_{UU_1})^{-1}\circ
i_{q,j}^\prime$. The fact that $\lambda$ satisfies $(4.4.1)$
implies that such defined $\lambda_1$ is independent of the
choices of $U$ and $i_{UU_i}$'s. One can also verify that
$\lambda_1$ satisfies
$\lambda_1(i)\circ\tilde{f}_{q,j}=\tilde{f}_{p,i}\circ i$ and
$(4.4.1)$. By repeating this process, we obtain the required
geodesic compatible system. \hfill $\Box$

    \vskip 0.1in
\noindent{\bf Remark 4.4.7: }{\it From the construction it is
easily seen that given a set of {\it finitely} many compatible
systems, we can assume that the corresponding induced compatible
systems have common $\W$ and $\W^\prime$.}\vskip 0.1in

\noindent{\bf Proposition 4.4.8: }{\it Two compatible systems are
isomorphic if and only if there are corresponding induced
compatible systems $\{\tilde{f}_{1,WW^\prime},\lambda_1\}$ and
$\{\tilde{f}_{2,WW^\prime},\lambda_2\}$, and automorphisms
$\{\delta_{V^\prime}\}$ of the uniformizing system
$(V^\prime,G^\prime, \pi^\prime)$ of each element
$W^\prime\in\W^\prime$, such that
$\delta_{V^\prime}\circ\tilde{f}_{1,WW^\prime}
=\tilde{f}_{2,WW^\prime}$, and for any injection
$i:(W_2,G_2,\pi_2)\rightarrow (W_1,G_1,\pi_1)$, we have
$\lambda_2(i)=\delta_{V^\prime_1}\circ\lambda_1(i)\circ
(\delta_{V^\prime_2}) ^{-1}$.}\vskip 0.1in

\noindent{\bf Proof:} First observe that the induced compatible
system defines identical pull-back orbibundles. Secondly, if
$\{\tilde{f}_{1,WW^\prime},\lambda_1\}$,
$\{\tilde{f}_{2,WW^\prime},\lambda_2\}$, and
$\{\delta_{V^\prime}\}$ satisfy
$\delta_{V^\prime}\circ\tilde{f}_{1,WW^\prime}=\tilde{f}_{2,WW^\prime}$,
and for any injection $i:(W_2,G_2,\pi_2)\rightarrow
(W_1,G_1,\pi_1)$,
$\lambda_2(i)=\delta_{V^\prime_1}\circ\lambda_1(i)\circ
(\delta_{V^\prime_2}) ^{-1}$, then for any orbifold bundle $E$
over $X^\prime$, if $E^\ast_i$, $i=1,2$, are the corresponding
pull-back orbibundles with $\bar{f}_i:E^\ast_i \rightarrow E$, we
have an isomorphism $\delta:E^\ast_1\rightarrow E^\ast_2$
determined by $\delta_{V^\prime}$ by $(4.1.3)$, which satisfies
$\bar{f}_2\circ\delta=\bar{f}_1$. Hence the proposition in one
direction.

On the other hand, for any two isomorphic compatible systems, we
take their induced compatible systems
$\{\tilde{f}_{1,WW^\prime},\lambda_1\}$ and
$\{\tilde{f}_{2,WW^\prime},\lambda_2\}$, which has common $W$ and
$W^\prime$. Then the isomorphism between the pull-back orbibundles
gives rise to a collection of automorphisms $\delta_{V^\prime}$
such that
$\tilde{f}_{2,WW^\prime}=\delta_{V^\prime}\circ\tilde{f}_{1,WW^\prime}$
and for any injection $i:(W_2,G_2,\pi_2)\rightarrow
(W_1,G_1,\pi_1)$, $$ g_{\lambda_2(i)}\circ\tilde{f}_{2,WW^\prime}=
\delta_{V^\prime_1}\circ g_{\lambda_1(i)}\circ
(\delta_{V^\prime_2})^{-1}\circ\tilde{f}_{2,WW^\prime}, $$ where
$\{g\}$ is the set of transition maps of orbifold bundle $E$ over
$X^\prime$. We take $E=TX^\prime$, then $g_i=d i$ for any
injection $i$, so we have
$\lambda_2(i)=\delta_{V^\prime_1}\circ\lambda_1(i)\circ
(\delta_{V^\prime_2}) ^{-1}$. \hfill $\Box$

    \vskip 0.1in
\noindent{\bf Definition 4.4.9: }{\it  It is easily seen that for
each $p\in X$, a compatible system determines a group homomorphism
$G_p\rightarrow G_{f(p)}$, and an isomorphism class of compatible
systems determines a conjugacy class of homomorphisms. We call
such a conjugacy class of group homomorphisms the {\it group
homomorphism} of the said isomorphism class of compatible systems
at $p$.} \vskip 0.1in

\noindent{\bf Definition 4.4.10: }{\it A $C^\infty$ map
$\tilde{f}:X\rightarrow X^\prime$ is called {\it regular} if the
underlying continuous map $f$ has the following property:
$f^{-1}(X_{reg}^\prime)$ is an open dense and connected subset of
$X$.}\vskip 0.1in

\noindent{\bf Lemma 4.4.11: }{\it  If $\tilde{f}$ is regular, then
$\tilde{f}$ is the unique germ of $C^\infty$ liftings of $f$.
Moreover, $\tilde{f}$ is good with a unique isomorphism class of
compatible systems.} \vskip 0.1in

\noindent{\bf Proof:}
We need to show that if $\tilde{f}_1, \tilde{f}_2: (V,G,\pi)\rightarrow
(V^\prime,G^\prime,\pi^\prime)$ are two local liftings of $f:\pi(V)\rightarrow
\pi^\prime(V^\prime)$, then there is an automorphism $\psi$ of
$(V^\prime,G^\prime,\pi^\prime)$ such that $\tilde{f}_2=\psi\circ\tilde{f}_1$.
Note that if for $p\in V$ such that $f(\pi(p))$ is regular in
$\pi^\prime(V^\prime)$, there is a unique automorphism
$\psi_p$ of $(V^\prime,G^\prime,\pi^\prime)$ such that
$\tilde{f}_2(p)=\psi_p\circ\tilde{f}_1(p)$. $\psi_p$ is locally constant in
$p$. So the fact that $f^{-1}(X_{reg}^\prime)$ is an open dense and
connected subset of $X$ implies that there is a unique $\psi$ such that
$\tilde{f}_2=\psi\circ\tilde{f}_1$. This proves the uniqueness of $\tilde{f}$
as the germ of $C^\infty$ liftings of $f$.

In order to show that $\tilde{f}$ is good, we need to construct a compatible
system of $\tilde{f}$. Suppose there is a compatible cover $\U$ of $X$, and
a collection of open subsets $\U^\prime$
of $X^\prime$ satisfying $(4.1.1a-c)$ and the following condition:
there is a 1:1 correspondence between elements of
$\U$ and $\U^\prime$, say $U\leftrightarrow U^\prime$, such that
$f(U)\subset U^\prime$ and
$U_2\subset U_1$ implies $U_2^\prime\subset U_1^\prime$.  Moreover, there is
a collection of local $C^\infty$ liftings  $\{\tilde{f}_{UU^\prime}\}$ of $f$
where $\tilde{f}_{UU^\prime}:(V,G,\pi)\rightarrow
(V^\prime,G^\prime,\pi^\prime)$ satisfy the following condition:
for any injection $i:(V_2,G_2,\pi_2)\rightarrow (V_1,G_1,\pi_1)$, there is
an injection $\lambda(i):(V_2^\prime,G_2^\prime,\pi^\prime_2)\rightarrow
(V^\prime_1,G_1^\prime,\pi_1^\prime)$ such that
$\tilde{f}_{U_1U^\prime_1}\circ i=\lambda(i)\circ\tilde{f}_{U_2U^\prime_2}$.
Then the condition that $\tilde{f}$ is regular implies that each
$\lambda(i)$ is unique to $i$ so that $(4.4.1)$ automatically holds
for $\lambda$. Now we can appeal to the general construction of
induced compatible systems to show that $\tilde{f}$ is good.

For any two compatible systems of $\tilde{f}$, we have two induced compatible
systems $\{\tilde{f}_{1,WW^\prime},\lambda_1\}$ and
$\{\tilde{f}_{2,WW^\prime},\lambda_2\}$. There is a collection of
automorphisms $\{\delta_{V^\prime}\}$ of the uniformizing system
$(V^\prime,G^\prime,\pi^\prime)$ of each element $W^\prime\in\W^\prime$,
such that $\delta_{V^\prime}\circ\tilde{f}_{1,WW^\prime}
=\tilde{f}_{2,WW^\prime}$. For any injection
$i:(W_2,G_2,\pi_2)\rightarrow (W_1,G_1,\pi_1)$, we define
$\tau_2(i)=\delta_{V^\prime_1}\circ\lambda_1(i)\circ
(\delta_{V^\prime_2})^{-1}$. Then we have a compatible system
$\{\tilde{f}_{2,WW^\prime},\tau_2\}$ which is isomorphic to
$\{\tilde{f}_{1,WW^\prime},\lambda_1\}$. Now $\tilde{f}$ is regular implies
that $\tau_2=\lambda_2$. Hence the lemma.
\hfill $\Box$

    \vskip 0.1in
\noindent{\bf Remark 4.4.12a: }{\it  Similar argument shows that
the exponential map $Exp$ is a good $C^\infty$ map with a unique
isomorphism class of compatible systems. The key point here is
that each line segment in $TX$ is mapped to a parametrized
geodesic and the set of singular lines is of codimension at least
two. See the construction of $Exp$ for details.}\vskip 0.1in

\noindent{\bf Remark 4.4.12b: }{\it  Here are some examples of
good $C^\infty$ maps. Let $E$ be an orbifold bundle over $X$. By
the definition, the projection is good. Moreover,  any $C^\infty$
section $\tilde{s}$ of $E$, as a $C^\infty$ map $X\rightarrow E$,
is also good. If $X$ is reduced, both of them are regular. Let $f$
be a good map. The map $f^*E\rightarrow E$ is again a good map.
}\vskip 0.1in

\noindent{\bf Remark 4.4.12c: }{\it  Here is another class of
regular $C^\infty$ maps. A $C^\infty$ map $\tilde{f}: X\rightarrow
X^\prime$ is called a $C^\infty$ embedding if each local lifting
$\tilde{f}_p:(V_p,G_p,\pi_p)\rightarrow
(V_{f(p)},G_{f(p)},\pi_{f(p)})$ is a $\lambda_p$-equivariant
embedding for some isomorphism $\lambda_p: G_p\rightarrow
G_{f(p)}$. If $X'$ is reduced, it is easily seen that
$f^{-1}(\Sigma X^\prime) =\Sigma X$ so that $\tilde{f}$ is
regular. For a $C^\infty$ embedding $\tilde{f}:X\rightarrow
X^\prime$, the normal bundle of $f(X)$ in $X^\prime$ as an
orbifold bundle is well-defined, which is isomorphic to the
quotient orbifold bundle $(TX^\prime)^\ast/TX$.}\vskip 0.1in

\noindent{\bf Example 4.4.13: }{\it  Let $\P^1=\{[z_0,z_1]\}$ be
the 1-dimensional complex projective space. We define a $\Z_2$
action on it by $x\cdot [z_0,z_1]=[xz_0,z_1]$. Let $X=\P^1/Z_2$ be
the orbifold as quotient space. Similarly, let $\P^2
=\{[z_0,z_1,z_2]\}$ be the 2-dimensional complex projective space,
with a $Z_2\oplus Z_2$ action on it, given by $(x,y)\cdot
[z_0,z_1,z_2] =[xz_0,yz_1,z_2]$. Let $X^\prime=\P^2/(Z_2\oplus
Z_2)$ be the orbifold as quotient space. We consider two sequences
of $C^\infty$ maps (actually they are holomorphic)
$\tilde{f}_n,\tilde{g}_n:X\rightarrow X^\prime$ defined by
$\tilde{f}_n([z_0,z_1])=[z_0,n^{-1}z_1,z_1]$ and
$\tilde{g}_n([z_0,z_1])=[z_0,n^{-1}z_0,z_1]$. It is easily seen
that both sequences consist of regular maps, so they are good
maps. As $n\rightarrow \infty$, both sequences converge. Let
$\tilde{f}= \lim\tilde{f}_n$ and $\tilde{g}=\lim\tilde{g}_n$. Then
both $\tilde{f}$ and $\tilde{g}$ are good maps (as we shall see),
and $\tilde{f}=\tilde{g}$ as $C^\infty$ maps, but with different
isomorphism classes of compatible systems. In fact, the group
homomorphism of $\tilde{f}$ which is from $\Z_2$ to $\Z_2\oplus
\Z_2$ is  given by $x\rightarrow (x,1)$ at $[0,1]$, and
$x\rightarrow (x,x)$ at $[1,0]$. For $\tilde{g}$, it is given by
$x\rightarrow (x,x)$ at $[0,1]$ and $x\rightarrow (1,x)$ at
$[1,0]$.}\vskip 0.1in \hfill $\Box$

The operation of composition is well-defined for good maps, as shown in the
following lemma.

    \vskip 0.1in
\noindent{\bf Lemma 4.4.14: }{\it  Let $\tilde{f}$, $\tilde{g}$ be
two good $C^\infty$ maps, then the composition
$\tilde{g}\circ\tilde{f}$ is also a good $C^\infty$ map, and any
isomorphism class of compatible systems of $\tilde{f}$ and
$\tilde{g}$ determines a unique isomorphism class of compatible
systems for the composition $\tilde{g}\circ\tilde{f}$.}\vskip
0.1in

\noindent{\bf Proof:}
We need to show that the composition of two good maps is again a good map.
The rest of the lemma follows easily from Definition 4.4.4.

Without loss of generality, we assume that the compatible systems
of $\tilde{f}$ and $\tilde{g}$ are geodesic, given by
$\{\tilde{f}_{UU^\prime},\tau\}$ and
$\{\tilde{g}_{U_1U_1^\prime},\lambda\}$, where $\U=\{U_{p,i}:
p\in\Lambda,i\in\Z_+\}$ and $\U_1=\{U_{p,i}:
p\in\Lambda^\prime,i=\Z_+\}$. In the construction of
$\{\tilde{g}_{U_1U_1^\prime},\lambda\}$ we can add $f(\Lambda)$
into $\Lambda^\prime$ trivially so we can assume that
$f(\Lambda)\subset\Lambda ^\prime$. Moreover, we can assume that
each $U^\prime_{p,i}\in\U^\prime$ is contained in a
$U_{f(p),j}\in\U_1$. However, there is one problem: an inclusion
between elements in $\U^\prime$ may not imply a corresponding
inclusion between elements in $\U_1$ and therefore an inclusion
between the corresponding elements in $\U_1^\prime$. We can
overcome this problem by changing the elements in $\U_1^\prime$ to
suitable finite intersections which are again geodesically convex
and star-shaped. Let $\tilde{h}=\tilde{g}\circ\tilde{f}$. We
define a compatible system $\{\tilde{h}_{UU_1^\prime}, \theta\}$
of $\tilde{h}$ by taking the compositions of
$\{\tilde{f}_{UU^\prime},\tau\}$ and
$\{\tilde{g}_{U_1U_1^\prime},\lambda\}$, in which the injection
for each pair $U^\prime_{p,i}\in\U^\prime \subset
U_{f(p),j}\in\U_1$ is chosen to be the geodesic one, i.e., the one
which induces the identity map on the tangent space of $f(p)$.
\hfill $\Box$

\vspace{2mm}

Now consider a good $C^\infty$ map $\tilde{f}:X\rightarrow
X^\prime$ with an isomorphism class of compatible systems $\xi$.
Then we have an isomorphism class of pull-back orbibundles
$(TX^\prime)^\ast_\xi$ over $X$ and the good $C^\infty$ map
$\bar{f}_\xi:(TX^\prime)^\ast_\xi\rightarrow TX^\prime$. For any
$C^\infty$ section $\tilde{s}$ of $(TX^\prime)^\ast_\xi$, we take
the composition
$\tilde{f}_{\xi,s}=Exp\circ\bar{f}_\xi\circ\tilde{s}$ from $X$
into $X^\prime$. Then $\tilde{f}_{\xi,s}$ is a good $C^\infty$ map
with an isomorphism class of compatible systems determined by
$\xi$. A natural question is: given a good map $\tilde{g}$ nearby
$\tilde{f}$ with an isomorphism class of compatible systems, is
there an isomorphism class of compatible systems $\xi$ of
$\tilde{f}$, and a $C^\infty$ section $\tilde{s}$ of the pull-back
orbifold bundle $(TX^\prime)^\ast_{\xi}$ such that $\tilde{g}$ is
realized as $\tilde{f}_{\xi,s}$ ?  If there is, will $\xi$ and
$\tilde{s}$ be unique? These questions seem to be non-trivial in
general, but as we shall see, they can be dealt with in certain
special cases, e.g., when $\tilde{f}$, $\tilde{g}$ are
pseudo-holomorphic maps from a complex orbicurve into an almost
complex orbifold.

Suppose $\tilde{\gamma}: I\rightarrow X$ is a $C^\infty$ path in
$X$. Then it is easy to see that $\tilde{\gamma}$ is a good map so
that for each compatible system of $\tilde{\gamma}$, and any
orbifold bundle $E$ over $X$, there is a pull-back orbifold bundle
$E^\ast$ over $I$, which is a smooth vector bundle on $I$. One can
verify that when $\tilde{\gamma}$ is a parametrized geodesic,
there is a unique isomorphism class of compatible systems of
$\tilde{\gamma}$, using the fact that geodesics satisfy the unique
continuity property. Therefore all of the pull-back bundles are
isomorphic. We consider especially the case when $E=TX$. Let
$\nabla$ be a Riemannian connection on $E$, and let $\nabla^\ast$
be the pull-back connection on $E^\ast$. The parallel transport of
$E^\ast$ along $I$, $Par_{ab}:E^\ast_a\rightarrow E^\ast_b$
defined using the connection $\nabla^\ast$, defines through the
exponential map a diffeomorphism between $V_{\gamma(a)}$ and
$V_{\gamma(b)}$, which is equivariant under any automorphism of
$V_{\gamma(a)}$ (and $V_{\gamma(b)}$) which fixes the tangent
vector of $\tilde{\gamma}$.

    \vskip 0.1in
\noindent{\bf Lemma 4.4.15: }{\it  Let $(TX^\prime)^\ast_{\xi,s}$
be the pull-back orbifold bundle defined by $\tilde{f}_{\xi,s}$
with the canonical isomorphism class of compatible systems. Then
the map $Par_t:(TX^\prime)^\ast_{\xi}\rightarrow
(TX^\prime)^\ast_{\xi,ts}$ defined by parallel transportation
along parametrized geodesics
$\tilde{\gamma}_x(t)=Exp\circ\bar{f}_{\xi}\circ (x,t\tilde{s}(x))$
in $X^\prime$ is an isomorphism for each $t\in [0,1]$.}\vskip
0.1in

\noindent{\bf Proof:} This follows from the fact that parallel
transport along geodesics is solely determined by the Riemannian
connection on $TX^\prime$ which is compatible with the transition
maps of $(TX^\prime)^\ast_{\xi,ts}$ for each $t\in [0,1]$. \hfill
$\Box$

    \vskip 0.1in
\noindent{\bf Definition 4.4.16: }{\it  A sequence of pairs
$(\tilde{f}_n,\xi_n)$ is said to {\it converge to
$(\tilde{f}_0,\xi_0)$ in the $C^\infty$ topology} if there is a
sequence of compatible systems
$\{\tilde{f}_{n,UU_n^\prime},\lambda_n\}$ of $\tilde{f}_n$
representing $\xi_n$ and a compatible system
$\{\tilde{f}_{0,UU_0^\prime},\lambda_0\}$ of $\tilde{f}_0$
representing $\xi_0$ with the following property: for each $U\in
\U$, there is an integer $n(U)>0$ such that for each $n\geq n(U)$,
there is an injection $\delta_{U^\prime_n}$ from the uniformizing
system of $U^\prime_n$ into that of $U^\prime_0$, such that
$\delta_{U^\prime_n}\circ\tilde{f}_{n,UU_n^\prime}$ converges to
$\tilde{f}_{0,UU_0^\prime}$ in $C^\infty$, and
$\delta_{U^\prime_{1,n}}\circ
\lambda_n(i)=\lambda_0(i)\circ\delta_{U^\prime_{2,n}}$ holds for
$n\geq \max(n(U_1),n(U_2))$ for any injection $i$ of inclusion
$U_2\rightarrow U_1$.}\vskip 0.1in

    Recall that in the smooth case, the pull-back bundle can be
    constructed by fiber product. Namely, if $p: E\rightarrow Y$ is a
    bundle and $f: X\rightarrow Y$ is a map, $f^*E=\{(x, v);
    f(x)=p(v)\}$. The latter is often denoted by $X\times_Y E$.
    In the orbifold case, $f^*E$ is not the ordinary fiber product.
    However, one can define an orbifold fiber product as follows.

    Suppose that $f_1:X\rightarrow Z, f_2: X_2\rightarrow Z$ are
    good maps. Let $(\tilde{f}^1_{W_1W},\lambda^1_{W_1 W}),
    (\tilde{f}^2_{W_2 W}, \lambda^2_{W_2 W})$ be
    compatible systems for the same cover $\{W\}$ of $Z$. Then,
    the orbifold fiber product (still denoted by $X_1\times_Z X_2$)
    is constructed by gluing the $(W_1\times_W W_2, G_1\times_G
    G_2)$ together, where $G_i$ (resp. $G$) are the local group of $X_i$
    (resp. $Z$). If $W_1\times_Z W_2$ is smooth, $X_1\times_Z X_2$ is an orbifold.
    It is obvious that $f^*E=X\times_Y E$. One can
    check that the projections $W_1\times_W W_2\rightarrow W_1,
    W_2$ define good maps $p_i: X_1\times_Z X_2\rightarrow X_1, X_2$.

\subsection{A canonical stratification of orbifolds}
In this subsection, we give a brief discussion on the structure of
the singular set $\Sigma X$ and describe a canonical
stratification of orbifolds. See [K1] for more details.

The singular set $\Sigma X$ of an orbifold $X$ is not an orbifold
in general. But we can consider $\Sigma X$ as an immersed image of
a disjoint union of orbifolds. More precisely, let $(1)=(H_p^0),
(H_p^1),..., (H_p^{n_p})$ be all the orbit types of a geodesic
chart $(V_p,G_p,\pi_p)$ at $p$. For $q\in U_p=\pi_p(V_p)$, we may
take $U_q$ small enough so that $U_q\subset U_p$. Then any
injection $\phi: V_q\rightarrow V_p$ induces a unique homomorphism
$\lambda_\phi:G_q\rightarrow G_p$, which gives a correspondence
$(H_q^i)\rightarrow \lambda_\phi(H_q^i)=(H_p^j)$. This
correspondence is independent of the choice of $\phi$. Consider
the set of pairs: $$ \widetilde{\Sigma U}_p=\{(q,(H_q^i))|q\in
\Sigma U_p, i\neq 0\}. $$ Take one representative $H_q^i\in
(H_q^i)$. Then the pair $(q,(H_q^i))$ determines exactly one orbit
$[\tilde{q}]$ in the fixed point set $V_p^{H_p^j}$ by the action
of the normalizer $N_{G_p}(H_p^j)$, where $\tilde{q}=\phi(q)$,
$(H_p^j)=\lambda_\phi(H_q^i)$. The correspondence
$(q,(H_q^i))\rightarrow [\tilde{q}]$ gives a homeomorphism $$
\widetilde{\Sigma U}_p\cong
\coprod_{j=1}^{n_p}V_p^{H_p^j}/N_{G_p}(H_p^j), \hspace{2mm} \mbox
{(disjoint union)}, $$ which gives $\widetilde{\Sigma
X}=\{(p,(H_p^j))|p\in\Sigma X, j\neq 0\}$ an orbifold structure $$
\{\pi_{p,j}:(V_p^{H_p^j}, N_{G_p}(H_p^j))\rightarrow
V_p^{H_p^j}/N_{G_p}(H_p^j): p\in X, j=1,\cdots,n_p.\}. $$ The
canonical map $\pi:\widetilde{\Sigma X}\rightarrow \Sigma X$
defined by $(p,(H_p^j))\rightarrow p$ is surjective, which has a
$C^\infty$ lifting $\tilde{\pi}$ given locally by embeddings
$V_p^{H_p^j}\rightarrow V_p$. Therefore $\widetilde{\Sigma X}$ can
be regarded as a resolution of the singular set $\Sigma X$, called
the canonical resolution.

\vspace{2mm}

A point $(p,(H_p^j))$ in $\widetilde{\Sigma X}$ is called
generic if $G_p=H_p^j$. The set $\widetilde{\Sigma X}_{gen}$ of all generic
points is open dense in $\widetilde{\Sigma X}$, and the map
$\pi|\widetilde{\Sigma X}_{gen}:\widetilde{\Sigma X}_{gen}
\rightarrow \Sigma X$ is bijective. Hence we have a partition of $X$
into a disjoint union of smooth manifolds:
$$
X=X_{reg}\cup\widetilde{\Sigma X}_{gen},
$$
which is called the canonical stratification of $X$.


\begin{thebibliography}{}

\bibitem [AR]{AR} A. Adem and Y. Ruan, {\em Twisted orbifold K-theory}, preprint
\bibitem [AV]{AV} D. Abramovich and A. Vistoli, {\em Compactifying
 the space of stable maps}, math.AG/9908167
\bibitem [BD]{BD} V.V. Batyrev and  D. Dais, {\em Strong McKay correspondence,
string-theoretic Hodge numbers and mirror symmetry}, Topology, {\bf 35} (1996),
901-929.


\bibitem [CR1]{CR1} W. Chen and Y. Ruan, {\em A new cohomology theory for
orbifold}, preprint, April 2000.

\bibitem [CR2]{CR2} W. Chen and Y. Ruan, {\em Orbifold Quantum Cohomology},
research announcement, May 2000.

\bibitem [DHVW]{DHVW} L. Dixon, J. Harvey, C. Vafa and E. Witten, {\em
Strings on orbifolds, I, II}, Nucl.Phys. B261(1985), 678, B274(1986), 285.

\bibitem [FO]{FO} K. Fukaya and K. Ono, {\em Arnold conjecture and
Gromov-Witten invariant}, Topology {\bf 38}(1999), no.5, 933-1048.

\bibitem [HW]{HW} P. Hartman and A. Wintner, {\em On the local behavior
of solutions of nonparabolic partial differential equations}, Amer. J. Math.
{\bf 75}(1953), 449-476.

\bibitem [HH]{HH} F. Hirzebruch and T. H\"{o}fer, {\em On the Euler number
of an orbifold}, Math. Ann. {\bf 286} (1990) 255-.


\bibitem [K1]{K1} T. Kawasaki, {\em The signature theorem for V-manifolds},
Topology, {\bf 17} (1978), 75-83.

\bibitem [Ke]{Ke} J. Kelley, {\em General topology}.


\bibitem [LT]{LT} J. Li and G. Tian, {\em Virtual moduli cycles and
Gromov-Witten invariants of general symplectic manifolds}, Topics in
symplectic 4-manifolds (Irvine, CA, 1996), 47-83, First Int. Press Lect.Ser.,I,
Internat. Press, Cambridge, MA, 1998.

\bibitem [McK]{Mc} J. McKay, {Graphs, singularities and finite groups},
Proc. Symp. in Pure Math., {\bf 37}(1980), 183-186.


\bibitem [Re]{Re} M. Reid, {\em McKay correspondence}, AG/9702016 v3.

\bibitem [Ru]{Ru} Y. Ruan, {\em Symplectic topology on algebraic
3-folds}, Jour. Diff. Geom., 39, 215-227(1994)
\bibitem [Ru1]{Ru1} Y. Ruan, {\em Virtual neighborhoods and pseudo-holomorphic
curves}, Turkish J. Math. {\bf 23}(1999), no. 1, 161-231.
\bibitem [Ru2]{Ru2}Y. Ruan, {\em Discrete torsion and twisted orbifold cohomology},
    math.AG/0005299
\bibitem [S]{S}         I. Satake,
                        {\em The Gauss-Bonnet theorem for V-manifolds},
                        J. Math. Soc. Japan {\bf 9} (1957), 464-492.


\bibitem [W]{W} E. Witten, {\em Two dimensional gravity and intersection theory
    on moduli space, } Survey in Diff. Geom, 1, 243-310 (1991).

\bibitem [Z]{Z} E. Zaslow, {\em Topological orbifold models and
quantum cohomology rings}, Comm. Math. Phys. {\bf 156} (1993), no.2, 301-331.

\end{thebibliography}
\end{document}